# Adaptive Inferential Method for Monotone Graph Invariants

Junwei Lu*     Matey Neykov*     Han Liu*†

**Abstract**

We consider the problem of undirected graphical model inference. In many applications, instead of perfectly recovering the unknown graph structure, a more realistic goal is to infer some graph invariants (e.g., the maximum degree, the number of connected subgraphs, the number of isolated nodes). In this paper, we propose a new inferential framework for testing nested multiple hypotheses and constructing confidence intervals of the unknown graph invariants under undirected graphical models. Compared to perfect graph recovery, our methods require significantly weaker conditions. This paper makes two major contributions: (i) Methodologically, for testing nested multiple hypotheses, we propose a skip-down algorithm on the whole family of monotone graph invariants (The invariants which are non-decreasing under addition of edges). We further show that the same skip-down algorithm also provides valid confidence intervals for the targeted graph invariants. (ii) Theoretically, we prove that the length of the obtained confidence intervals are optimal and adaptive to the unknown signal strength. We also prove generic lower bounds for the confidence interval length for various invariants. Numerical results on both synthetic simulations and a brain imaging dataset are provided to illustrate the usefulness of the proposed method.

## 1 Introduction

Graphical models are widely used for modeling complex networks such as gene regulatory networks and brain connectivity networks (Luscombe et al., 2004; Rubinov and Sporns, 2010). We use an undirected graph $G = (V, E)$ to represent the conditional dependency structure of a $d$-dimensional random vector $\boldsymbol{X} = (X_1, \ldots, X_d)^T$, where each vertex in $V = \{1, \ldots, d\}$ corresponds to a component of the random vector $\boldsymbol{X}$. Specifically, two nodes $j$ and $k$ are connected if and only if $X_j$ and $X_k$ are conditionally dependent given the other variables.

*Department of Operations Research and Financial Engineering, Princeton University, Princeton NJ 08544, USA; e-mail: {junweil,mneykov,hanliu}@princeton.edu

†Research supported by NSF DMS1454377-CAREER; NSF IIS 1546482-BIGDATA; NIH R01MH102339; NSF IIS1408910; NIH R01GM083084.

Learning the structure of the graph in a graphical model has been widely studied in the literature. Many theoretical studies focus on achieving a perfect recovery of the true graph, i.e., they aim to find a consistent graph estimator $\widehat{G}$ satisfying $\mathbb{P}(\widehat{G} = G) \to 1$ as the sample size goes to infinity. For the Gaussian graphical model, many works estimate the graph through the precision matrix $\boldsymbol{\Theta}$ satisfying $\boldsymbol{\Theta}_{jk} \neq 0$ if and only if the edge $(j, k) \in E$ (Meinshausen and Bühlmann, 2006; Yuan and Lin, 2007; Friedman et al., 2008; Peng et al., 2009; Lam and Fan, 2009; Ravikumar et al., 2011; Cai et al., 2011; Shen et al., 2012). For the Ising model, the true graph is determined by the sparsity pattern of the interaction parameters and we can estimate the parameters by sparse logistic regression (Ravikumar et al., 2010). A number of authors relax the Gaussian or Ising assumption by allowing the node-conditional distribution to belong to a univariate exponential family (see, e.g., Lee et al., 2006; Höfling and Tibshirani, 2009; Jalali et al., 2011; Anandkumar et al., 2012; Yang et al., 2012; Allen and Liu, 2012; Yang et al., 2013). This line of research extents to the mixed graphical model framework in which the conditional distributions of two nodes can belong to two different distributions from a univariate exponential family (see, e.g., Fellinghauer et al., 2013; Lee and Hastie, 2015; Yang et al., 2014a; Chen et al., 2015). Semiparametric extensions using copulas are also developed in Liu and Wasserman (2009), Liu et al. (2012a), Xue and Zou (2012), and Liu et al. (2012b), and extended for mixed data in Fan et al. (2014). Danaher et al. (2014), Qiu et al. (2013), Mohan et al. (2014) consider joint estimation of multiple graphical models. However, in order to achieve the perfect graph recovery in these works, some strong conditions especially the minimal signal strength assumption are usually needed for the CLIME (Cai et al., 2011), the neighborhood selection (Meinshausen and Bühlmann, 2006; Zhou et al., 2009), the graphical Lasso (Lam and Fan, 2009) and the transelliptical graphical model (Liu et al., 2012b). The minimal signal strength assumption imposes that every edge has strong enough signal which is unlikely to hold in reality.

In this paper, we consider a statistical problem which imposes weaker assumptions compared to minimal signal strength type of conditions. In particular, we aim to infer the topological structure of the graph without imposing strong conditions required for perfect recovery. The topological structure of a graph is an important feature in many real applications. For instance, in a gene regulatory network, it is scientifically interesting to infer the number of metabolic cycles where cells move through in the genomic network (Luscombe et al., 2004). In neuroscience applications, one critical problem is to infer the degrees of various cerebral areas in the brain network during certain cognitive processes (Hagmann et al., 2008). In this paper, we propose a unified inferential framework on the topological quantities related to graphical models. In particular, we focus on topological quantities which are invariant under graph isomorphism. We call such quantities graph invariants, denoted as $\mathcal{I}$. Examples of graph invariants we consider include: maximum degree, maximum clique size, chromatic number, girth, the number of connected subgraphs, etc. We are interested in testing hypotheses with a nested structure:

$$H_{0k} : \mathcal{I} \leq k, \text{ for } k \text{ being all possible values of } \mathcal{I}. \tag{1.1}$$

Our goal is to propose a test with family-wise error control and power optimality. Moreover, we aim to derive a confidence interval for $\mathcal{I}$ utilizing the method for nested multiple hypotheses. Our paper will focus on the Gaussian graphical models whose graph structure is encoded by the sparsity of precision matrices, but we will also discuss how the proposed method applies to other more flexible graphical models.

A related research area on graphical models are concerned with inferential methods on a single edge in the graphical model (Yang et al., 2014b; Gu et al., 2015; Ren et al., 2015; Janková and van de Geer, 2016)). A few authors have also proposed to test multiple edges simultaneously (among others, Liu, 2013; Cai et al., 2014; Lu et al., 2015). However, their proposals require a user-specified set of edges to be tested simultaneously, and thus cannot be directly applied to infer the topological structure of a graph. In another related research direction, several authors (Arias-Castro et al., 2012, 2015, 2011; Addario-Berry et al., 2010) have considered testing whether a graph is empty or it has certain topological structure and discussed lower bounds. One of the drawbacks of these approaches is that the null hypothesis is constrained to be an empty graph. The most related work to our paper is Neykov et al. (2016). They propose a general theoretical framework for analyzing lower bounds for testing hypothesis of whether a graph satisfies a certain graph property. To match the lower bound, a data splitting procedure is required. However, our method does not need to split the dataset.

Our contributions to the inference of graph invariants for graphical models are three-fold in methodology, theory of inference and fundamental limits.

- **Methodology.** We propose a fully data-driven method for the nested hypotheses on a wide family of graph invariants $\mathcal{I}$, which also derives confidence intervals for invariants. Unlike Neykov et al. (2016)'s data splitting method, which tests a single hypothesis $H_0: \mathcal{I} \leq k$ with a specified $k$, our method does not need to pre-determine the value of $k$. In practice, it is usually hard to decide which $k$ to test. On the other hand, the nested hypotheses in (1.1) initiate a fully data-driven approach to choosing the value of $k$. Nested multiple hypotheses problem is considered by Bauer and Hackl (1987) for continuous parameters test. We propose a skip-down algorithm for the nested hypotheses on combinatorial quantities, which iteratively screens critical edge sets for the graph invariant of interest. We also derive confidence intervals from the skip-down algorithm. Constructing confidence intervals for combinatorial quantities is challenging as the standard confidence interval theory on continuous and smooth parameters does not directly apply. We show that our proposed confidence interval $[\widehat{\mathcal{I}}_L, \widehat{\mathcal{I}}_U]$ is asymptotically honest for all monotone graph invariants.

- **Theory of Inference.** For the Gaussian graphical model, we demonstrate that the length of our proposed confidence interval $[\widehat{\mathcal{I}}_L, \widehat{\mathcal{I}}_U]$ is adaptive to the signal strength of the precision matrix $\boldsymbol{\Theta}$. To elaborate in detail, we define the significant edge set as $E_{\mathrm{Sig}}(\boldsymbol{\Theta}) := \{(j,k) \,|\, |\boldsymbol{\Theta}_{jk}| \geq C\sqrt{\log d/n}\}$, which only keeps edges with signal strength larger than $C\sqrt{\log d/n}\}$ for some sufficiently large constant $C$. We show that the expected length of confidence interval $[\widehat{\mathcal{I}}_L, \widehat{\mathcal{I}}_U]$ is adaptive to the value of the invariant $\mathcal{I}(E_{\mathrm{Sig}}(\boldsymbol{\Theta}))$, which denotes the value of the graph invariant based on edge set $E_{\mathrm{Sig}}(\boldsymbol{\Theta})$. Specifically, under some regularity and scaling conditions, we show that for any monotone invariant $\mathcal{I}$, $\mathbb{E}_{\boldsymbol{\Theta}}[\widehat{\mathcal{I}}_U - \widehat{\mathcal{I}}_L]$ is decreasing when $\mathcal{I}(E_{\mathrm{Sig}}(\boldsymbol{\Theta})))$ becomes larger. More details regarding the adaptivity will be shown in Theorem 4.5.

- **Fundamental Limits.** We provide a general theoretical result for the lower bound of the confidence interval length for graph invariants. We provide a sufficient condition under which the confidence interval length is optimal for a large family of invariants. Our sufficient condition is uniquely characterized by the geometry of the graph invariant. This makes it easier to show the lower bound of confidence interval length for many graph invariants. In comparison, the lower bound results in Neykov et al. (2016) for graph properties hypothesis test rely on probabilistic type conditions such as negative association (Joag-Dev and Proschan, 1983), which are difficult to verify in practice. Our result also bares similarity to results in Cai and Guo (2015), who show a lower bound on the confidence interval length for the high dimensional linear model. However, our confidence interval is for combinatorial quantities where classical parametric and continuous statistical theory cannot be directly applied.

## 1.1 Notation

We denote the graph induced by a matrix $\boldsymbol{\Theta}$ as $G(\boldsymbol{\Theta}) = (V, E(\boldsymbol{\Theta}))$, where $(j,k) \in E(\boldsymbol{\Theta})$ if and only if $\boldsymbol{\Theta}_{jk} \neq 0$. Given any edge set $E$, we let $V(E) = \{v \,|\, \exists\, u \in V, \text{ s.t. } (u,v) \in E\}$ be the set of vertices in $E$. We denote the cardinality of a set $V$ as $|V|$. Given two sequences $\{a_n\}_{n=1}^\infty$ and $\{b_n\}_{n=1}^\infty$, if $a_n \leq Cb_n$ for some finite positive constant $C$ for all $n$ large enough, we denote it as $a_n = O(b_n)$ or $b_n = \Omega(a_n)$. We use the notation $a_n \asymp b_n$ if $a_n = O(b_n)$ and $b_n = O(a_n)$. We also write $a_n = o(b_n)$ if $a_n/b_n \to 0$ as $n$ goes infinity. We denote $a \vee b = \max(a,b)$ and $a \wedge b = \min(a,b)$ for any two scalars $a$ and $b$. Given a scalar $x$, $\lfloor x \rfloor$ denotes the largest integer which is smaller than $x$. Let $[n]$ denote the set $\{1, \ldots, n\}$. For a matrix $\mathbf{A} \in \mathbb{R}^{d_1 \times d_2}$, we denote the max-norm $\|\mathbf{A}\|_{\max} = \max_{i \in [d_1], j \in [d_2]} |\mathbf{A}_{ij}|$, the $\ell_1$-norm $\|\mathbf{A}\|_1 = \max_{j \in [d_2]} \sum_{i \in [d_1]} |\mathbf{A}_{ij}|$, and the operator norm $\|\mathbf{A}\|_2 = \sup_{\|\mathbf{v}\|_2 = 1} \|\mathbf{A}\mathbf{v}\|_2$. Throughout the paper, we let $C, C_1, C_2, \ldots$ be generic constants whose values may change in different places.

## 1.2 Paper Organization

The rest of our paper is organized as follows. In Section 2, we introduce preliminary concepts of graph theory including graph invariants and graph properties. In Section 3, we provide the skip-down algorithm for graph invariant confidence intervals and tests on graph properties. In Section 4, we show theoretical results on the optimality for our confidence intervals and tests. In Section 5, we generalize our method to more flexible graphical models. We provide numerical results on both synthetic simulations and brain imaging applications in Section 7.

## 2 Preliminaries on Graph Theory

Let $V = \{1, \ldots, d\}$ be the set of nodes and the edge set $E$ be a subset of $V \times V$. $G = (V, E)$ denotes an undirected graph with vertices in $V$ and edges in $E$. We say $G = (V, E)$ is isomorphic to $G' = (V', E')$ if there exists a one-to-one mapping $\sigma : V \to V'$, such that $(j, k) \in E$ if and only if $(\sigma(j), \sigma(k)) \in E'$. We also write $G \preceq G'$ if $V \subseteq V'$ and $E \subseteq E'$. Let $\mathcal{G}_d$ be the set of all graphs with $d$ vertices. A graph invariant is a function $\mathcal{I} : \mathcal{G}_d \to \mathbb{Z}$ such that $\mathcal{I}(G) = \mathcal{I}(G')$ if $G$ is isomorphic to $G'$. In other words, a graph invariant is a geometric characterization of the graph, which is invariant to vertex permutations. In this paper, we are interested in a special family of graph invariants called monotone invariants defined as follows.

**Definition 2.1.** A graph invariant $\mathcal{I}$ is *monotone* whenever $\mathcal{I}(G) \leq \mathcal{I}(G')$ for all $G \preceq G'$.

Specifically, if a graph invariant takes binary values, i.e., $\mathcal{P} : \mathcal{G}_d \to \{0, 1\}$ such that $\mathcal{P}(G) = \mathcal{P}(G')$ if $G$ is isomorphic to $G'$, we call this invariant $\mathcal{P}$ a graph property. We also say a graph $G$ satisfies a property $\mathcal{P}$ if $\mathcal{P}(G) = 1$. Similarly to Definition 2.1, we define the monotone graph property as follows.

**Definition 2.2.** A graph property $\mathcal{P}$ is *monotone* whenever $\mathcal{P}(G) \leq \mathcal{P}(G')$ for all $G \preceq G'$.

In this paper, we will always take the vertex set $V = \{1, \ldots, d\}$. Therefore, for simplicity of notation, given a graph $G = (V, E)$, we also write an invariant or a property as a function of edge set, i.e., $\mathcal{I}(E) = \mathcal{I}(G)$ and $\mathcal{P}(E) = \mathcal{P}(G)$. From the definitions above, we can see that a monotone graph invariant is non-decreasing under addition of edges. Similarly, if a graph $G$ satisfies a monotone property, the property is preserved under addition of edges to $G$.

Many extensively used graph invariants are monotone. For example, we can easily check that the following invariants are monotone:

- The maximum degree of $G$, which is the largest number of vertices connected to a single vertex in the graph (see Figure 1(a) for an example of a graph with maximum degree 5);

- The size of the longest chain in $G$, where a chain is a set of edges connecting a sequence of distinct vertices consecutively, and the size of a chain is the number of edges it contains (see Figure 1(b) for an example of a graph with longest chain of size 5);

- The size of the largest clique in $G$, where a clique is a subgraph such that any pair of its vertices are connected and the size of a clique is the number of vertices it contains (see Figure 1(c) for an example of a graph with largest clique size 5);

- The chromatic number of $G$, which is the smallest number of colors needed to color all vertices, so that no vertices sharing the same color are adjacent (see Figure 1 (d) for an example of a graph with chromatic number 3);

- The negative number of isolated nodes of $G$, where a isolated node is a vertex with no neighbor. Here we take the negative number because the number of isolated nodes is non-decreasing under the addition of edges (see Figure 1 (e) for an example of a graph with 5 isolated nodes);

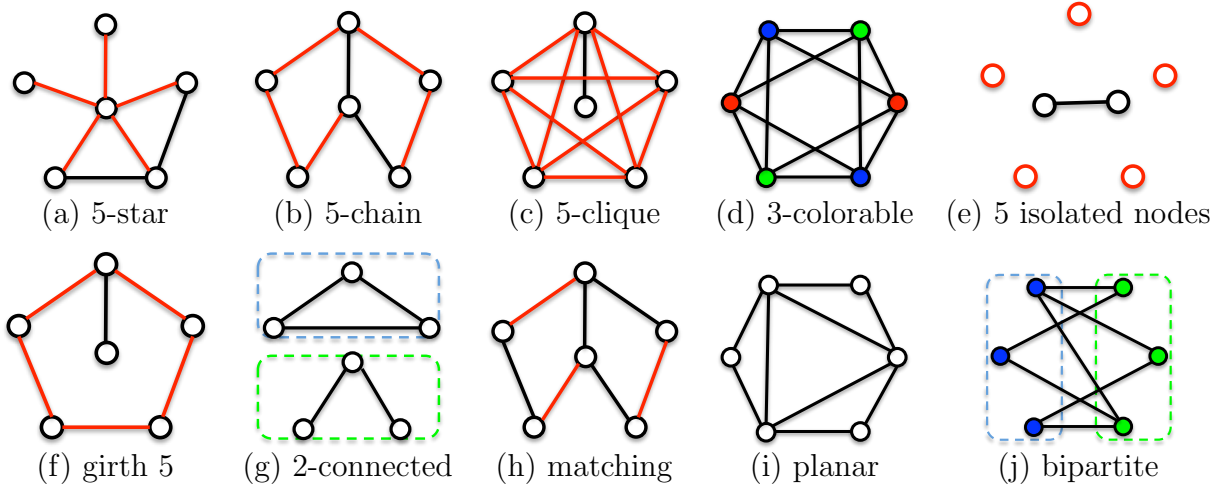

Figure 1: Examples of graph invariants and properties.

- The negative girth of $G$, where the girth of a graph is the length of the shortest cycle and the girth equals infinity when $G$ has no cycle. We take the negative girth also because the girth itself is non-increasing when adding edges (see Figure 1 (f) for an example of a graph with girth 5);

- The negative number of connected subgraphs of $G$, where a subgraph is connected if any pair of its vertices are connected by a chain (see Figure 1 (g) for an example of a graph with connected subgraphs). Notice that the number of connected subgraphs is again non-increasing under the addition of edges and thus the negative number of connected subgraphs is monotone.

In the last three examples of monotone invariants listed above, we consider the negative values as more natural quantities. The negative sign is introduced since the invariants – number of isolated nodes, girth and number of connected subgraphs are non-increasing under the addition of edges, i.e., they are "monotone decreasing" instead of monotone increasing. In order to unify our methodology, we focus only on monotone increasing invariants in this paper. Any monotone decreasing invariant, can be simply converted to a monotone increasing one by taking its negative value.

Any graph invariant $\mathcal{I}$, induces a graph property $\mathcal{P}_{\mathcal{I},k}$ defined as $\mathcal{P}_{\mathcal{I},k}(G) = 0$ if $\mathcal{I}(G) \leq k$ and $\mathcal{P}_{\mathcal{I},k}(G) = 1$ if $\mathcal{I}(G) > k$. It is easy to see that $\mathcal{P}_{\mathcal{I},k}$ is a monotone graph property if $\mathcal{I}$ is monotone. For example, when $\mathcal{I}$ is the chromatic number, the induced property $\mathcal{P}_{\mathcal{I},k}(G) = 0$ if and only if $G$ is $k$-colorable, which means we can assign each vertex a color from $k$ colors such that no two adjacent vertices share the same color, and in particular, $\mathcal{P}_{\mathcal{I},2}(G) = 0$ if and only if $G$ is bipartite[1]. See Figure 1(e) for an example of a 3-colorable graph and Figure 1(j) for an example of a bipartite graph. Another example is when $\mathcal{I}$ is the negative girth.

---

[1] Recall that a bipartite graph is a graph whose vertices can be separated into two disjoint sets so that every edge connects vertices from one set to another.

The induced property $\mathcal{P}_{\mathcal{I},-\infty}(G) = 0$ if and only if $G$ is a forest and the induced property $\mathcal{P}_{\mathcal{I},-4}(G) = 0$ if and only if $G$ is triangle-free.

Therefore, the above examples of monotone invariants naturally imply corresponding examples of monotone properties. In addition to these natural examples, we also have the following examples of monotone properties:

- $G$ has a perfect matching, i.e. $G$ has a subset of edges without common vertices such that each vertex of $G$ is an endpoint of one of these edges (see Figure 1(h) for an example of a graph having perfect matching);

- $G$ is not planar, where a graph is planar if it can be drawn on the plane in such a way so that its edges only intersect at the vertices of the graph (see Figure 1(i) for an example of a planar graph);

- $G$ has a subgraph which is isomorphic to a given graph $H$.

The last property can actually derive a family of monotone graph properties given different subgraph pattern $H$. If the given graph $H$ is a $k$-star, the last property becomes the property that $G$ contains a star of size equal to or larger than $k$, which is equivalent to the induced property $\mathcal{P}_{\mathcal{I},k}$ when $\mathcal{I}$ is the maximum degree. Similarly, we can also set $H$ to be a chain of size $k$ or a clique of size $k$ which corresponds to the induced property $\mathcal{P}_{\mathcal{I},k}$ for $\mathcal{I}$ to be the size of the longest chain or the size of the largest clique. Another example is when $H$ is a graph with $d/2$ disjoint edges, which is equivalent to the existence of perfect matching. We visualize the aforementioned examples of $H$'s with the red edges in Figures 1(a) - 1(d) and 1(h).

## 3 Inferential Methods for Graph Invariants

In this section, we introduce the general framework for testing nested hypotheses for monotone graph invariants and then derive confidence intervals from the tests. Theoretical results on the validity of the tests and confidence intervals are also provided.

We begin with formulating our inferential problems on graph invariants. In order to illustrate the main idea of our testing procedure, we focus on the Gaussian graphical model first and discuss further extensions to other models in Section 5. Let $\mathbf{X} = (\mathbf{X}_1, \ldots, \mathbf{X}_n)$ be i.i.d. observations of the random vector $\boldsymbol{X} \sim N(\boldsymbol{0}, \boldsymbol{\Sigma})$. The precision matrix $\boldsymbol{\Theta} = \boldsymbol{\Sigma}^{-1}$ encodes the underlying conditional independence graph $G = (V, E)$. In particular, the edge $(j, k) \in E$ if and only if $\boldsymbol{\Theta}_{jk} \neq 0$. We consider the parameter space of precision matrices

$$\mathcal{U}_s = \left\{ \boldsymbol{\Theta} \in \mathbb{R}^{d \times d} \,\Big|\, \lambda_{\min}(\boldsymbol{\Theta}) \geq 1/\rho, \max_{j \in [d]} \|\boldsymbol{\Theta}_j\|_0 \leq s, \|\boldsymbol{\Theta}\|_1 \leq M \right\}. \tag{3.1}$$

The above class requires the precision matrix to have at most $s$ nonzero entries for each column and therefore the graphs induced by such matrices have maximum degrees at most $s$.

We denote the range of a graph invariant $\mathcal{I}$ as $[I_L^*, I_U^*]$. Throughout the paper, $[I_L^*, I_U^*]$ usually represents the default range of a graph invariant. For example, for any graph property $\mathcal{P}$, the default range is $[0, 1]$. If $\mathcal{I}$ is the negative number of isolated nodes, the default range is $[-d, 0]$ and for $\mathcal{I}$ being the negative number of connected subgraphs, the range is $[-d, -1]$. Sometimes we may have prior knowledge on the precision matrix which can make the default range smaller. For example, if we know that $\Theta \in \mathcal{U}_s$, the range of maximum degree can be set to $[0, s]$. In this paper, we mainly focus on two inferential problems on graphical models. Let $G$ be the graph which the graphical model is Markov with respect to. The first problem is testing multiple hypotheses with a nested structure

$$H_{0k} : \mathcal{I}(G) \leq k \text{ versus } H_{1k} : \mathcal{I}(G) > k, \text{ for } k \in [I_L^*, I_U^*). \tag{3.2}$$

If $\mathcal{I}(G) = k^*$, we have $H_{0k}$ is true for $k \geq k^*$ and $H_{1k}$ is true for $k < k^*$. For each $k \in [I_L^*, I_U^*)$, we aim to propose a test $\psi_k \in \{0, 1\}$ for $H_{0k}$, such that the family-wise type-I error is controlled, i.e.,

$$\lim_{n \to \infty} \mathbb{P}_{\mathcal{I}(G) = k^*}\big(\text{There exists a type-I error}\big) = \lim_{n \to \infty} \mathbb{P}_{\mathcal{I}(G) = k^*}\big(\exists k \geq k^* \text{ s.t. } \psi_k = 1\big) \leq \alpha.$$

When an invariant is a property, since it only takes binary values, we have a special single hypothesis on graph property

$$H_0 : \mathcal{P}(G) = 0 \text{ versus } H_1 : \mathcal{P}(G) = 1,$$

where $\mathcal{P}$ is a monotone property of interest. We aim to propose a test $\psi \in \{0, 1\}$ at significance level $\alpha$ such that

$$\lim_{n \to \infty} \mathbb{P}_{H_0}(\psi = 1) \leq \alpha.$$

The second problem is constructing a confidence interval of a monotone invariant $\mathcal{I}$. We aim to construct a confidence interval $[\widehat{\mathcal{I}}_L, \widehat{\mathcal{I}}_U]$ at significance level $\alpha$ such that

$$\lim_{n \to \infty} \mathbb{P}\big(\mathcal{I}(G) \in [\widehat{\mathcal{I}}_L, \widehat{\mathcal{I}}_U]\big) \geq 1 - \alpha.$$

## 3.1 A Generic Skip-Down Algorithmic Framework

In this subsection, we first describe the proposed skip-down algorithm (See Algorithm 1) for testing the hypothesis in (3.2). We then show that the same algorithm can also deliver a valid confidence interval for the graph invariant being tested.

To describe the algorithm, we first introduce a concept called "the critical edge set". Intuitively, for any graph invariant $\mathcal{I}$ to be tested, we need to find a set of critical edges which may potentially change the value of $\mathcal{I}$. A formal definition is as follows.

**Definition 3.1** (Critical edge set). For any monotone graph invariant $\mathcal{I}$ and an edge set $E_0 \subseteq V \times V$, we define the critical edge set of $\mathcal{I}$ with respect to $E_0$ as

$$\mathcal{C}_\mathcal{I}(E_0) := \big\{e \in E_0^c \,|\, \exists E' \supseteq E_0 \text{ s.t. } \mathcal{I}(E') > \mathcal{I}(E' \setminus \{e\})\big\}. \tag{3.3}$$

The definition of $\mathcal{C}_\mathcal{I}(E_0)$ is quite abstract. To understand its intuition, we first consider a special case when $\mathcal{I}$ is a monotone graph property (denoted by $\mathcal{P}$). As $\mathcal{P}$ only takes binary values, $\mathcal{C}_\mathcal{P}(E_0)$ can be equivalently written as

$$\mathcal{C}_\mathcal{P}(E_0) := \{e \in E_0^c \,|\, \exists E' \supseteq E_0 \text{ s.t. } \mathcal{P}(E') = 1 \text{ and } \mathcal{P}(E' \backslash \{e\}) = 0\}. \tag{3.4}$$

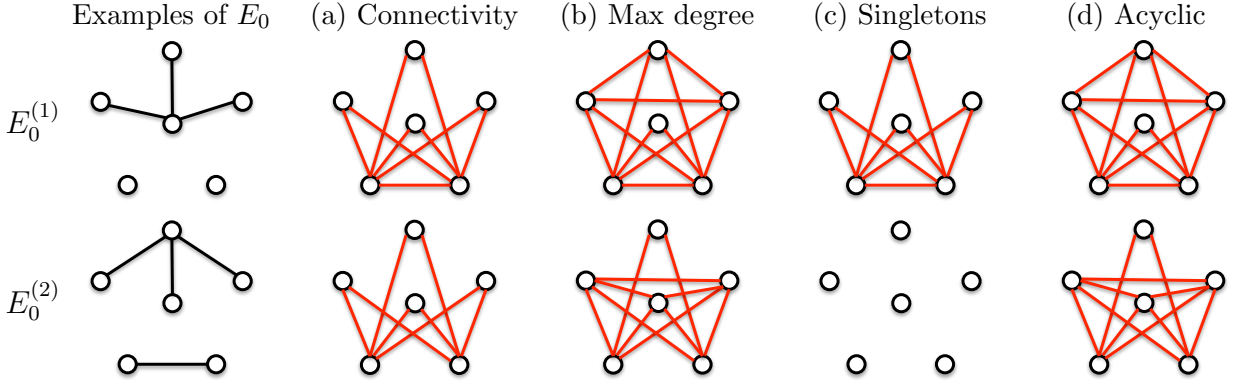

Figure 2: Two examples of the critical edge sets $\mathcal{C}_\mathcal{I}(E_0)$ (in red edges) for (a) $\mathcal{I}$ is the negative number of connected subgraphs, (b) $\mathcal{I}$ is the maximum degree, (c) $\mathcal{I}$ is the negative number of isolated nodes, and (d) $\mathcal{I} = 0$ when the graph is a forest and $\mathcal{I} = 1$ otherwise. The first row is for $E_0 = E_0^{(1)}$ and the second row is for the second example when $E_0 = E_0^{(2)}$.

An edge $e$ is critical for $\mathcal{P}$ with respect to $E_0$, if there exits a set of edges $\{e_1, \ldots, e_k\}$ that do not belong to $E_0 \cup \{e\}$ such that adding them to $E_0$ makes $(V, E_0 \cup \{e_1, \ldots, e_k, e\})$ satisfies $\mathcal{P}$ but $(V, \{e_1, \ldots, e_k\})$ does not. This implies that $e$ has the "potential" to change the property of $E_0$ through edge addition. In Figure 2, we provide examples of critical edge sets for various invariants. From Figure 2, we see that $\mathcal{C}_\mathcal{I}(E_0)$ can be strictly smaller than $E_0^c$ for some invariants. This implies that we can ignore some edges and leads to a more powerful test. For a more detailed discussion of critical edges, see Section 3.2.

We begin by introducing the high level idea of testing the nested hypotheses in (3.2). At first, we do not have any information on the true graph. Therefore, we test the existence of edges in $\mathcal{C}_\mathcal{I}(\varnothing)$. Suppose we have rejected edges in some edge set $E_0 \subseteq \mathcal{C}_\mathcal{I}(\varnothing)$, then we can update our knowledge that $E_0$ is a subgraph of the true graph. Since $\mathcal{I}$ is a monotone invariant, we also know that $\mathcal{I}(G) \geq \mathcal{I}(E_0)$. In the next step, we can further update our knowledge of the true graph by testing the existence of edges in $\mathcal{C}_\mathcal{I}(E_0)$. We can repeat this procedure for multiple times and iteratively refine our knowledge on $\mathcal{I}(G)$.

An essential component step of the above idea is to test the existence of a certain edge set $E$. This can be achieved by a generic precision matrix estimator $\widehat{\Theta}^d$ and a $(1 - \alpha)$ quantile estimator $c(\alpha, E)$ for the statistic $T_E = \max_{e \in E} \sqrt{n}(\widehat{\Theta}_e^d - \Theta_e)$ given any edge set $E \subseteq V \times V$ satisfying

$$\lim_{n \to \infty} \sup_{\Theta \in \mathcal{U}_s} \left| \mathbb{P}\left( \max_{e \in E} \sqrt{n}(\widehat{\Theta}_e^d - \Theta_e) > c(\alpha, E) \right) - \alpha \right| = 0. \tag{3.5}$$

We will provide concrete examples of $\widehat{\Theta}^d$ and $c(\alpha, E)$ satisfying (3.5) in Section 5.

**Algorithm 1** Skip-Down Method for Inferring Graph Invariant $\mathcal{I}$
---
**Input:** $\{\widehat{\Theta}_e^d\}_{e \in V \times V}$ and the range $[I_L^*, I_U^*]$.
Initialize $t = 0, E_0 = \varnothing$ and the interval $[\widehat{\mathcal{I}}_L, \widehat{\mathcal{I}}_U] = [I_L^*, I_U^*]$.
**repeat**
   $t \leftarrow t + 1$;
   Select the screening set: $\mathcal{A} \leftarrow \mathcal{C}_\mathcal{I}(E_{t-1})$;
   Update the rejected set: $E_t \leftarrow E_{t-1} \cup \{e \in \mathcal{A} \,|\, \sqrt{n} \cdot |\widehat{\Theta}_e^d| > c(\alpha, \mathcal{A})\}$;
   Update the interval $[\widehat{\mathcal{I}}_L, \widehat{\mathcal{I}}_U] \leftarrow [\mathcal{I}(E_t), I_U^*]$;
**until** $\mathcal{I}(E_t) \geq I_U^*$ or $E_t = E_{t-1}$
**Output:**
• **Nested hypotheses:** Let $\widehat{\mathcal{I}}_L = \min(\max(\mathcal{I}(E_t), I_L^*), I_U^*)$. At the significance level $\alpha$, we reject $H_{0k}$ for $k < \widehat{\mathcal{I}}_L$ and not reject $H_{0k}$ for $k \in [\widehat{\mathcal{I}}_L, I_U^*]$.
• **Confidence interval:** The $1 - \alpha$ confidence interval $[\widehat{\mathcal{I}}_L, I_U^*]$.
---

Algorithm 1 implements the high level idea above to test hypotheses in (3.2) and construct a confidence interval for the monotone invariant $\mathcal{I}$. We call this algorithm the skip-down method. The skip-down method is motivated by the step-down method (Romano and Wolf, 2005), which is designed for controlling the family-wise error of general multiple hypothesis tests. Our skip-down method explicitly exploits the nested structure of multiple hypotheses and has stronger power. Particularly, in each step, the step-down method tests all hypotheses which have not been previously rejected, while the skip-down method "skips" edges not belonging to the critical edge set, and hence improves the power of the test.

To construct confidence intervals, we choose the lower side as $\widehat{\mathcal{I}}_L$ in Algorithm 1 while the upper side is simply chosen to be $I_U^*$. This asymmetry is due to the monotonicity of graph invariants. In fact, we cannot construct a better upper bound for maximum degree according to the following theorem.

**Theorem 3.1.** Denote $\mathcal{I}_{\text{Deg}}(G)$ as the maximum degree of $G$. Define the family of confidence upper bound of $\mathcal{I}_{\text{Deg}}$ as

$$U(\mathcal{I}_{\text{Deg}}, \alpha) = \left\{\widehat{U}(\cdot) : \mathbb{R}^{d \times n} \to [0, s] \,\Big|\, \inf_{\Theta \in \mathcal{U}_s} \mathbb{P}_\Theta\big(\mathcal{I}(\Theta) \leq \widehat{U}(\mathbf{X})\big) \geq 1 - \alpha\right\}. \qquad (3.6)$$

Let $\mathcal{U}_s(\theta) = \big\{\Theta \in \mathcal{U}_s \,\big|\, \min_{e \in E(\Theta)} |\Theta_e| \geq \theta\big\}$. If $s = o(d^{1/2})$ and $\theta \leq C\sqrt{\log d/n}$ for some sufficiently small positive constant $C$, we have

$$\liminf_{n \to \infty} \inf_{\widehat{U} \in U(\mathcal{I}_{\text{Deg}}, \alpha)} \sup_{\Theta \in \mathcal{U}_s(\theta)} \mathbb{E}_\Theta\big[\widehat{U} - \mathcal{I}_{\text{Deg}}(\Theta)\big] \geq s(1 - \alpha). \qquad (3.7)$$

We can achieve the lower bound (3.7) by a naive upper bound $\widehat{U} = s$ with probability $1 - \alpha$ and $\widehat{U} = 0$ with probability $\alpha$. Theorem 3.1 explains why we simply choose the confidence upper bound to be $I_U^*$ for maximum degree. Similar results hold for many other invariant, we refer to Theorems 4.3 for details.

Since graph properties are a special type of graph invariants with range $[0, 1]$, Algorithm 1 induces a test for $H_0 : \mathcal{P}(G) = 0$ versus $H_1 : \mathcal{P}(G) = 1$. Using the input $\{\widehat{\boldsymbol{\Theta}}_e^d\}_{e \in V \times V}$ and the range $[0, 1]$ for Algorithm 1, we construct the test at a significance level $\alpha$ as $\psi_\alpha = 0$ if the output $\widehat{\mathcal{I}}_L = 0$ and $\psi_\alpha = 1$ otherwise. Following this idea, Algorithm 2 summarizes the procedure for testing graph properties.

An important step in the implementation of Algorithms 1 and 2 is scanning through the critical edge set $\mathcal{C}_\mathcal{I}(E_0)$ for any given $E_0$. One general procedure for doing this is to exhaustively search all edges in $\mathcal{C}_\mathcal{I}(E_0)$. Specifically, for all edge set $E' \supseteq E_0$, we search for every edge $e \in E' \backslash E_0$ and if $\mathcal{I}(E') > \mathcal{I}(E' \backslash \{e\})$, we add the edge $e$ into $\mathcal{C}_\mathcal{I}(E_0)$. Although this is a valid procedure for any monotone invariant or property, the computational complexity is exponential to the dimension. On the other hand, for many specific invariants and properties, the structure of $\mathcal{C}_\mathcal{I}(E_0)$ is simple and fast algorithms finding $\mathcal{C}_\mathcal{I}(E_0)$ exist. In Section 3.2, we will give specific algorithms for finding the critical edge sets for three invariants and one property: the maximum degree, the negative number of isolated nodes, the negative number of connected subgraphs and the property that a graph is not a forest.

---

**Algorithm 2** Skip-Down Method for Testing a Graph Property $\mathcal{P}$

---

**Input:** $\{\widehat{\boldsymbol{\Theta}}_e^d\}_{e \in V \times V}$
Initialize $t = 0, E_0 = \varnothing$.
**repeat**
$\quad t \leftarrow t + 1$;
$\quad$ Select the screening set: $\mathcal{A} \leftarrow \mathcal{C}_\mathcal{P}(E_{t-1})$;
$\quad$ Update the rejected set: $E_t \leftarrow E_{t-1} \cup \{e \in \mathcal{A} \mid \sqrt{n} \cdot |\widehat{\boldsymbol{\Theta}}_e^d| > c(\alpha, \mathcal{A})\}$;
**until** $\mathcal{P}(E_t) = 1$ or $E_t = E_{t-1}$
**Output:** $\psi_\alpha = 0$ if $\mathcal{P}(E_t) = 0$ and $\psi_\alpha = 1$ otherwise.

---

## 3.2 Case Study for Skip-Down Algorithm

In this section, we provide examples of graph invariants and show how to implement the skip-down algorithm for testing nested hypotheses and constructing confidence intervals of the following graph invariants:

- $\mathcal{I}_{\text{Conn}}(G) = $ the negative number of connected subgraphs of $G$;

- $\mathcal{I}_{\text{Deg}}(G) = $ the maximum degree of $G$;

- $\mathcal{I}_{\text{Iso}}(G) = $ the negative number of isolated nodes of $G$.

Recall that the induced graph property $\mathcal{P}_{\mathcal{I},k}(G) = 0$ if $\mathcal{I}(G) \leq k$ and $\mathcal{P}_{\mathcal{I},k}(G) = 1$ if $\mathcal{I}(G) > k$. The above four invariants therefore naturally induce the following properties:

- $\mathcal{P}_{\text{Conn},-k}(G) = 1$ if and only if $G$ has less than $k$ connected subgraphs;

- $\mathcal{P}_{\text{Deg},k}(G) = 1$ if and only if $G$ has maximum degree larger than $k$;

- $\mathcal{P}_{\text{Sig},-k}(G) = 1$ if and only if $G$ has less than $k$ isolated nodes.

We also consider a fourth example of graph property:

- $\mathcal{P}_{\text{Cycle}}(G) = 0$ if and only if $G$ is a forest.

In order to apply Algorithms 1 and 2 for the above invariants and properties, there are two important steps; we need to first specify critical edge set $\mathcal{C}_{\mathcal{I}}(E_0)$ given any edge set $E_0$, and second to calculate the invariant for any given graph. The second step can directly utilize existing algorithms for deterministic graphs. As mentioned previously, the critical edge set can be found via exhaustive search but the computation is not efficient in general. However, for the examples we give above, the critical edge sets can be selected explicitly. In the following proposition, we show the explicit forms of the critical edge sets for the above four properties.

**Proposition 3.2** (Critical edge sets for induced properties). Given any graph $G_0 = (V, E_0)$, we have the following concrete forms of the critical edge sets.

- **Connected subgraphs.** Denote all connected subgraphs of $G_0$ by $\{G_{0\ell} = (V_{0\ell}, E_{0\ell})\}_{\ell=1}^{k'}$. If $\mathcal{P}_{\text{Conn},-k}(G_0) = 1$, i.e., $k' < k$, we have $\mathcal{C}_{\mathcal{P}_{\text{Conn},-k}}(E_0) = \varnothing$, otherwise

$$\mathcal{C}_{\mathcal{P}_{\text{Conn},-k}}(E_0) = \big\{(u,v) \in E_0^c \,\big|\, u \in V_{0\ell}, v \in V_{0\ell'}, \ell \neq \ell'\big\}. \tag{3.8}$$

Thus the critical edge set consists of all edges that link the connected subgraphs of $G_0$.

- **Maximum degree.** If $\mathcal{P}_{\text{Deg},k}(G_0) = 1$, i.e., the maximum degree of $G_0$ is larger than $k$, we have $\mathcal{C}_{\mathcal{P}_{\text{Deg},k}}(E_0) = \varnothing$, otherwise $\mathcal{C}_{\mathcal{P}_{\text{Deg},k}}(E_0) = E_0^c$.

- **Singletons.** Denote the set of all isolated nodes of $G_0$ by $V_{\text{Sig}}$. If $\mathcal{P}_{\text{Sig},-k}(G_0) = 1$, i.e., $|V_{\text{Sig}}| < k$, we have $\mathcal{C}_{\mathcal{P}_{\text{Sig},-k}}(E_0) = \varnothing$, otherwise

$$\mathcal{C}_{\mathcal{P}_{\text{Sig},-k}}(E_0) = \{(u,v) \in E_0^c \,|\, u \in V_{\text{Sig}} \text{ or } v \in V_{\text{Sig}}\}. \tag{3.9}$$

Therefore, the critical edge set contains the edges that connect the isolated nodes.

- **Acyclic.** If $\mathcal{P}_{\text{Cycle}}(G) = 1$, we have $\mathcal{C}_{\mathcal{P}_{\text{Cycle}}}(E_0) = \varnothing$, otherwise when $G_0$ is a forest, $\mathcal{C}_{\mathcal{P}_{\text{Cycle}}}(E_0) = E_0^c$.

The critical edge set for an invariant can be obtained from the one of its induced property. In fact, there exists a direct connection between the critical edge set of an invariant and its induced property:

$$\mathcal{C}_{\mathcal{I}}(E_0) = \bigcup_{k=I_L^*}^{I_U^*} \mathcal{C}_{\mathcal{P}_{\mathcal{I},k}}(E_0). \tag{3.10}$$

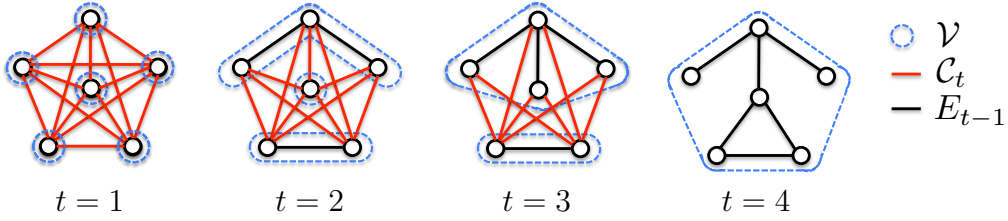

Figure 3: An example of Algorithm 3 for $\mathcal{I}_{\text{Conn}}$ from the 1st to 4th iteration. The blue dashed circles are the connected node sets $\mathcal{V}$, the red edges are the critical edge set $\mathcal{C}_t$ and the black edges are $E_{t-1}$. In the 4th iteration, since $|\mathcal{V}| = 1$, we stop the algorithm.

This can be directly derived from Definition 3.1. On the one hand, for any $e \in \mathcal{C}_\mathcal{I}(E_0)$, there exists $E' \supseteq E_0$ such that $\mathcal{I}(E') > \mathcal{I}(E'\backslash\{e\})$, therefore $e \in \mathcal{C}_{\mathcal{P}_{\mathcal{I},\mathcal{I}(E'\backslash\{e\})}}(E_0)$ by (3.4). On the other hand, if $e \in \mathcal{C}_{\mathcal{P}_{\mathcal{I},k_0}}(E_0)$ for some $k_0 \in [I_L^*, I_U^*]$, by (3.4) and the definition of induced property, there exists $E' \supseteq E_0$ such that $\mathcal{I}(E') > k_0 \geq \mathcal{I}(E'\backslash\{e\})$, which implies that $e \in \mathcal{C}_\mathcal{I}(E_0)$. Therefore, by (3.10), we have a corollary of Proposition 3.2 on the critical edge set of invariants.

**Corollary 3.3** (Critical edge sets for invariants). Given a graph $G_0 = (V, E_0)$, as in Proposition 3.2, we denote all the connected subgraphs of $G_0$ as $\{G_{0\ell} = (V_{0\ell}, E_{0\ell})\}_{\ell=1}^{k'}$ and all the set of isolated nodes of $G_0$ as $V_{\text{Sig}}$. The critical edge sets for the four invariants are as follows.

- **Connected subgraphs:** $\mathcal{C}_{\mathcal{I}_{\text{Conn}}}(E_0) = \{(u,v) \in E_0^c \,|\, u \in V_{0\ell}, v \in V_{0\ell'}, \ell \neq \ell'\}$;
- **Maximum degree:** $\mathcal{C}_{\mathcal{I}_{\text{Deg}}}(E_0) = E_0^c$;
- **Singletons:** $\mathcal{C}_{\mathcal{I}_{\text{Iso}}}(E_0) = \{(u,v) \in E_0^c \,|\, u \in V_{\text{Sig}} \text{ or } v \in V_{\text{Sig}}\}$.

The examples of the critical edges for invariants are visualized in Figure 2. Having the explicit forms of critical edge sets, we are now ready to implement Algorithms 1 and 2 for the above examples.

**Example 3.1** (Number of connected subgraphs). By Corollary 3.3, we need to partition the node set $V$ into disjoint connected node sets $\mathcal{V} = \{V_1, V_2, \ldots, V_{k'}\}$ representing connected subgraphs in each iteration. We update the connected node sets $V_1, V_2, \ldots, V_{k'}$ in each iteration of skip-down algorithm as follows. If an edge $(u,v)$ is rejected and $u \in V_1$ and $v \in V_2$, we take the union of $V_1$ and $V_2$ and update the connected node sets as $\mathcal{V} = \{V_1 \cup V_2, V_3, \ldots, V_{k'}\}$. A detailed description of the algorithm for the confidence interval of $\mathcal{I}_{\text{Conn}}$ with a default range $[-d, -1]$ is shown in Algorithm 3. We visualize the algorithm in Figure 3. Similarly, we can also test $\mathcal{P}_{\text{Conn},-k}$ by modifying Algorithm 3.

**Example 3.2** (Maximum degree). According to Corollary 3.3, the critical edge set is simply the complement of the rejected edge set. Therefore, for $\mathcal{I}_{\text{Deg}}$, in the $t$-th iteration of Algorithm 1, we select the screening set $\mathcal{A} = E_{t-1}^c$. In order to obtain $\mathcal{I}(E_t)$, we can directly count the maximal number of neighbors for each node. Therefore, the algorithm for maximum degree confidence interval can directly implement Algorithm 1 by plugging in

**Algorithm 3** Skip-Down Method for $\mathcal{I}_{\text{Conn}}$

**Input:** $\{\widehat{\boldsymbol{\Theta}}_e^d\}_{e \in V \times V}$.
Initialize $t = 0, E_0 = \varnothing$ and connected node sets $\mathcal{V} = \{\{1\}, \ldots, \{d\}\}$.
**repeat**
  $t \leftarrow t + 1$;
  Find the critical edge set $\mathcal{C}_t = \left\{(u, v) \in E_{t-1}^c \,\middle|\, u \in S, v \in T, S \neq T \text{ and } S, T \in \mathcal{V}\right\}$.
  Update the rejected set: $\mathcal{R}_t = \{e \in \mathcal{C}_t \,|\, \sqrt{n} \cdot |\widehat{\boldsymbol{\Theta}}_e^d| > c(\alpha, \mathcal{C}_t)\}$;
  $E_t \leftarrow E_{t-1} \cup \mathcal{R}_t$;
  Update the connected node sets $\mathcal{V}$:
  **for** $(u, v) \in \mathcal{R}_t$ **do**
    **if** $u, v$ belong to different node sets $S, T$ in $\mathcal{V}$, i.e., $u \in S, v \in T$ and $S \neq T$ **then**
      $\mathcal{V} \leftarrow (\mathcal{V} \backslash \{S, T\}) \cup \{S \cup T\}$
    **end if**
  **end for**
**until** $|\mathcal{V}| = 1$ or $\mathcal{R}_t = \varnothing$
**Output:**
• **Nested hypotheses:** We reject $H_{0k}$ for $k < -|\mathcal{V}|$ and not reject $H_{0k}$ for $k \in [-|\mathcal{V}|, -1]$.
• **Confidence interval:** The $1 - 2\alpha$ confidence interval $[\widehat{\mathcal{I}}_L, \widehat{\mathcal{I}}_U]$, where $\widehat{\mathcal{I}}_L = -|\mathcal{V}|$ and $\widehat{\mathcal{I}}_U = -1$ with probability $1 - \alpha$ and $\widehat{\mathcal{I}}_U = -|\mathcal{V}|$ with probability $\alpha$.

the explicit forms of $\mathcal{A} = E_{t-1}^c$ and $\mathcal{I}(E_t)$. Similar methods can be applied to test $\mathcal{P}_{\text{Deg},k}$ by plugging in the set $\mathcal{C}_{\mathcal{P}_{\text{Deg},k}}(E_0)$ defined in Proposition 3.2 into Algorithm 2.

**Example 3.3** (Singletons). The critical edge set for $\mathcal{C}_{\mathcal{I}_{\text{Sig},-k}}(E_0)$ is shown in (3.9). It is apparent from this explicit form that we need to keep the track of the set of isolated nodes for each iteration. If an edge $(u, v)$ is rejected, we simply delete $u$ and $v$ from the set of isolated nodes if they belong to the set. A detailed description of step-down algorithm for $\mathcal{I}_{\text{Sig},-k}$ with a default range $[-d, 0]$ is shown in Algorithm 4.

**Example 3.4** (Acyclic graph). The critical edge set is also the complement of the rejected edge set for the cyclicity property. Similarly to the maximum degree test, in the $t$th iteration, we also select the screening set $\mathcal{A} = E_{t-1}^c$. The procedure for checking whether $E_t \in \mathcal{P}_{\text{Cycle}}$, i.e., $E_t$ is contains a cycle is similar the one of detecting connectivity in Example 3.1. A detailed implementation of the test is shown in Algorithm 5.

# 4 Theory of Skip-down Method

In this section, we first prove the validity of obtained tests and confidence intervals described in Algorithm 1. We then show the optimality of the length of the confidence intervals for a family of monotone invariants. We show that the length of the confidence interval is adaptive to the signal strength.

**Algorithm 4** Skip-Down Method for $\mathcal{I}_{\text{Iso}}$
___
**Input:** $\{\widehat{\Theta}_e^d\}_{e \in V \times V}$.
Initialize $t = 0, E_0 = \varnothing$ and isolated node set $V_{\text{Sig}} = \{1, \ldots, d\}$.
**repeat**
   $t \leftarrow t + 1$;
   Find the critical edge set $\mathcal{C}_t = \{(u,v) \in E_{t-1}^c \,|\, u \in V_{\text{Sig}} \text{ or } v \in V_{\text{Sig}}\}$.
   Update the rejected set: $\mathcal{R}_t = \{e \in \mathcal{C}_t \,|\, \sqrt{n} \cdot |\widehat{\Theta}_e^d| > c(\alpha, \mathcal{C}_t)\}$;
   $E_t \leftarrow E_{t-1} \cup \mathcal{R}_t$;
   Update the connected node sets $V_{\text{Sig}}$:
   **for** $(u,v) \in \mathcal{R}_t$ **do**
     $V_{\text{Sig}} \leftarrow V_{\text{Sig}} \setminus \{u,v\}$;
   **end for**
**until** $|V_{\text{Sig}}| = 0$ or $\mathcal{R}_t = \varnothing$
**Output:**
- **Nested hypotheses:** We reject $H_{0k}$ for $k < -|V_{\text{Sig}}|$ and not reject $H_{0k}$ for $k \in [-|\mathcal{V}|, 0]$.
- **Confidence interval:** The $1 - 2\alpha$ confidence interval $[\widehat{\mathcal{I}}_L, \widehat{\mathcal{I}}_U]$, where $\widehat{\mathcal{I}}_L = -|V_{\text{Sig}}|$ and $\widehat{\mathcal{I}}_U = 0$ with probability $1 - \alpha$ and $\widehat{\mathcal{I}}_U = -|V_{\text{Sig}}|$ with probability $\alpha$.
___

## 4.1 Validity of Tests and Confidence Intervals

Given the precision matrix $\Theta$, recall that $G(\Theta)$ is the graph induced by the support of $\Theta$. We also use the shorthand $\mathcal{I}(\Theta) = \mathcal{I}(G(\Theta))$ for the corresponding invariant for $G(\Theta)$ and similarly put $\mathcal{P}(\Theta) = \mathcal{P}(G(\Theta))$ for the corresponding property. Given a graph invariant $\mathcal{I}$, we define the parameter space

$$\mathcal{U}_s(I_L^*, I_U^*) = \{\Theta \in \mathcal{U}_s \,|\, \mathcal{I}(\Theta) \in [I_L^*, I_U^*]\}.$$

Here we assume $I_L^*$ and $I_U^*$ are known.

The following theorem shows the family-wise error of nested hypotheses and the asymptotic coverage probability of the confidence interval given by the skip-down method in Algorithm 1.

**Theorem 4.1** (Asymptotic coverage probability). *Suppose $\Theta \in \mathcal{U}_s$ and (3.5) is satisfied. Given any monotone invariant $\mathcal{I}$, the multiple test given by Algorithm 1 has*

$$\limsup_{n \to \infty} \sup_{\Theta \in \mathcal{U}_s(I_L^*, I_U^*)} \mathbb{P}_\Theta\big(\exists k \geq \mathcal{I}(\Theta) \text{ such that } H_{0k} \text{ is rejected}\big) \leq \alpha, \quad (4.1)$$

*and the confidence interval $[\widehat{\mathcal{I}}_L, I_U^*]$ satisfies*

$$\liminf_{n \to \infty} \inf_{\Theta \in \mathcal{U}_s(I_L^*, I_U^*)} \mathbb{P}_\Theta\big(\mathcal{I}(\Theta) \in [\widehat{\mathcal{I}}_L, I_U^*]\big) \geq 1 - \alpha. \quad (4.2)$$

Since graph properties are a special types of graph invariants, Theorem 4.1 actually implies the validity of the test $\psi_\alpha$ in Algorithm 2.

**Algorithm 5** Skip-Down Method for $\mathcal{P}_{\text{Cycle}}$

---

**Input:** $\{\widehat{\boldsymbol{\Theta}}_e^d\}_{e \in V \times V}$.
Initialize $t = 0, E_0 = \varnothing$ and path node sets $\mathcal{V}_{\text{Paths}} = \{\{1\}, \ldots, \{d\}\}$.
**repeat**
   $t \leftarrow t + 1$;
   Update the rejected set: $\mathcal{R}_t = \{e \in E_t^c \,|\, \sqrt{n} \cdot |\widehat{\boldsymbol{\Theta}}_e^d| > c(\alpha, E_t^c)\}$;
   $E_t \leftarrow E_{t-1} \cup \mathcal{R}_t$;
   Detect the cycles in the graph:
   **for** $(u, v) \in \mathcal{R}_t$ **do**
     **if** $u, v$ belong to different node sets $S, T$ in $\mathcal{V}_{\text{Paths}}$, i.e., $u \in S, v \in T$ and $S \neq T$ **then**
       $\mathcal{V}_{\text{Paths}} \leftarrow (\mathcal{V}_{\text{Paths}} \backslash \{S, T\}) \cup \{S \cup T\}$
     **else**
       $\mathcal{V}_{\text{Paths}} \leftarrow \varnothing$, **Break;**
     **end if**
   **end for**
**until** $\mathcal{R}_t = \varnothing$
**Output:** $\psi_\alpha = 0$ if $\mathcal{V}_{\text{Paths}} \neq \varnothing$ and $\psi_\alpha = 1$ otherwise.

---

**Corollary 4.2** (Uniform asymptotic validitiy). Given any monotone property $\mathcal{P}$, we define the parameter space

$$\mathcal{G}_0 := \{\boldsymbol{\Theta} \in \mathcal{U}_s \,|\, \mathcal{P}(\boldsymbol{\Theta}) = 0\}. \tag{4.3}$$

Under the same conditions of Theorem 4.1, the test $\psi_\alpha$ given by Algorithm 2 has

$$\liminf_{n \to \infty} \inf_{\boldsymbol{\Theta} \in \mathcal{G}_0} \mathbb{P}_{\boldsymbol{\Theta}}(\psi_\alpha = 0) \geq 1 - \alpha. \tag{4.4}$$

Both results of Theorem 4.1 and Corollary 4.2 are uniform over the parameter space $\mathcal{U}_s(I_L^*, I_U^*)$. Besides, the parameter space $\mathcal{U}_s(I_L^*, I_U^*)$ does not impose any restriction on the signal strength. The signal strength assumption is usually required to obtain a consistent graph estimator (Ravikumar et al., 2011; Cai et al., 2011; Liu et al., 2012b), such that $\widehat{G}$ satisfies $\mathbb{P}(\widehat{G} = G) \to 1$ as $n \to \infty$. The consistency of graph recovery has been shown for many graphical model estimators including CLIME (Cai et al., 2011), neighborhood selection (Meinshausen and Bühlmann, 2006; Zhou et al., 2009), graphical Lasso (Lam and Fan, 2009) and transellipitical graphical model (Liu et al., 2012b) and all these results need the minimal signal strength condition. Given any consistent graph estimator $\widehat{G}$, there is a plug-in invariant estimator $\mathcal{I}(\widehat{G})$ such that $\mathbb{P}(\mathcal{I}(\widehat{G}) = \mathcal{I}(G)) \to 1$. Similarly, to test the properties, one can also construct a plug-in test $\psi := \mathcal{P}(\widehat{G})$. In comparison with the plug-in methods, the confidence intervals and tests obtained by Algorithms 1 and 2 have two major advantages: (1) they do not require any signal strength conditions and (2) one can choose different significance levels and thus have controllable uncertainty assessment for our inferential procedures.

Neykov et al. (2016) propose a data-splitting method to test the graph properties. Their method splits the data into two disjoint parts $\mathcal{D}_1$ and $\mathcal{D}_2$. The method uses the first data split $\mathcal{D}_1$ to estimate the inverse covariance $\widehat{\Theta}(\mathcal{D}_1)$ and order the estimators $|\widehat{\Theta}_{[1]}| \geq \ldots \geq |\widehat{\Theta}_{[d(d-1)/2]}|$. It adds the edges from the largest to the smallest until they detect witness pattern $\widehat{E}$ for the property. For example, for the property that a graph has maximum degree larger than $k$, the pattern $\widehat{E}$ is the first star of size $k+1$. The method then applies the step-down algorithm on the second data split $\mathcal{D}_2$ to test on the existence of edges in $\widehat{E}$. If all edges in $\widehat{E}$ are rejected, the method rejects the hypothesis and does not reject otherwise. Compared to the data splitting test proposed in Neykov et al. (2016), the skip-down algorithm has three advantages. First, it is a data-driven procedure without any pre-specified value of a graph invariant to be tested. Second, it is a unified procedure valid for all monotone graph invariants, which are non-decreasing under edge addition, while the data splitting method needs to consider different invariants case by case. Third, the data splitting method needs the first split to screen the critical edge set and conducts the tests on the critical edge set by the second split. In comparison, the skip-down algorithm avoids data splitting and conducts both screening and inference using the entire dataset.

## 4.2 Optimality and Adaptivity of the Confidence Intervals

Theorem 4.1 proves that the confidence intervals $[\widehat{\mathcal{I}}_L, I_U^*]$ constructed by Algorithm 1 is asymptotically honest uniformly over the parameter space $\mathcal{U}_s(I_L^*, I_U^*)$. In this subsection, we show that the averaged length of confidence interval obtained in Algorithm 1 is optimal in the minimax sense.

Define a parameter space with minimal signal strength $\theta$ as

$$\mathcal{U}_\mathcal{I}(I_L^*, I_U^*; \theta) = \Big\{ \Theta \in \mathcal{U}_s \;\Big|\; \mathcal{I}(\Theta) \in [I_L^*, I_U^*], \min_{e \in E(\Theta)} |\Theta_e| \geq \theta \Big\}. \tag{4.5}$$

Let the family of confidence upper bound of $\mathcal{I}$ be

$$U(\mathcal{I}, \alpha) = \Big\{ \widehat{U}(\cdot): \mathbb{R}^{d \times n} \to [I_L^*, I_U^*] \;\Big|\; \inf_{\Theta \in \mathcal{U}_s(I_L^*, I_U^*)} \mathbb{P}_\Theta \big( \mathcal{I}(\Theta) \leq \widehat{U}(\mathbf{X}) \big) \geq 1 - \alpha \Big\}.$$

Theorem 3.1 shows we cannot construct non-trivial confidence upper bound for maximum degree. The following theorem shows that a similar result holds for $\mathcal{I}_{\text{Conn}}$ and $\mathcal{I}_{\mathcal{I}_{\text{Iso}}}$ as well.

**Theorem 4.3** (Lower bound of confidence upper bounds). Suppose $I_L^* \geq d/2$ and $I_U^* - I_L^* = o(d^{1/2})$. If $\theta \leq C\sqrt{\log d/n}$ for some sufficiently small positive constant $C$, we have

$$\liminf_{n \to \infty} \inf_{\widehat{U} \in U(\mathcal{I}, \alpha)} \sup_{\Theta \in \mathcal{U}_\mathcal{I}(I_L^*, I_U^*; \theta)} \mathbb{E}_\Theta \big[ \widehat{U} - \mathcal{I}(\Theta) \big] \geq (I_U^* - I_L^*)(1 - \alpha), \tag{4.6}$$

for $\mathcal{I} = \mathcal{I}_{\text{Conn}}$ or $\mathcal{I} = \mathcal{I}_{\text{Iso}}$.

Similar to Theorem 3.1, the lower bound in (4.6) can be achieved by a naive upper side $\widehat{U} = I_U^*$ with probability $1 - \alpha$ and $\widehat{U} = I_L^*$ with probability $\alpha$. This explains why in Algorithm 1, we only simply choose the upper bound as $I_U^*$. For the result of generic invariant confidence upper bound, we refer to Theorem 6.2.

Although we cannot construct a non-trivial confidence upper bound for many invariants, we can still construct good enough confidence intervals. Define the family of honest confidence intervals as

$$I(\mathcal{I}, \alpha) = \left\{ [\widehat{L}, \widehat{U}] \mid \widehat{L}(\mathbf{X}) \leq \widehat{U}(\mathbf{X}) \text{ a.s.}, \inf_{\boldsymbol{\Theta} \in \mathcal{U}_s(I_L^*, I_U^*)} \mathbb{P}_{\boldsymbol{\Theta}}\big(\mathcal{I}(\boldsymbol{\Theta}) \in [\widehat{L}(\mathbf{X}), \widehat{U}(\mathbf{X})]\big) \geq 1 - \alpha \right\}. \quad (4.7)$$

The following theorem gives the lower bounds of confidence interval length for $\mathcal{I} = \mathcal{I}_{\text{Deg}}$, $\mathcal{I}_{\text{Conn}}$ and $\mathcal{I}_{\text{Iso}}$.

**Theorem 4.4** (Lower bound of confidence intervals). Given a precision matrix $\boldsymbol{\Theta} \in \mathcal{U}_{\mathcal{I}}(I_L^*, I_U^*; \theta)$, we define its significant edge set as

$$E_{\text{Sig}}(\boldsymbol{\Theta}) := \big\{ (j, k) \mid |\boldsymbol{\Theta}_{jk}| \geq C\sqrt{\log d/n} \big\},$$

where $C$ is some sufficiently large constant. We then define the oracle length as

$$\text{Oracle Length}(\boldsymbol{\Theta}) = I_U^* - \mathcal{I}(E_{\text{Sig}}(\boldsymbol{\Theta})). \quad (4.8)$$

For $\mathcal{I} = \mathcal{I}_{\text{Deg}}$, $\mathcal{I}_{\text{Conn}}$ or $\mathcal{I}_{\text{Iso}}$, if $\theta \leq C'\sqrt{\log d/n}$ for some sufficiently small positive constant $C'$, we have

$$\liminf_{n \to \infty} \inf_{[\widehat{L}, \widehat{U}] \in I(\mathcal{I}, \alpha)} \sup_{\boldsymbol{\Theta} \in \mathcal{U}_s(I_L^*, I_U^*; \theta)} \frac{\mathbb{E}_{\boldsymbol{\Theta}}[\widehat{U} - \widehat{L}]}{\text{Oracle Length}(\boldsymbol{\Theta})} \geq 1 - 2\alpha. \quad (4.9)$$

By the definition in (4.8), Oracle Length($\boldsymbol{\Theta}$) becomes smaller, if there are more entries in $\boldsymbol{\Theta}$ with signal strength larger than $C\sqrt{\log d/n}$. Namely, the oracle length is adaptive to the number of edges to strong signal strength. Therefore, (4.9) implies that the lower bound of the confidence interval length is adaptive to the significant edge set $E_{\text{Sig}}(\boldsymbol{\Theta})$. For the lower bound of generic invariants, we refer to Theorem 6.1 in Section 6.

From (4.9), it is straightforward to derive that

$$\liminf_{n \to \infty} \inf_{[\widehat{L}, \widehat{U}] \in I(\mathcal{I}, \alpha)} \sup_{\boldsymbol{\Theta} \in \mathcal{U}_{\mathcal{I}}(I_L^*, I_U^*; \theta)} \mathbb{E}_{\boldsymbol{\Theta}}[\widehat{U} - \widehat{L}] \geq (1 - 2\alpha)(I_U^* - I_L^*).$$

This implies that if there is no edge with strong enough signal strength, we only have a trivial rate $O(I_U^* - I_L^*)$ for the confidence interval length.

Now we discuss the upper bound the confidence interval length from Algorithm 1. The following theorem shows that it achieves the lower bound in Theorem 4.4.

**Theorem 4.5** (Size of confidence interval). Suppose $\boldsymbol{\Theta} \in \mathcal{U}_s$ and (3.5) is satisfied. For any monotone invariant $\mathcal{I}$ with range $[I_L^*, I_U^*]$, if $I_U^* - I_L^* = O(d^2)$, for any $\alpha \in (0, 1)$ and $\theta > 0$,

$$\lim_{n \to \infty} \sup_{\boldsymbol{\Theta} \in \mathcal{U}_{\mathcal{I}}(I_L^*, I_U^*; \theta)} \frac{\mathbb{E}_{\boldsymbol{\Theta}}[\widehat{\mathcal{I}}_U - \widehat{\mathcal{I}}_L]}{\text{Oracle Length}(\boldsymbol{\Theta}) + 1} \leq 1. \quad (4.10)$$

We add one in the denominator Oracle Length($\boldsymbol{\Theta}$) + 1 of (4.10) just to avoid singularity when $I_U^* = \mathcal{I}(E_{\text{Sig}}(\boldsymbol{\Theta}))$. We can see that the length is adaptive to the value $\mathcal{I}(E_{\text{Sig}}(\boldsymbol{\Theta}))$. As we argued above, the oracle length in (4.10) shows the first level of adaptivity for the skip-down algorithm: the length of our confidence interval is smaller if there are more edges with strong enough signal strength. The assumption that $I_U^* - I_L^* = O(d^2)$ is satisfied for all examples in Section 2. In fact, this assumption is mild in the sense that for monotone $\mathcal{I}$, there are at most $d(d-1)/2$ possible values and we can easily rescale $\mathcal{I}$ such that $I_U^* - I_L^* = O(d^2)$.

Theorem 4.3 shows that it is impossible to construct an adaptive upper side of confidence interval. In fact, the following theorem shows that the oracle length in (4.10) mainly comes from the lower side.

**Theorem 4.6** (Size of confidence lower bound). Suppose $\boldsymbol{\Theta} \in \mathcal{U}_s$. If (3.5) is satisfied, for any monotone invariant $\mathcal{I}$ with range $[I_L^*, I_U^*]$, if $I_U^* - I_L^* = O(d^2)$, for any $\alpha \in (0,1)$ and $\theta > 0$, we have

$$\lim_{n\to\infty} \sup_{\boldsymbol{\Theta} \in \mathcal{U}_\mathcal{I}(I_L^*, I_U^*; \theta)} \frac{\mathbb{E}_{\boldsymbol{\Theta}}[\mathcal{I}(\boldsymbol{\Theta}) - \widehat{\mathcal{I}}_L]}{\mathcal{I}(\boldsymbol{\Theta}) - \mathcal{I}(E_{\text{Sig}}(\boldsymbol{\Theta})) + 1} \leq 1. \tag{4.11}$$

Similar to (4.10), we add one in the denominator $I_U^* - \mathcal{I}(E_{\text{Sig}}(\boldsymbol{\Theta})) + 1$ of (4.11) just to avoid singularity. We remark that (4.11) also gives the type II error analysis for nested hypotheses in (3.2). In Algorithm 1, we do not reject $H_{0k}$ for $k \in (\widehat{\mathcal{I}}_L, I_U^*]$, thus the number of type II errors is $\max\{\mathcal{I}(\boldsymbol{\Theta}) - \widehat{\mathcal{I}}_L, 0\}$. In fact, the proof of Theorem 4.6 shows that the expected number of type II errors has

$$\lim_{n\to\infty} \mathbb{E}_{\boldsymbol{\Theta}}[\mathcal{I}(\boldsymbol{\Theta}) - \widehat{\mathcal{I}}_L] \leq \mathcal{I}(\boldsymbol{\Theta}) - \mathcal{I}(E_{\text{Sig}}(\boldsymbol{\Theta})).$$

Therefore, when the minimal signal strength satisfies $\min_{(j,k)\in E} |\boldsymbol{\Theta}_{jk}| \geq C\sqrt{\log d/n}$ for sufficiently large $C$, there is asymptotically no type II error.

## 5 Extension to Non-Gaussian Graphical Models

In this section, we show that the skip-down method for nested hypotheses and confidence intervals can be applied to a general family of graphical models and estimators.

We first introduce a few notations for general graphical models. We say a symmetric matrix $\boldsymbol{\Omega}$ is an interaction matrix for a graphical model Markov to graph $G^* = (V, E^*)$, if for any two different $j, k \in V$, $\boldsymbol{\Omega}_{jk} = 0$ if and only if $(j,k) \notin E^*$. For the Gaussian graphical model, the precision matrix $\boldsymbol{\Theta}$ is an interaction matrix by definition. For Gaussian copula graphical model (Liu et al., 2012a; Xue and Zou, 2012) and the transelliptical graphical model (Liu et al., 2012b), the interaction matrix can be the latent correlation matrix. For the semiparametric exponential family graphical model (Yang et al., 2014b), the interaction matrix can be the interactive parameters in the canonical form.

In the skip-down algorithm for Gaussian graphical model in Algorithm 1, we only use the debiased precision matrix estimators $\widehat{\boldsymbol{\Theta}}^d$ in the step of updating the rejected edge set. In fact, for general graphical models, it suffices to validate Algorithm 1 if we can construct some estimator $\widehat{\boldsymbol{\Omega}}^d$ for the interaction matrix such that the quantile of statistic $\max_{(j,k)\in E} \sqrt{n}(\widehat{\boldsymbol{\Omega}}^d_{jk} - \boldsymbol{\Omega}_{jk})$ can be well estimated for any edge set $E \subseteq V \times V$. In specific, we have the following theorem for general graphical models.

**Theorem 5.1.** Let $\boldsymbol{\Omega}$ be an interaction matrix for a graphical model. Suppose there exist an interaction matrix estimator $\widehat{\boldsymbol{\Omega}}^d$ and a quantile estimator $c(\alpha, E)$ for the statistic $\max_{(j,k)\in E} \sqrt{n}(\widehat{\boldsymbol{\Omega}}^d_{jk} - \boldsymbol{\Omega}_{jk})$ given any edge set $E \subseteq V \times V$ satisfying

$$\lim_{n\to\infty} \sup_{\boldsymbol{\Omega}\in\mathcal{U}_s} \left| \mathbb{P}_{\boldsymbol{\Omega}}\Big( \max_{e\in E} \sqrt{n}(\widehat{\boldsymbol{\Omega}}^d_e - \boldsymbol{\Omega}_e) > c(\alpha, E) \Big) - \alpha \right| = 0.$$

For any monotone invariant $\mathcal{I}$, if we plug such estimators $\widehat{\boldsymbol{\Omega}}^d$ and $c(\alpha, E)$ into Algorithm 1, the output test is honest, i.e.,

$$\limsup_{n\to\infty} \sup_{\boldsymbol{\Omega}\in\mathcal{U}_s(I_L^*, I_U^*)} \mathbb{P}_{\boldsymbol{\Omega}}\big(\exists k \geq \mathcal{I}(\boldsymbol{\Omega}) \text{ such that } H_{0k} \text{ is rejected}\big) \leq \alpha, \tag{5.1}$$

and the output confidence interval $[\widehat{\mathcal{I}}_L, I_U^*]$ satisfies

$$\liminf_{n\to\infty} \inf_{\boldsymbol{\Omega}\in\mathcal{U}_s(I_L^*, I_U^*)} \mathbb{P}_{\boldsymbol{\Omega}}\big(\mathcal{I}(\boldsymbol{\Omega}) \in [\widehat{\mathcal{I}}_L, I_U^*]\big) \geq 1 - \alpha. \tag{5.2}$$

The above theorem generalizes Theorem 4.1 to non-Gaussian graphical models. In fact, it provides a sufficient condition (3.5) to guarantee that the test and confidence interval achieved from Algorithm 1. In fact, after applying a debiasing step, (3.5) is satisfied for many estimators for Gaussian graphical model as long as non-Gaussian graphical models including Gaussian copula graphical model, transelliptical graphical model, semiparametric exponential family graphical model and so on. For the details of debiasing step and approaching (3.5), we refer Section S.5 in the appendix.

## 6 A Generic Framework of Lower Bounds

In this section, we show the lower bound of the confidence interval length for a generic family of invariants characterized by a concept called "hollow graphs", whose formal definition is given below.

**Definition 6.1** (Hollow graph)**.** A graph $G = (V, E)$ is called $R$-hollow, if

$$\max_{\varnothing \neq F \subseteq E} \frac{|F|}{|V(F)| - 1} \leq R. \tag{6.1}$$

The quantity $R$ in (6.1) measures the maximal "density" of edges. A $k$-clique is $k$-hollow as the edges are fully connected and a $k$-chain is $(1 - 1/k)$-hollow as it is relatively sparse. Therefore, we can see that any graph with maximal degree $s$ is at most $s$-hollow. Definition 6.1 is a classical definition proposed by Nash-Williams (1964). If a graph is $R$-hollow for some constant $R$ independent of the graph size, we say the graph is hollow.

Before presenting the theorem on the lower bound, we introduce a few notations on the graph. The maximum degree of a graph $G = (V, E)$ is denoted by $d_{\max}(G)$. The union of two graphs $G_1 = (V, E_1)$ and $G_2 = (V, E_2)$ is $G_1 \cup G_2 := (V, E_1 \cup E_2)$. We say $G' = (V, E')$ is an isomorphic copy of $G = (V, E)$ if $G'$ is isomorphic to $G$ and $V(E') \cap V(E) = \varnothing$. Furthermore, we say $G_1, \ldots, G_N$ are different isomorphic copies of $G$ if each $G_j$ is an isomorphic copy of $G$ and $V(E_j) \cap V(E_{j'}) = \varnothing$ for any $1 \leq j \neq j' \leq N$.

Now we briefly explain the intuition behind the proof of the lower bound of confidence interval length. Given an invariant $\mathcal{I}$, the idea of the proof is to reduce the problem to a lower bound in a certain hypothesis test (Cai and Guo, 2015). In order to obtain sharp bounds, we compare a graph with invariant equal to $I_L^*$ with multiple graphs with invariant equal to $I_U^*$. The construction of the alternative graphs with invariant $I_U^*$ relies on the existence of isomorphic copies of a certain graph.

**Theorem 6.1** (Lower bound of confidence interval length). Given any monotone invariant $\mathcal{I}$ in the range $[I_L^*, I_U^*]$, suppose there exist two graphs $G_L = (V, E_L)$ and $G_U = (V, E_U)$ with $G_L \preceq G_U$ and they satisfy $|V(E_L)| = O(1)$, $|V(E_U)| = o(d^{1/2})$ and $G_U$ is hollow. Given some $N$ satisfying $d^{1/2} \leq N \leq d/(2|V(E_L)|)$, we assume there exist $N$ different isomorphic copies of $G_L$ denoted as $G_{L,1}, \ldots, G_{L,N}$ such that

$$\mathcal{I}(\cup_{j=1}^N G_{L,j} \cup G_L) = I_L^* \text{ and } \mathcal{I}(\cup_{j=1}^N G_{L,j} \cup G_U) = I_U^*. \tag{6.2}$$

If there exist constants $C_1$ and $C_2$ such that

$$\theta \leq C_1 \sqrt{\log d/n} \text{ and } d_{\max}(G_U)\sqrt{\log d/n} \leq C_2, \tag{6.3}$$

we have the following lower bound on the confidence interval length

$$\liminf_{n \to \infty} \inf_{[\widehat{L}, \widehat{U}] \in I(\mathcal{I}, \alpha)} \sup_{\boldsymbol{\Theta} \in \mathcal{U}_s(I_L^*, I_U^*; \theta)} \frac{\mathbb{E}_{\boldsymbol{\Theta}}[\widehat{U} - \widehat{L}]}{\text{Oracle Length}(\boldsymbol{\Theta})} \geq 1 - 2\alpha, \tag{6.4}$$

where Oracle Length($\boldsymbol{\Theta}$) is defined in (4.8).

We now sketch the high level idea behind the proof of Theorem 6.1. The first step reduces the minimax result in (6.4) to a lower bound of testing $H_0 : \mathcal{I}(\boldsymbol{\Theta}) = I_L^*$ versus $H_1 : \mathcal{I}(\boldsymbol{\Theta}) = I_U^*$. Next we further reduce the above test to the test $H_0 : \boldsymbol{\Theta} = \boldsymbol{\Theta}_0$ versus $H_1 : \boldsymbol{\Theta} \in \{\boldsymbol{\Theta}_1, \ldots, \boldsymbol{\Theta}_M\}$, where $\mathcal{I}(\boldsymbol{\Theta}_0) = I_L^*$ and $\mathcal{I}(\boldsymbol{\Theta}_j) = I_L^*$ for $1 \leq j \leq M$. In order to obtain sharp bounds one needs to construct a maximally challenging set of matrices $\{\boldsymbol{\Theta}_0, \boldsymbol{\Theta}_1, \ldots, \boldsymbol{\Theta}_M\}$ for hypothesis test. Condition (6.2) enables our construction. We choose a $\boldsymbol{\Theta}_0$ with $G(\boldsymbol{\Theta}_0) = \cup_{j=1}^N G_{L,j} \cup G_L$ and select $\boldsymbol{\Theta}_1, \ldots, \boldsymbol{\Theta}_M$ so that each $G(\boldsymbol{\Theta}_j)$ for $1 \leq j \leq M$ is isomorphic to $\cup_{j=1}^N G_{L,j} \cup G_U$. The reason the graph $\cup_{j=1}^N G_{L,j} \cup G_U$ is used is to reproduce multiple isomorphic $G(\boldsymbol{\Theta}_j)$'s in the alternative which makes it challenging to tell them from the graph $G(\boldsymbol{\Theta}_0)$ under the null. Assumptions on the sizes of $G_L$ and $G_U$, as well as the range of $N$ are imposed to ensure the existence of sufficiently many graphs in the alternative. The second condition of (6.3) guarantees the positive definiteness of the precision matrix. It is satisfied when the condition $s\sqrt{\log d/n} = o(1)$ holds.

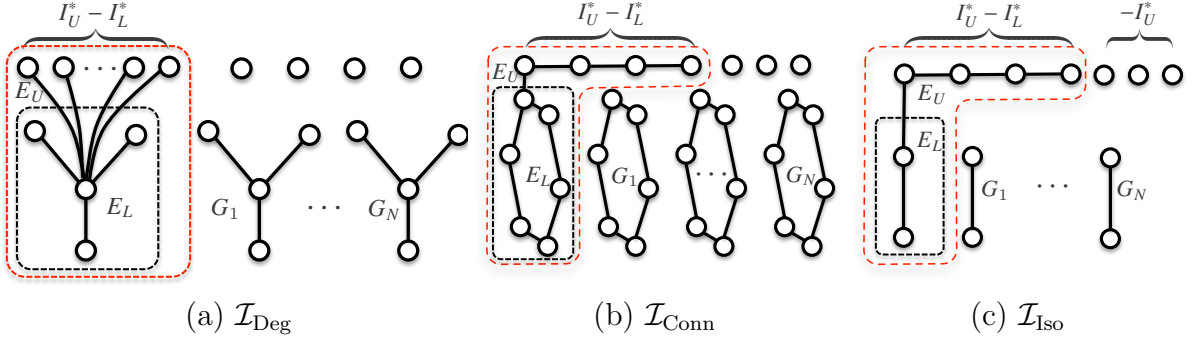

(a) $\mathcal{I}_{\text{Deg}}$  (b) $\mathcal{I}_{\text{Conn}}$  (c) $\mathcal{I}_{\text{Iso}}$

Figure 4: A visualization of a construction satisfying the conditions of Theorem 6.1 for the invariants $\mathcal{I}$: (a) the maximum degree, (b) the negative number of connected subgraphs and (c) the negative number of isolated nodes.

The following theorem gives a generic lower bound of confidence upper bound, i.e., the right endpoint of the confidence interval.

**Theorem 6.2** (Negative result of confidence upper bound)**.** Given any monotone invariant $\mathcal{I}$ in the range $[I_L^*, I_U^*]$. Under the same conditions as Theorem 6.1, if (6.3) is satisfied, we have

$$\liminf_{n\to\infty} \inf_{\widehat{U}\in U(\mathcal{I},\alpha)} \sup_{\boldsymbol{\Theta}\in\mathcal{U}_\mathcal{I}(I_L^*,I_U^*;\theta)} \mathbb{E}_{\boldsymbol{\Theta}}\big[\widehat{U} - \mathcal{I}(\boldsymbol{\Theta})\big] \geq (I_U^* - I_L^*)(1-\alpha). \tag{6.5}$$

In the final part of this section, we give concrete examples of invariants satisfying the conditions of Theorem 6.1. Specifically, the following examples show how to derive Theorems 4.3 and 4.4 from Theorems 6.1 and 6.2.

**Example 6.1.** The maximum degree $\mathcal{I}_{\text{Deg}}$ satisfies the conditions of Theorem 6.1 for $I_L^* = O(1)$ and $I_U^* = o(d^{1/2})$. We visualize the construction in Figure 4(a). We let $G_L$ be an $I_L^*$-star and $G_U$ be an $I_U^*$-star. By definition we immediately have that $G_U$ is hollow. Since $I_U^* = o(d^{1/2})$ and $I_L^* = O(1)$, we have $|V(E_L)| = O(1)$ and $|V(E_U)| = o(d^{1/2})$. We let $N = d/(2(I_L^* + 1)) = \Omega(d^{1/2})$ and reproduce $N$ different isomorphic copies of the $I_L^*$-star, as $G_{L,1}, \ldots, G_{L,N}$ shown in Figure 4(a). One can easily verify that this construction satisfies (6.2).

**Example 6.2.** The negative number of connected subgraphs $\mathcal{I}_{\text{Conn}}$ satisfies the conditions of Theorem 6.1 for any pair of $I_L^*$ and $I_U^*$ such that $I_L^* = -\lfloor \gamma d \rfloor$ for some $\gamma \in (1/2, 1]$ and $I_U^* - I_L^* = o(d^{1/2})$. We let $G_L$ be a $1/(2\gamma - 1)$-loop[2] and $G_U$ be a graph connecting an $(I_U^* - I_L^*)$-chain to $G_L$, as is shown in Figure 4(b). It is easy to see that $G_U$ is hollow, $|V(E_L)| = O(1)$ and $|V(E_U)| = o(d^{1/2})$. Choose $N = (2k-1)d/2$ and reproduce $N$ different isomorphic copies of the $1/(2\gamma - 1)$-loop. It can be checked that (6.2) is satisfied.

---

[2]Without loss of generality, we assume $1/(2\gamma - 1)$ is an integer and same for $N$.

**Example 6.3.** The negative number of isolated nodes $\mathcal{I}_{\text{Iso}}$ satisfies the conditions of Theorem 6.1 for any pair of $I_L^*$ and $I_U^*$ satisfying $I_L^* = -\lfloor \gamma d \rfloor$ for some $\gamma \in (0, 1)$ and $I_U^* - I_L^* = o(d^{1/2})$. We let $G_L$ be a 2-chain and $G_U$ be a graph connecting an $(I_U^* - I_L^*)$-chain to $G_L$, as is shown in Figure 4(c). By construction $G_U$ is hollow, $|V(E_L)| = O(1)$ and $|V(E_U)| = o(d^{1/2})$. Choose $N = (d - I_L^*)/2$ and reproduce $N$ different isomorphic copies of the 2-chain. One can easily check that (6.2) is satisfied.

# 7 Numerical Results

In this section we show the numerical performance for the proposed confidence interval to three graph invariants: the negative number of connected subgraphs, the maximum degree and the negative number of isolated nodes. In addition we apply our method to a neuroimaging dataset.

## 7.1 Synthetic Data

We first implement Algorithm 1 for synthetic simulations. Here we use the confidence interval to illustrate the performance of the skip-down algorithm. We consider three invariants: $\mathcal{I}_{\text{Conn}}$, $\mathcal{I}_{\text{Deg}}$ and $\mathcal{I}_{\text{Iso}}$. Let $\boldsymbol{X}_1, \ldots, \boldsymbol{X}_n$ be i.i.d. $n$ samples generated from $N(0, \boldsymbol{\Theta}^{-1})$, where the precision matrix $\boldsymbol{\Theta} \in \mathbb{R}^{d \times d}$. In order to illustrate Theorem 4.5 and show how the confidence interval length changes adaptively with the value of $\mathcal{I}(E_{\text{Sig}}(\boldsymbol{\Theta}))$, we choose two parameters $k \geq k_\mu$ and generate the precision matrix $\boldsymbol{\Theta}$ satisfying $\mathcal{I}(\boldsymbol{\Theta}) = k$ and $\mathcal{I}(E_{\text{Sig}}(\boldsymbol{\Theta})) = k_\mu$. To generate such a precision matrix, we set the minimal signal strength $\theta = 0.01\sqrt{\log d/n}$ and we set part of the non-zero entries of $\boldsymbol{\Theta}$ as $\mu > \theta$. To calculate the average coverage probability and average confidence interval length, we repeat the simulation 500 times. We allow the precision matrix to change in different repetitions to show that our confidence interval is honest.

Specifically, we construct the precision matrices following different scenarios for three different invariants as follows.

- $\mathcal{I}_{\text{Conn}}$: Denote the adjacency matrix of a $d$-chain with $\mathbf{A}_{\text{chain}(d)} \in \mathbb{R}^{d \times d}$. Formally, $[\mathbf{A}_{\text{chain}(d)}]_{s,t} = 1$ for any $1 \leq s \neq t \leq d$ satisfying $|s - t| = 1$ and $[\mathbf{A}_{\text{chain}(d)}]_{s,t} = 0$ otherwise. Recall that $\mathcal{I}_{\text{Conn}}$ denotes the negative number of connected subgraphs, and hence the parameter $k$ takes negative values. We construct our graph with $-k$ connected subgraphs by removing $-k - 1$ edges from the $d$-chain. Given any $k \in [-d, -1]$, let the adjacency matrix for the edges to cut be $\mathbf{A}_{\text{cut}(k)}$ such that $[\mathbf{A}_{\text{cut}(k)}]_{s,t} = 1$ if and only if $|s - t| = 1$ and $s = j\lfloor -d/k \rfloor$ for $j = 1, 2, \ldots, -k - 1$ and otherwise $[\mathbf{A}_{\text{cut}(k)}]_{s,t} = 0$. We construct the precision matrix as $\boldsymbol{\Theta}_{\mu,k} = \mathbf{I}_d + \mu(\mathbf{A}_{\text{chain}(d)} - \mathbf{A}_{\text{cut}(k)})$. Given $k$ and $k_\mu$ to be the values of $\mathcal{I}_{\text{Conn}}(\boldsymbol{\Theta})$ and $\mathcal{I}_{\text{Conn}}(E_{\text{Sig}}(\boldsymbol{\Theta}))$ respectively, for each repetition, we uniformly sample $k - k_\mu$ edges from $E(\boldsymbol{\Theta}_{\mu,k})$ and change the values of entries on $\boldsymbol{\Theta}_{\mu,k}$ corresponding to these edges from $\mu$ to $\theta$. Denote this new precision matrix as $\widetilde{\boldsymbol{\Theta}}$ and we can see that $\mathcal{I}_{\text{Conn}}(\widetilde{\boldsymbol{\Theta}}) = k$ and $\mathcal{I}_{\text{Conn}}(E_{\text{Sig}}(\widetilde{\boldsymbol{\Theta}})) = k_\mu$. We generate $\boldsymbol{X}_1, \ldots, \boldsymbol{X}_n$

Table 1: The estimated coverage probability (the column Prob.) and averaged confidence interval length (the column Length) for $\mathcal{I}_{\text{Conn}}$. We set the dimension $d = 100$, the sample size $n \in \{400, 600\}$, the values of the invariant $k = -25$, $k_\mu \in \{-25, -26, -27, -28\}$ and the signal strength $\mu \in \{0.2, 0.4, 0.6, 0.8\}$. The results are calculated based on 500 repetitions.

|   |   | $\mu = 0.2$ | | $\mu = 0.4$ | | $\mu = 0.6$ | | $\mu = 0.8$ | |
| --- | --- | --- | --- | --- | --- | --- | --- | --- | --- |
| $n$ | $k_\mu$ | Prob. | Length | Prob. | Length | Prob. | Length | Prob. | Length |
| | -25 | 0.978 | 43.81 | 0.848 | 23.86 | 0.802 | 23.26 | 0.920 | 23.39 |
| | -26 | 0.984 | 44.63 | 0.962 | 24.97 | 0.970 | 24.39 | 0.976 | 24.47 |
| 400 | -27 | 0.970 | 45.13 | 0.970 | 25.65 | 0.968 | 25.01 | 0.970 | 25.08 |
| | -28 | 0.976 | 46.25 | 0.976 | 26.74 | 0.976 | 26.09 | 0.976 | 26.25 |
| | -25 | 0.976 | 35.91 | 0.756 | 23.16 | 0.826 | 23.25 | 0.896 | 23.33 |
| | -26 | 0.974 | 36.72 | 0.940 | 24.08 | 0.972 | 24.15 | 0.966 | 24.25 |
| 600 | -27 | 0.976 | 37.69 | 0.970 | 25.13 | 0.976 | 25.20 | 0.976 | 25.25 |
| | -28 | 0.986 | 38.88 | 0.984 | 26.39 | 0.986 | 26.41 | 0.986 | 26.50 |

i.i.d. from $N(0, \widetilde{\mathbf{\Theta}}^{-1})$ and construct $[\widehat{\mathcal{I}}_L, I_U^*]$ from Algorithm 3. For $\mathcal{I}_{\text{Conn}}$, we consider $k = -25$, $k_\mu = -25, -26, -27$ and $-28$ and set $[I_L^*, I_U^*] = [-d, -1]$.

- $\mathcal{I}_{\text{Deg}}$: Let the adjacency matrix of a $k$-star graph be $\mathbf{A}_{\text{star}(k)} \in \mathbb{R}^{(k+1)\times(k+1)}$ such that $[\mathbf{A}_{\text{star}(k)}]_{1t} = 1$ for $1 \leq t \leq k+1$ and $[\mathbf{A}_{\text{star}(k)}]_{st} = 0$ if $s \neq 1$. We construct the precision matrix with signal strength $\mu$ for the graph assembling $\lfloor d/(k+1) \rfloor$ number of $k$-stars as $\mathbf{\Theta}_{\mu, k} = \mathbf{I}_d + \mu \cdot \text{diag}(\mathbf{A}_{\text{star}(k)}, \ldots, \mathbf{A}_{\text{star}(k)}, \mathbf{I}_{d-k\lfloor d/(k+1) \rfloor})$. Given $k$ and $k_\mu$, for each repetition, denote $\widetilde{\mathbf{A}} = \text{diag}(\mathbf{A}_{\text{star}(k-k_\mu)}, \mathbf{I}_{k-k_\mu})$ and we construct the precision matrix as

$$\widetilde{\mathbf{\Theta}} = \mathbf{I}_d + \mu \text{diag}(\mathbf{A}_{\text{star}(k)}, \ldots, \mathbf{A}_{\text{star}(k)}, \mathbf{I}_{d-k\lfloor d/(k+1) \rfloor}) + (\theta - \mu) \text{diag}(\widetilde{\mathbf{A}}, \ldots, \widetilde{\mathbf{A}}, \mathbf{I}_{d-k\lfloor d/(k+1) \rfloor}).$$

  We can see that $\mathcal{I}_{\text{Deg}}(\widetilde{\mathbf{\Theta}}) = k$ and $\mathcal{I}_{\text{Deg}}(E_{\text{Sig}}(\widetilde{\mathbf{\Theta}})) = k_\mu$. We generate $\mathbf{X}_1, \ldots, \mathbf{X}_n$ i.i.d. from $N(0, \widetilde{\mathbf{\Theta}}^{-1})$ and construct $[\widehat{\mathcal{I}}_L, I_U^*]$ from Algorithm 1 for $\mathcal{I} = \mathcal{I}_{\text{Deg}}$. For $\mathcal{I}_{\text{Deg}}$, we consider $k = 5$, $k_\mu = 5, 4, 3$ and $2$ and set $[I_L^*, I_U^*] = [0, 20]$.

- $\mathcal{I}_{\text{Iso}}$: Since $\mathcal{I}_{\text{Iso}}$ is the negative number of isolated nodes, the parameters $k$ and $k_\mu$ are negative. If $d + k$ is even, we construct the precision matrix with $-k$ isolated nodes as $\mathbf{\Theta}_{\mu, k} = \mu \text{diag}(\mathbf{I}_{-k}, \mathbf{A}_{\text{chain}(1)}, \ldots, \mathbf{A}_{\text{chain}(1)})$. The precision matrix represents the graph containing $-k$ isolated nodes and $(d+k)/2$ disconnected single edges. If $d + k$ is even we construct the precision matrix as $\mathbf{\Theta}_{\mu, k} = \mathbf{I}_d + \mu \text{diag}(\mathbf{I}_{-k}, \mathbf{A}_{\text{chain}(1)}, \ldots, \mathbf{A}_{\text{chain}(1)}, \mathbf{A}_{\text{chain}(2)})$. Since $d + k$ is odd, we let the last chain in the graph contain 2 edges. Given $k$ and $k_\mu$ such that $k - k_\mu$ is even, for each repetition, we uniformly sample $(k - k_\mu)/2$ edges from the single edge chain in $E(\mathbf{\Theta}_{\mu, k})$ and change the weights on these edges from $\mu$ to $\theta$. Denote this new precision matrix as $\widetilde{\mathbf{\Theta}}$ and we can see that $\mathcal{I}_{\text{Iso}}(\widetilde{\mathbf{\Theta}}) = -k$ and $\mathcal{I}_{\text{Iso}}(E_{\text{Sig}}(\widetilde{\mathbf{\Theta}})) = -k_\mu$. We generate $\mathbf{X}_1, \ldots, \mathbf{X}_n$ i.i.d. from $N(0, \widetilde{\mathbf{\Theta}}^{-1})$ and construct $[\widehat{\mathcal{I}}_L, I_U^*]$ from Algorithm 4. For $\mathcal{I}_{\text{Iso}}$, we consider $k = -3$, $k_\mu = -3, -5, -7$ and $-9$ and set $[I_L^*, I_U^*] = [-d, 0]$.

Table 2: The estimated coverage probability (the column Prob.) and averaged confidence interval length (the column Length) for $\mathcal{I}_{\text{Deg}}$. We set the dimension $d = 100$, the sample size $n \in \{400, 600\}$, the values of the invariant $k = 5$, $k_\mu \in \{5, 4, 3, 2\}$ and the signal strength $\mu \in \{0.2, 0.4, 0.6, 0.8\}$. The results are calculated based on 500 repetitions.

| | | $\mu = 0.2$ | | $\mu = 0.4$ | | $\mu = 0.6$ | | $\mu = 0.8$ | |
|---|---|---|---|---|---|---|---|---|---|
| $n$ | $k_\mu$ | Prob. | Length | Prob. | Length | Prob. | Length | Prob. | Length |
| 400 | 5 | 0.962 | 15.68 | 0.964 | 15.38 | 0.962 | 15.39 | 0.962 | 15.39 |
| | 4 | 0.962 | 15.68 | 0.964 | 15.38 | 0.962 | 15.39 | 0.962 | 15.39 |
| | 3 | 0.962 | 16.44 | 0.962 | 16.31 | 0.962 | 16.34 | 0.962 | 16.35 |
| | 2 | 0.962 | 17.31 | 0.962 | 17.28 | 0.962 | 17.29 | 0.962 | 17.31 |
| 600 | 5 | 0.978 | 15.63 | 0.978 | 15.63 | 0.978 | 15.65 | 0.978 | 15.65 |
| | 4 | 0.980 | 15.63 | 0.978 | 15.64 | 0.978 | 15.65 | 0.978 | 15.65 |
| | 3 | 0.978 | 16.60 | 0.978 | 16.59 | 0.978 | 16.62 | 0.978 | 16.62 |
| | 2 | 0.978 | 17.57 | 0.978 | 17.55 | 0.978 | 17.56 | 0.978 | 17.59 |

Given the data $\boldsymbol{X}_1, \ldots, \boldsymbol{X}_n$, we estimate the precision matrix by the CLIME estimator

$$\widehat{\boldsymbol{\Theta}}_j = \underset{\boldsymbol{\beta} \in \mathbb{R}^d}{\arg\min} \|\boldsymbol{\beta}\|_1 \quad \text{s.t.} \quad \|\widehat{\boldsymbol{\Sigma}}\boldsymbol{\beta} - \mathbf{e}_j\|_\infty \leq \lambda \quad (7.1)$$

where $\mathbf{e}_j$ is the $j$-th canonical basis in $\mathbb{R}^d$ for any $j = 1, \ldots, d$. The tuning parameter $\lambda$ in (7.1) is chosen through minimizing a 3-fold cross validation

$$\text{CV}(\lambda) = \sum_{k=1}^{3} \|\widehat{\boldsymbol{\Sigma}}^{(k)} \widehat{\boldsymbol{\Theta}}_\lambda^{(-k)} - \mathbf{I}_d\|_{\text{F}}^2, \quad (7.2)$$

where $\widehat{\boldsymbol{\Sigma}}^{(k)}$ is the sample covariance matrix only using the $k$-th fold of the dataset and $\widehat{\boldsymbol{\Theta}}_\lambda^{(-k)}$ is the CLIME estimator using the remaining data. In the simulations for all three invariants, we set the dimension $d = 100$ and sample size $n = 400$ and $600$. We set $\mu = 0.2, 0.4, 0.6$ and $0.8$. We choose the significance level of confidence intervals as 5%.

The estimated coverage probability and the averaged confidence interval length calculated through 500 repetitions are reported in Tables 1, 2 and 3 for $\mathcal{I}_{\text{Conn}}$, $\mathcal{I}_{\text{Deg}}$ and $\mathcal{I}_{\text{Iso}}$ respectively. From these results, we can see that when the value $\mu$ is relatively small, the confidence interval lengths are larger in order to guarantee the confidence interval cover the true invariant under small signal strength. When $\mu$ becomes larger, the confidence interval lengths converge to the optimal rate $O(I_U^* - \mathcal{I}(E_{\text{Sig}}(\boldsymbol{\Theta})))$ shown in (4.10). This illustrates that the proposed confidence interval is adaptive to $k_\mu$.

Table 3: The estimated coverage probability (the column Prob.) and averaged confidence interval length (the column Length) for $\mathcal{I}_{\text{Iso}}$. We set the dimension $d = 100$, the sample size $n \in \{400, 600\}$, the values of the invariant $k = -3$, $k_\mu \in \{-3, -5, -7, -9\}$ and the signal strength $\mu \in \{0.2, 0.4, 0.6, 0.8\}$. The results are calculated based on 500 repetitions.

|       |       | $\mu = 0.2$ |        | $\mu = 0.4$ |        | $\mu = 0.6$ |        | $\mu = 0.8$ |        |
|-------|-------|-------|--------|-------|--------|-------|--------|-------|--------|
| $n$   | $k_\mu$ | Prob. | Length | Prob. | Length | Prob. | Length | Prob. | Length |
|       | -3    | 0.962 | 32.97  | 0.964 | 4.32   | 0.958 | 2.84   | 0.938 | 2.82   |
| 400   | -5    | 0.972 | 33.49  | 0.972 | 6.25   | 0.972 | 4.78   | 0.972 | 4.77   |
|       | -7    | 0.970 | 51.42  | 0.970 | 8.28   | 0.970 | 6.70   | 0.970 | 6.66   |
|       | -9    | 0.964 | 36.63  | 0.964 | 10.15  | 0.964 | 8.58   | 0.964 | 8.51   |
|       | -3    | 0.978 | 21.16  | 0.976 | 2.93   | 0.972 | 2.91   | 0.962 | 2.89   |
| 600   | -5    | 0.974 | 22.63  | 0.974 | 4.87   | 0.974 | 4.83   | 0.974 | 4.81   |
|       | -7    | 0.968 | 24.15  | 0.968 | 6.74   | 0.968 | 6.71   | 0.968 | 6.69   |
|       | -9    | 0.974 | 26.19  | 0.974 | 8.70   | 0.974 | 8.68   | 0.974 | 8.66   |

## 7.2 Neuroscience Application

We apply our inferential method to the brain imaging dataset studied in Simony et al. (2016). The dataset contains fMRI scans from 36 subjects taken while the subjects were listening to the stimuli generated from a seven-minute story *Pieman* (told by Jim O'Grady at the "The Moth" storytelling event). The 36 subjects also listened to a word-scrambled version of the story. In particular, the story was segmented into 608 short words and their order was scrambled randomly. The raw functional data was preprocessed to correct head motion, time slicing, spatial smoothing and temporal filtering in Simony et al. (2016). For both the intact story and the word scrambled settings, each subject had 300 fMRI measurements and the measurements were taken every 1.4 seconds.

The original fMRI dataset has the dimensional 271,633 representing 271,633 3-mm isotropic voxels. We reduce the dimension to 172 regions of interest (ROIs) introduced by Baldassano et al. (2015) through averaging the voxels in the same ROI. Therefore, for each subject, we have the data with dimension $d = 172$ and sample size $n = 300$. We average the data across the 36 subject to obtain a single $300 \times 172$ dataset and standardize each ROI such that they have mean zero and standard deviation one. We apply the Gaussian graphical model to the dataset so that the brain network is induced by the precision matrix $\boldsymbol{\Theta}$ and each ROI corresponds to a node in the network. Our goal is to infer two combinatorial quantities of the brain network: the number of connected subgraphs and the maximum degree. In fact, to grasp more detailed structural information of the network, we aim to infer the above two invariants for the precision matrix at different filtration levels. In specific, given a precision matrix $\boldsymbol{\Theta}$ and a filtration level $\mu > 0$, we define the thresholded matrix $[\mathcal{T}_\mu(\boldsymbol{\Theta})]_{jk} = \boldsymbol{\Theta}_{jk} \mathbb{1}\{|\boldsymbol{\Theta}_{jk}| \geq \mu\}$ for all $1 \leq j, k \leq d$. We want to construct confidence intervals for invariants $\mathcal{I}_{\text{Conn}}(\mathcal{T}_\mu(\boldsymbol{\Theta}))$ and $\mathcal{I}_{\text{Deg}}(\mathcal{T}_\mu(\boldsymbol{\Theta}))$ for different levels of $\mu$ under both the intact story and word scrambled settings. It is easy to check that $\mathcal{I}_{\text{Conn}}(\mathcal{T}_\mu(\boldsymbol{\Theta}))$ and $\mathcal{I}_{\text{Deg}}(\mathcal{T}_\mu(\boldsymbol{\Theta}))$

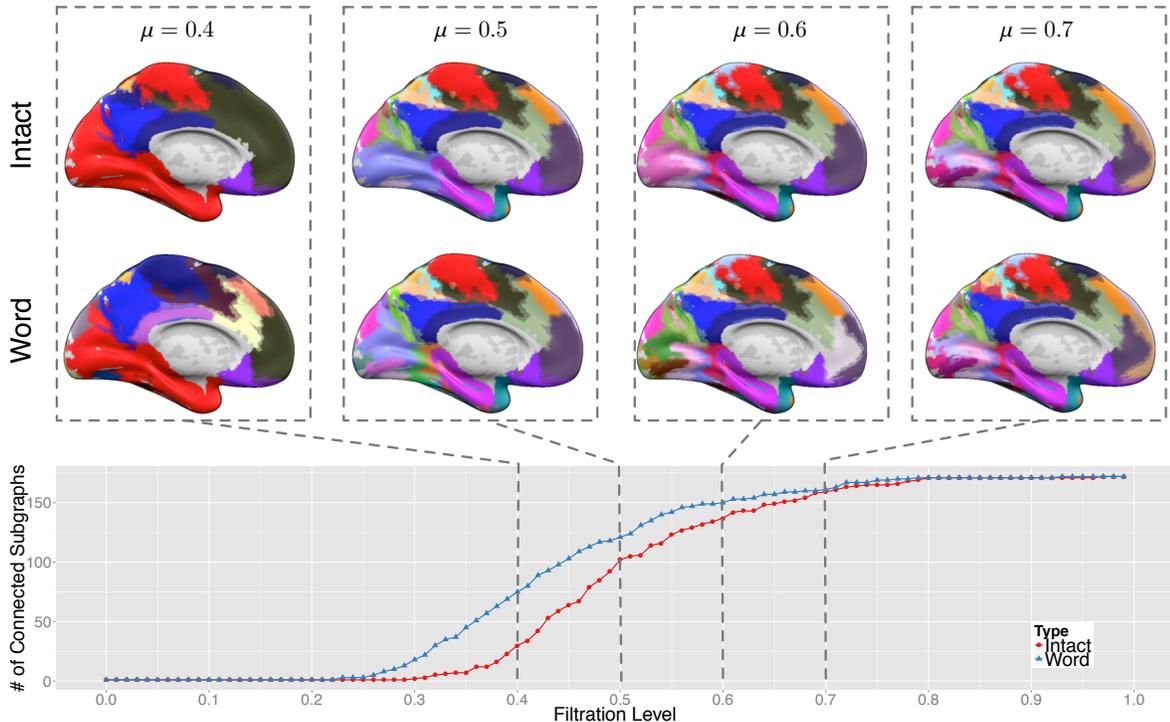

Figure 5: The connected subgraphs in brain networks changing with the filtration level $\mu$. The upper panel illustrates the connected regions of interest by same colors. The lower panel shows how the number of connected subgraphs change with the filtration level $\mu$ for both intact story (in red) and word scrambled (in blue).

are also monotone invariants, therefore the skip-down algorithm can be applied.

As before, the tuning parameter of the CLIME estimator is chosen via 3-fold cross validation using the risk in (7.2). The significance level of the confidence interval is set to $\alpha = 5\%$ and since there is no prior information on how large the maximum degree of the network is, or how many connected subgraphs the network contains, we consider the largest possible range for both invariants.

We implement Algorithm 1 to construct the confidence intervals. Figures 5 and 6 visualize how the lower endpoints of the confidence intervals, of the number of connected subgraphs and the maximum degree respectively, change with the filtration level $\mu$. In addition, we visualize the structural information of the output edge set $E_{t^*}$ generated from the skip-down algorithm, where $t^*$ is the number of iterations needed for Algorithm 1 to conclude. Figures 5 illustrates the different connected ROIs in $G_{t^*} = (V, E_{t^*})$ by different colors and Figure 6 illustrates the degree of each ROI in $G_{t^*}$. For the number of connected subgraphs, Figure 5 shows that the brain network has fewer connected subgraphs when the subject is listening to the intact story compared to when the subject is listening to the scrambled word version. In particular, we can see in Figure 7(a) that compared to the word scrambled setting, the dorsolateral prefrontal cortex (DL-PFC) is connected to the inferior frontal gyrus (IFG) under the intact story setting. The inferior frontal gyrus is the brain area responsible for language processing and comprehension (Grewe et al., 2005; Caplan, 2006) and the dorsolateral prefrontal cortex is responsible for working memory tasks (Barbey et al., 2013). The fact that these two areas are connected in the intact story setting suggests that both language processing and memory are working together in the procedure of understanding the intact Pieman story.

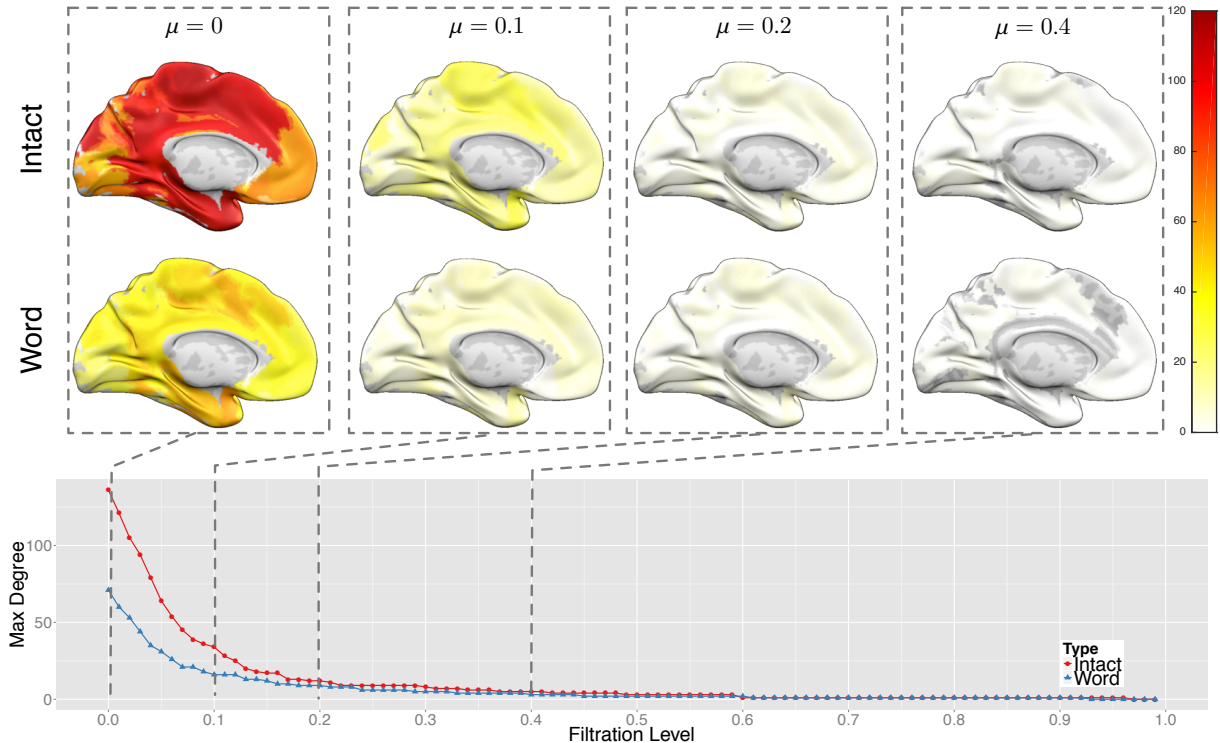

Figure 6: The degree of region of interest in brain networks changing with filtration level $\mu$. The upper panel illustrates the degree of each region of interest. The lower panel shows how the maximum degree changes with the filtration level $\mu$ for both intact story (in red) and word scrambled (in blue).

Regarding maximum degree, Figure 6 shows that the brain network has higher maximum degree when the subject is listening to the intact story compared to when the subject is listening to scrambled words. In Figure 7(b) we also observe that the precuneus area, which is known to be involved with understanding high-level concepts in stories (Lerner et al., 2011; Ames et al., 2015), has a higher degree under the intact story setting compared to the word scrambled setting.

# 8 Discussion

In this manuscript we propose generic inferential methods for monotone invariants and monotone properties under graphical models. Our skip-down algorithm is based on screening critical edge sets specific to the invariant of interest. We provide fast algorithms critical edge set search for several graph invariants. An interesting research direction is to provide an efficient critical edge set searching algorithm for a general family of invariants. Moreover, we can explore the inferential methods for the combinatorial structures beyond graphs in the future researches, including ranking, hypergraph and partitions.

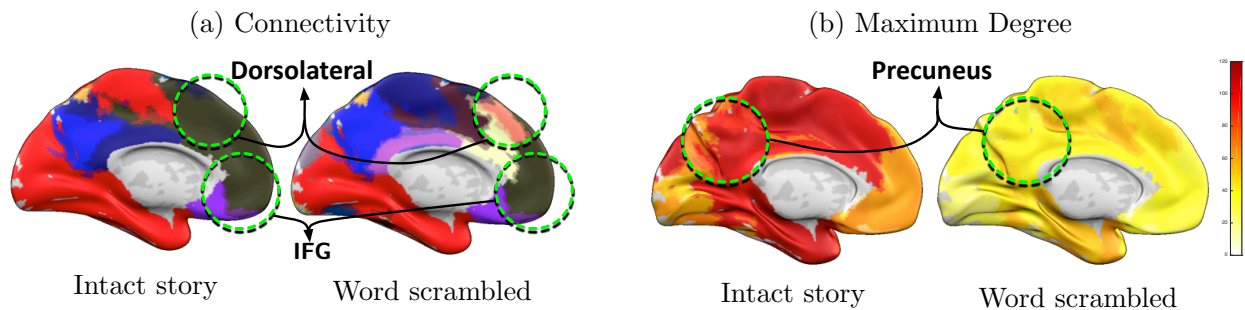

Figure 7: Structural differences in cerebral cortices between the intact story and word scrambled settings. (a) shows the connected regions of interest in brain networks under the filtration level $\mu = 0.4$. (b) shows the degree of regions of interest in brain networks under the filtration level $\mu = 0$.

For the lower bound on the confidence interval length, we show that the signal strength at the rate $\sqrt{\log d / n}$ is critical for a large family of invariants. A standing question is to find out whether invariants whose confidence interval length has lower bounds depending on smaller signal strength exist, and to construct optimal confidence intervals for such invariants.

# References

ADDARIO-BERRY, L., BROUTIN, N., DEVROYE, L. and LUGOSI, G. (2010). On combinatorial testing problems. *The Annals of Statistics* **38** 3063–3092.

ALLEN, G. I. and LIU, Z. (2012). A log-linear graphical model for inferring genetic networks from high-throughput sequencing data. In *Bioinformatics and Biomedicine (BIBM), 2012 IEEE International Conference on*. IEEE.

AMES, D. L., HONEY, C. J., CHOW, M. A., TODOROV, A. and HASSON, U. (2015). Contextual alignment of cognitive and neural dynamics. *Journal of Cognitive Neuroscience* **27** 655–664.

ANANDKUMAR, A., TAN, V. Y., HUANG, F., WILLSKY, A. S. ET AL. (2012). High-dimensional structure estimation in Ising models: Local separation criterion. *The Annals of Statistics* **40** 1346–1375.

ARIAS-CASTRO, E., BUBECK, S. and LUGOSI, G. (2012). Detection of correlations. *The Annals of Statistics* **40** 412–435.

ARIAS-CASTRO, E., BUBECK, S., LUGOSI, G. and VERZELEN, N. (2015). Detecting Markov random fields hidden in white noise. *arXiv preprint arXiv:1504.06984* .

ARIAS-CASTRO, E., CANDÈS, E. J. and DURAND, A. (2011). Detection of an anomalous cluster in a network. *The Annals of Statistics* **39** 278–304.

Baldassano, C., Beck, D. M. and Fei-Fei, L. (2015). Parcellating connectivity in spatial maps. *PeerJ* **3** e784.

Barbey, A. K., Koenigs, M. and Grafman, J. (2013). Dorsolateral prefrontal contributions to human working memory. *Cortex* **49** 1195–1205.

Bauer, P. and Hackl, P. (1987). Multiple testing in a set of nested hypotheses. *Statistics: A Journal of Theoretical and Applied Statistics* **18** 345–349.

Cai, T., Liu, W. and Luo, X. (2011). A constrained $\ell_1$ minimization approach to sparse precision matrix estimation. *Journal of the American Statistical Association* **106** 594–607.

Cai, T. T. and Guo, Z. (2015). Confidence intervals for high-dimensional linear regression: Minimax rates and adaptivity. *arXiv preprint arXiv:1506.05539* .

Cai, T. T., Liu, W. and Xia, Y. (2013). Two-sample covariance matrix testing and support recovery in high-dimensional and sparse settings. *Journal of the American Statistical Association* **108** 265–277.

Cai, T. T., Liu, W. and Xia, Y. (2014). Two-sample test of high dimensional means under dependence. *Journal of the Royal Statistical Society. Series B. Statistical Methodology* **76** 349–372.

Caplan, D. (2006). Why is broca's area involved in syntax? *Cortex* **42** 469–471.

Chen, S., Witten, D. M. and Shojaie, A. (2015). Selection and estimation for mixed graphical models. *Biometrika* **102** 47–64.

Chernozhukov, V., Chetverikov, D. and Kato, K. (2013). Gaussian approximations and multiplier bootstrap for maxima of sums of high-dimensional random vectors. *The Annals of Statistics* **41** 2786–2819.

Danaher, P., Wang, P. and Witten, D. M. (2014). The joint graphical Lasso for inverse covariance estimation across multiple classes. *Journal of the Royal Statistical Society. Series B. Statistical Methodology* **76** 373–397.

Fan, J., Liu, H., Ning, Y. and Zou, H. (2014). High dimensional semiparametric latent graphical model for mixed data. *ArXiv e-prints, arXiv:1404.7236* .

Fellinghauer, B., Bühlmann, P., Ryffel, M., von Rhein, M. and Reinhardt, J. D. (2013). Stable graphical model estimation with random forests for discrete, continuous, and mixed variables. *Comput. Stat. Data Anal.* **64** 132–152.

Friedman, J. H., Hastie, T. J. and Tibshirani, R. J. (2008). Sparse inverse covariance estimation with the graphical Lasso. *Biostatistics* **9** 432–441.

Grewe, T., Bornkessel, I., Zysset, S., Wiese, R., von Cramon, D. Y. and Schlesewsky, M. (2005). The emergence of the unmarked: A new perspective on the language-specific function of broca's area. *Human brain mapping* **26** 178–190.

Gu, Q., Cao, Y., Ning, Y. and Liu, H. (2015). Local and global inference for high dimensional gaussian copula graphical models. *arXiv preprint arXiv:1502.02347* .

Hagmann, P., Cammoun, L., Gigandet, X., Meuli, R., Honey, C. J., Wedeen, V. J. and Sporns, O. (2008). Mapping the structural core of human cerebral cortex. *PLoS Biol* **6** e159.

Han, F. and Liu, H. (2014). Distribution-free tests of independence with applications to testing more structures. *arXiv preprint arXiv:1410.4179* .

Höfling, H. and Tibshirani, R. J. (2009). Estimation of sparse binary pairwise Markov networks using pseudo-likelihoods. *J. Mach. Learn. Res.* **10** 883–906.

Jalali, A., Ravikumar, P. D., Vasuki, V. and Sanghavi, S. (2011). On learning discrete graphical models using group-sparse regularization. In *International Conference on Artificial*

*Intelligence and Statistics*.

JANKOVÁ, J. and VAN DE GEER, S. (2016). Honest confidence regions and optimality in high-dimensional precision matrix estimation. *TEST* 1–20.

JAVANMARD, A. and MONTANARI, A. (2014). Confidence intervals and hypothesis testing for high-dimensional regression. *J. Mach. Learn. Res.* **15** 2869–2909.

JOAG-DEV, K. and PROSCHAN, F. (1983). Negative association of random variables with applications. *The Annals of Statistics* 286–295.

KRUSKAL, J. B. (1956). On the shortest spanning subtree of a graph and the traveling salesman problem. *Proceedings of the American Mathematical society* **7** 48–50.

LAM, C. and FAN, J. (2009). Sparsistency and rates of convergence in large covariance matrix estimation. *The Annals of Statistics* **37** 4254–4278.

LECAM, L. (1973). Convergence of estimates under dimensionality restrictions. *The Annals of Statistics* **1** 38–53.

LEE, J. D. and HASTIE, T. J. (2015). Learning the structure of mixed graphical models. *Journal of Computational and Graphical Statistics* **25** 230–253.

LEE, S.-I., GANAPATHI, V. and KOLLER, D. (2006). Efficient structure learning of Markov networks using $\ell_1$-regularization. In *Advances in Neural Information Processing Systems*.

LERNER, Y., HONEY, C. J., SILBERT, L. J. and HASSON, U. (2011). Topographic mapping of a hierarchy of temporal receptive windows using a narrated story. *The Journal of Neuroscience* **31** 2906–2915.

LIU, H., HAN, F., YUAN, M., LAFFERTY, J. D. and WASSERMAN, L. A. (2012a). High-dimensional semiparametric Gaussian copula graphical models. *The Annals of Statistics* **40** 2293–2326.

LIU, H., HAN, F. and ZHANG, C.-H. (2012b). Transelliptical graphical models. In *Proc. of NIPS* (P. Bartlett, F. Pereira, C. Burges, L. Bottou and K. Weinberger, eds.). 809–817.

LIU, H. and WASSERMAN, J. D. L. L. A. (2009). The nonparanormal: Semiparametric estimation of high dimensional undirected graphs. *J. Mach. Learn. Res.* **10** 2295–2328.

LIU, W. (2013). Gaussian graphical model estimation with false discovery rate control. *The Annals of Statistics* **41** 2948–2978.

LU, J., KOLAR, M. and LIU, H. (2015). Post-regularization inference for dynamic nonparanormal graphical models. *arXiv preprint arXiv:1512.08298* .

LUSCOMBE, N. M., BABU, M. M., YU, H., SNYDER, M., TEICHMANN, S. A. and GERSTEIN, M. (2004). Genomic analysis of regulatory network dynamics reveals large topological changes. *Nature* **431** 308–312.

MEINSHAUSEN, N. and BÜHLMANN, P. (2006). High dimensional graphs and variable selection with the Lasso. *The Annals of Statistics* **34** 1436–1462.

MOHAN, K., LONDON, P., FAZEL, M., WITTEN, D. M. and LEE, S.-I. (2014). Node-based learning of multiple gaussian graphical models. *J. Mach. Learn. Res.* **15** 445–488.

NASH-WILLIAMS, C. S. J. (1961). Edge-disjoint spanning trees of finite graphs. *Journal of the London Mathematical Society* **1** 445–450.

NASH-WILLIAMS, C. S. J. (1964). Decomposition of finite graphs into forests. *Journal of the London Mathematical Society* **1** 12–12.

NEYKOV, M., LU, J. and LIU, H. (2016). Combinatorial inference for graphical models. *arXiv preprint arXiv:1608.03045* .

NEYKOV, M., NING, Y., LIU, J. S. and LIU, H. (2015). A unified theory of confidence regions and testing for high dimensional estimating equations. *arXiv preprint arXiv:1510.08986* .

NING, Y. and LIU, H. (2014). A general theory of hypothesis tests and confidence regions for sparse high dimensional models. *arXiv preprint arXiv:1412.8765* .

PENG, J., WANG, P., ZHOU, N. and ZHU, J. (2009). Partial correlation estimation by joint sparse regression models. *J. Am. Stat. Assoc.* **104** 735–746.

QIU, H., HAN, F., LIU, H. and CAFFO, B. S. (2013). Joint estimation of multiple graphical models from high dimensional time series. *ArXiv e-prints, arXiv:1311.0219* .

RAVIKUMAR, P., WAINWRIGHT, M. J. and LAFFERTY, J. D. (2010). High-dimensional Ising model selection using $\ell_1$-regularized logistic regression. *The Annals of Statistics* **38** 1287–1319.

RAVIKUMAR, P., WAINWRIGHT, M. J., RASKUTTI, G. and YU, B. (2011). High-dimensional covariance estimation by minimizing $\ell_1$-penalized log-determinant divergence. *Electron. J. Stat.* **5** 935–980.

REN, Z., SUN, T., ZHANG, C.-H. and ZHOU, H. H. (2015). Asymptotic normality and optimalities in estimation of large gaussian graphical models. *Ann. Statist.* **43** 991–1026.

ROMANO, J. and WOLF, M. (2005). Exact and approximate stepdown methods for multiple hypothesis testing. *Journal of the American Statistical Association* **100** 94–108.

RUBINOV, M. and SPORNS, O. (2010). Complex network measures of brain connectivity: uses and interpretations. *Neuroimage* **52** 1059–1069.

SHEN, X., PAN, W. and ZHU, Y. (2012). Likelihood-based selection and sharp parameter estimation. *J. Am. Stat. Assoc.* **107** 223–232.

SIMONY, E., HONEY, C. J., CHEN, J., LOSITSKY, O., YESHURUN, Y., WIESEL, A. and HASSON, U. (2016). Dynamic reconfiguration of the default mode network during narrative comprehension. *Nature Communications* **7** 12141.

VAN DE GEER, S. A., BÜHLMANN, P., RITOV, Y. and DEZEURE, R. (2014). On asymptotically optimal confidence regions and tests for high-dimensional models. *The Annals of Statistics* **42** 1166–1202.

VAN DER VAART, A. W. and WELLNER, J. A. (1996). *Weak Convergence and Empirical Processes: With Applications to Statistics*. Springer.

XUE, L. and ZOU, H. (2012). Regularized rank-based estimation of high-dimensional nonparanormal graphical models. *The Annals of Statistics* **40** 2541–2571.

YANG, E., ALLEN, G. I., LIU, Z. and RAVIKUMAR, P. (2012). Graphical models via generalized linear models. In *Advances in Neural Information Processing Systems 25*. 1358–1366.

YANG, E., RAVIKUMAR, P., ALLEN, G. I., BAKER, Y., WAN, Y.-W. and LIU, Z. (2014a). A general framework for mixed graphical models. *arXiv preprint arXiv:1411.0288* .

YANG, E., RAVIKUMAR, P., ALLEN, G. I. and LIU, Z. (2013). On poisson graphical models. In *Advances in Neural Information Processing Systems 26*. 1718–1726.

YANG, Z., NING, Y. and LIU, H. (2014b). On semiparametric exponential family graphical models. *arXiv preprint arXiv:1412.8697* .

YUAN, M. and LIN, Y. (2007). Model selection and estimation in the gaussian graphical model. *Biometrika* **94** 19–35.

ZHANG, C.-H. and ZHANG, S. S. (2013). Confidence intervals for low dimensional parameters in high dimensional linear models. *J. R. Stat. Soc. B* **76** 217–242.

ZHOU, S., VAN DE GEER, S. and BÜHLMANN, P. (2009). Adaptive Lasso for high dimensional regression and Gaussian graphical modeling. *arXiv preprint arXiv:0903.2515* .

*Technical supplementary material to*

# Adaptive Inferential Method for Monotone Graph Invariants

Junwei Lu, Matey Neykov and Han Liu

This document contains the supplementary material to the paper "Adaptive Inferential Method for Monotone Graph Invariants". We mainly provide technical details of proving confidence interval lower bound here.

## S.1 Proofs of Results on the Skip-Down Method

In this section, we prove theoretical results on skip-down algorithm. We will prove Theorems 4.1, 4.5 and 4.6 and Proposition 3.2.

### S.1.1 Proof of Theorem 4.1

Under (3.5), it suffices to prove (4.1) and (4.2) by showing that

$$\liminf_{n\to\infty} \inf_{\boldsymbol{\Theta}\in\mathcal{U}_s(I_L^*,I_U^*)} \mathbb{P}_{\boldsymbol{\Theta}}\big(\mathcal{I}(\boldsymbol{\Theta}) \leq \widehat{\mathcal{I}}_L\big) \geq 1-\alpha. \tag{S.1}$$

We denote the true edge set as $E^* = E(\boldsymbol{\Theta})$. We aim to bound the probability of the event $\{\mathcal{I}(\boldsymbol{\Theta}) < \widehat{\mathcal{I}}_L\}$. Assume that the output of Algorithm 1 satisfies $\mathcal{I}(\boldsymbol{\Theta}) < \widehat{\mathcal{I}}_L$ and that the algorithm has stopped at the $t = K$-th step. Under these assumptions the final rejected edge set $E_K$ satisfies $\mathcal{I}(E_K) > \mathcal{I}(\boldsymbol{\Theta})$. Our high-level outline to prove this theorem is to show the following two facts:

(1) At least one rejected edge belongs to $\mathcal{C}_{\mathcal{I}}(E^*)$;

(2) If the $\ell$-th iteration is the step before which no edge in $\mathcal{C}_{\mathcal{I}}(E^*)$ is rejected, then $\mathcal{C}_{\mathcal{I}}(E^*) \subseteq \mathcal{C}_{\mathcal{I}}(E_{\ell-1})$.

Suppose the above two claims have been proved. We denote the first edge rejected in $\mathcal{C}_{\mathcal{I}}(E^*)$ as $\bar{e}$. According to the steps in Algorithm 1, we have

$$\max_{e\in\mathcal{C}_{\mathcal{I}}(E^*)} \sqrt{n}\widehat{\boldsymbol{\Theta}}_e^d \geq \sqrt{n}\widehat{\boldsymbol{\Theta}}_{\bar{e}}^d \geq c(\alpha, \mathcal{C}_{\mathcal{I}}(E_{\ell-1})) \geq c(\alpha, \mathcal{C}_{\mathcal{I}}(E^*)), \tag{S.2}$$

where the first inequality is by $\bar{e} \in \mathcal{C}_{\mathcal{I}}(E^*)$, the second inequality is by the mechanism of updating the rejected set and the third inequality is by $\mathcal{C}_{\mathcal{I}}(E^*) \subseteq \mathcal{C}_{\mathcal{I}}(E_{\ell-1})$. Therefore, by (3.5) and the fact that $\boldsymbol{\Theta}_e = 0$ for any $e \in \mathcal{C}_{\mathcal{I}}(E^*) \subseteq (E^*)^c$, we have

$$\sup_{\boldsymbol{\Theta}\in\mathcal{U}_s(I_L^*,I_U^*)} \mathbb{P}_{\boldsymbol{\Theta}}(\mathcal{I}(\boldsymbol{\Theta}) < \widehat{\mathcal{I}}_L) \leq \sup_{\boldsymbol{\Theta}\in\mathcal{U}_s(I_L^*,I_U^*)} \mathbb{P}_{\boldsymbol{\Theta}}\Big(\max_{e\in\mathcal{C}_{\mathcal{I}}(E^*)} \sqrt{n}\widehat{\boldsymbol{\Theta}}_e^d \geq c(\alpha, \mathcal{C}_{\mathcal{I}}(E^*))\Big) \leq \alpha + o(1).$$

So in the rest part of the proof, we only need to show the above two facts.

We first prove that at least one rejected edge belongs to $\mathcal{C}_\mathcal{I}(E^*)$. To begin with, there must exist at least one edge in $(E^*)^c$ which is rejected, otherwise $E_K \subseteq E^*$ which implies a contradictory result that $\widehat{\mathcal{I}}_L = \mathcal{I}(E_K) \leq \mathcal{I}(E^*) = \mathcal{I}(\Theta)$. To show the stronger result that $E_K \cap \mathcal{C}_\mathcal{I}(E^*) \neq \varnothing$, we consider the edge set $E_K \cup E^*$. Since $\mathcal{I}(E_K) > \mathcal{I}(E^*)$ and $\mathcal{I}$ is monotone, we have $\mathcal{I}(E_K \cup E^*) > \mathcal{I}(E^*)$ as well. As we have shown that $E_K \backslash E^* \neq \varnothing$, we denote the edges in this set as $E_K \backslash E^* = \{e_1, \ldots, e_m\}$. Here the order of the edges' indices is arbitrary. Consider the nested sequence

$$E^* \subseteq E^* \cup \{e_1\} \subseteq E^* \cup \{e_1, e_2\} \subseteq \cdots \subseteq E^* \cup \{e_1, \ldots, e_m\} = E_K \cup E^*,$$

where $\mathcal{I}(E_K \cup E^*) > \mathcal{I}(E^*)$.

Notice that the invariant of the first edge set in the above sequence is strictly smaller than the one of the last one. The fact that $\mathcal{I}$ is monotone implies that there exists an integer $m_0 \in [0, m]$ such that $\mathcal{I}(E^* \cup \{e_1, \ldots, e_{m_0}\}) < \mathcal{I}(E^* \cup \{e_1, \ldots, e_{m_0+1}\})$, where for $m_0 = 0$, we denote $E^* \cup \{e_1, \ldots, e_{m_0}\} = E^*$. In particular, we have the sequence

$$\underbrace{E^* \subseteq \cdots \subseteq E^* \cup \{e_1, \ldots, e_{m_0}\}}_{\mathcal{I}(\bullet) \leq \mathcal{I}(E^* \cup \{e_1, \ldots, e_{m_0}\})} \subseteq \underbrace{E^* \cup \{e_1, \ldots, e_{m_0+1}\} \subseteq \cdots \subseteq E_K \cup E^*}_{\mathcal{I}(E^* \cup \{e_1, \ldots, e_{m_0+1}\}) \leq \mathcal{I}(\bullet)}. \tag{S.3}$$

Comparing to (3.3), we have $e_{m_0+1} \in \mathcal{C}_\mathcal{I}(E^*)$ and by the construction, $e_{m_0+1} \in E_K$. Therefore, we prove that $E_K \cap \mathcal{C}_\mathcal{I}(E^*) \neq \varnothing$.

Recall $\ell$ denotes the first iteration of Algorithm 1 in which and edge from $\mathcal{C}_\mathcal{I}(E^*)$ has been rejected. Now, we are going to show that $\mathcal{C}_\mathcal{I}(E^*) \subseteq \mathcal{C}_\mathcal{I}(E_{\ell-1})$. We have $E_{\ell-1} \cap \mathcal{C}_\mathcal{I}(E^*) = \varnothing$ and $\mathcal{I}(E_{\ell-1}) \leq \mathcal{I}(E^*)$. For any $e' \in \mathcal{C}_\mathcal{I}(E^*)$, denote $E'$ is the set satisfying (3.3) for $\mathcal{C}_\mathcal{I}(E^*)$ such that $E' \supseteq E^*$ and $\mathcal{I}(E') > \mathcal{I}(E' \backslash \{e'\})$. The monotone invariant implies that $\mathcal{I}(E' \cup E_{\ell-1}) \geq \mathcal{I}(E')$. To prove $e' \in \mathcal{C}_\mathcal{I}(E_{\ell-1})$, it suffices to show that

$$\mathcal{I}(E' \cup E_{\ell-1}) > \mathcal{I}\big((E' \cup E_{\ell-1}) \backslash \{e'\}\big).$$

Combining $E_{\ell-1} \cap \mathcal{C}_\mathcal{I}(E^*) = \varnothing$ with $e' \in \mathcal{C}_\mathcal{I}(E^*)$, it is equivalent to show

$$\mathcal{I}(E' \cup E_{\ell-1}) > \mathcal{I}((E' \backslash \{e'\}) \cup (E_{\ell-1} \backslash E')).$$

This is obviously true if $E_{\ell-1} \backslash E' = \varnothing$ as $\mathcal{I}(E') > \mathcal{I}(E' \backslash \{e'\})$. For the case $E_{\ell-1} \backslash E' \neq \varnothing$, we prove by contradiction. If $\mathcal{I}((E' \backslash \{e'\}) \cup (E_{\ell-1} \backslash E')) \geq \mathcal{I}(E' \cup E_{\ell-1})$, we have

$$E' \backslash \{e'\} \subseteq (E' \backslash \{e'\}) \cup E_{\ell-1} \text{ and } \mathcal{I}(E' \backslash \{e'\}) < \mathcal{I}(E' \cup E_{\ell-1}) \leq \mathcal{I}((E' \backslash \{e'\}) \cup E_{\ell-1}).$$

Similar to (S.3), we denote $E_{\ell-1} \backslash E' = \{e'_1, \ldots, e'_{m'}\}$ and there exists an integer $m'_0 \in [0, m']$ such that

$$\underbrace{(E' \backslash \{e'\}) \cup \{e'_1, \ldots, e'_{m'_0}\}}_{\mathcal{I}(\bullet) \geq \mathcal{I}(E' \backslash \{e'\})} \subseteq \underbrace{(E' \backslash \{e'\}) \cup \{e'_1, \ldots, e'_{m'_0+1}\}}_{\mathcal{I}((E' \backslash \{e'\}) \cup E_{\ell-1}) \leq \mathcal{I}(\bullet)}. \tag{S.4}$$

Since $e'_{m'_0+1} \in E_{\ell-1}\backslash E'$ and $E^* \subseteq E'\backslash\{e'\}$, we have $e'_{m'_0+1} \notin E^*$. Combining with (S.4), the definition in (3.3) is satisfied for $e'_{m'_0+1}$ and thus $e'_{m'_0+1} \in \mathcal{C}_\mathcal{I}(E^*)$. This contradicts the fact that $E_{\ell-1} \cap \mathcal{C}_\mathcal{I}(E^*) = \emptyset$. In summary, we show that $\mathcal{I}(E' \cup E_{\ell-1}) > \mathcal{I}((E' \cup E_{\ell-1})\backslash\{e'\})$ and therefore, due to the previous discussion, prove that $\mathcal{C}_\mathcal{I}(E^*) \subseteq \mathcal{C}_\mathcal{I}(E_{\ell-1})$.

## S.1.2 Proof of Theorems 4.5 and 4.6

We first define a few notations before presenting the proof. For a positive number $\mu > 0$, recall that the thresholded matrix $[\mathcal{T}_\mu(\Theta)]_{jk} = \Theta_{jk}\mathbb{1}\{|\Theta_{jk}| \geq \mu\}$ for all $1 \leq j, k \leq d$. We define $\mathcal{T}_\mu(\Theta)$ in order to show the adaptivity of our confidence interval. We can see $\mathcal{I}(\mathcal{T}_\mu(\Theta)) = \mathcal{I}(E_{\text{Sig}}(\Theta))$ when $\mu = C\sqrt{\log d/n}$.

We also define the parameter space

$$\mathcal{U}_s(I_L^*, I_U^*; \theta, \mu) = \left\{\Theta \in \mathcal{U}_s \;\middle|\; \mathcal{I}(\Theta) \leq I_U^*, \min_{e \in E(\Theta)} |\Theta_e| \geq \theta, \mathcal{I}(\mathcal{T}_\mu(\Theta)) \geq I_L^*\right\}. \tag{S.5}$$

The proof of both Theorems 4.5 and S.7 can be directly derived from the following lemma.

**Lemma S.1.** *Suppose $\Theta \in \mathcal{U}_s$ and $(\log(dn))^6/n + s^2(\log dn)^4/n = o(1)$. For any monotone invariant $\mathcal{I}$ with range $[I_L^*, I_U^*]$, there exists a positive constant $C$ such that if $\mu \geq C\sqrt{\log d/n}$, for any $\alpha \in (0,1)$ and $\theta > 0$,*

$$\sup_{\Theta \in \mathcal{U}_\mathcal{I}(I_L^*, I_U^*; \theta)} \mathbb{P}_\Theta(\mathcal{I}(\mathcal{T}_\mu(\Theta)) \leq \widehat{\mathcal{I}}_L) = 1 - O(1/d^3). \tag{S.6}$$

We defer the proof of Lemma S.1 to Section S.4 in the technical supplementary. We first give the proof of Theorem 4.5. By the construction of $[\widehat{\mathcal{I}}_L, I_U^*]$ in Algorithm 1, we have

$$\mathbb{E}_\Theta[I_U^* - \widehat{\mathcal{I}}_L] \leq \mathbb{E}_\Theta[I_U^* - \mathcal{I}(\mathcal{T}_\mu(\Theta)) \vee I_L^* \,|\, \mathcal{I}(\mathcal{T}_\mu(\Theta)) \leq \widehat{\mathcal{I}}_L]\mathbb{P}(\mathcal{I}(\mathcal{T}_\mu(\Theta)) \leq \widehat{\mathcal{I}}_L)$$
$$+ (I_U^* - I_L^*)\mathbb{P}(\mathcal{I}(\mathcal{T}_\mu(\Theta)) > \widehat{\mathcal{I}}_L)$$
$$\leq [I_U^* - \mathcal{I}(\mathcal{T}_\mu(\Theta)) \vee I_L^* + 1]\mathbb{P}(\mathcal{I}(\mathcal{T}_\mu(\Theta)) \leq \widehat{\mathcal{I}}_L)$$
$$+ (I_U^* - I_L^*)\mathbb{P}(\mathcal{I}(\mathcal{T}_\mu(\Theta)) > \widehat{\mathcal{I}}_L).$$

Therefore, by (S.6), we have

$$\sup_{\Theta \in \mathcal{U}_\mathcal{I}(I_L^*, I_U^*; \theta)} \frac{\mathbb{E}_\Theta[\widehat{\mathcal{I}}_U - \widehat{\mathcal{I}}_L]}{I_U^* - \mathcal{I}(\mathcal{T}_\mu(\Theta)) + 1} \leq 1 + o(1) + \frac{I_U^* - I_L^*}{I_U^* - \mathcal{I}(\mathcal{T}_\mu(\Theta)) + 1} \cdot O(1/d^3).$$

Since $I_U^* - I_L^* = O(d^2)$, we have the second term on the right hand side above to be $o(1)$. Therefore, we finish the proof of (4.10).

Similarly, we can also have

$$\mathbb{E}_\Theta[\mathcal{I}(\Theta) - \widehat{\mathcal{I}}_L] \leq [\mathcal{I}(\Theta) - \mathcal{I}(\mathcal{T}_\mu(\Theta)) \vee I_L^*]\mathbb{P}(\mathcal{I}(\mathcal{T}_\mu(\Theta)) \leq \widehat{\mathcal{I}}_L) + [\mathcal{I}(\Theta) - I_L^*]\mathbb{P}(\widehat{I}_L > \mathcal{I}(\mathcal{T}_\mu(\Theta))).$$

By (S.6), following the same argument as the proof of (4.10) above, we prove (4.11).

Theorem S.7 is a direct corollary of (S.6). For a monotone property $\mathcal{P}$, we have the range $I_L^* = 0$ and $I_U^* = 1$. Moreover, if $\boldsymbol{\Theta} \in \mathcal{G}_1(\mu; \mathcal{P})$ defined in (S.20), we have $\mathcal{T}_\mu(\boldsymbol{\Theta}) = \boldsymbol{\Theta}$. Since Algorithm 2 is derived from Algorithm 1 in the sense that the test $\psi_\alpha = 1$ is equivalent to a $1 - 2\alpha$ confidence lower side $\widehat{\mathcal{I}}_L = 1$. Thus, under the alternative that $\mathcal{P}(\boldsymbol{\Theta}) = 1$, we have

$$\{\psi_\alpha = 1\} = \{\mathcal{P}(\boldsymbol{\Theta}) \le \widehat{\mathcal{I}}_L\} = \{\mathcal{P}(\mathcal{T}_\mu(\boldsymbol{\Theta})) \le \widehat{\mathcal{I}}_L\}.$$

Therefore, (S.21) can be directly derived from (S.6).

### S.1.3 Proof of Proposition 3.2

We first prove a general result on when the critical edge set is empty. We claim that given any monotone property $\mathcal{P}$, if $\mathcal{P}(G_0) = 1$ then $\mathcal{C}_\mathcal{P}(E_0) = \varnothing$. This is straightforward from Definition 3.1. If $\mathcal{P}(G_0) = 1$, by the monotonicity of $\mathcal{P}$, for any $E' \supseteq E_0$, we have $\mathcal{P}(E') = 1$. Therefore, there is no edge $e \in E_0^c$ such that $\mathcal{P}(E' \backslash \{e\}) = 0$. This implies that $\mathcal{C}_\mathcal{P}(E_0) = \varnothing$. Therefore, we prove the four results in Proposition 3.2 when the critical edge sets are empty. In the following part of the proof, we will discuss four properties case by case. An illustration of the proof are shown in Figure 8.

• **Connected subgraphs.** Now we prove the results on the connected subgraphs. We can consider the case $\mathcal{P}_{\text{Conn},-k}(G_0) = 0$. We choose the $k'$ nodes, each of which from one of the connected subgraphs $\{G_{0\ell} = (V_{0\ell}, E_{0\ell})\}_{\ell=1}^{k'}$. Namely, we choose arbitrary $v_\ell \in V_{0\ell}$ for $1 \le \ell \le k'$. For any edge $(u, v) \in E_0$ such that $j \in V_{0\ell}, k \in V_{0\ell'}, \ell \neq \ell'$. We arbitrarily select $k' - k - 1$ edges from the set $\{(v_s, v_{s'}) | (s, s') \neq (\ell, \ell') \text{ or } (s, s') \neq (\ell', \ell)\}$ and denote this set $\widetilde{E}$. We construct the set

$$E' = E_0 \cup \{e\} \cup \widetilde{E}.$$

We can find an illustration of the construction in Figure 8(a). Notice that $G' = (V, E')$ has $k-1$ connected subgraphs and thus $\mathcal{P}_{\text{Conn},-k}(E') = 1$. We can also check that $E' \backslash \{e\} = E_0 \cup \widetilde{E}$ with has $k$ connected subgraphs. By (3.3), we have $(u,v) \in \mathcal{C}_{\mathcal{P}_{\text{Conn},-k}}(E_0)$ and thus

$$\mathcal{C}_{\mathcal{P}_{\text{Conn},-k}}(E_0) \supseteq \{(u,v) \in E_0 \,|\, j \in V_{0\ell}, k \in V_{0\ell'}, \ell \neq \ell'\}.$$

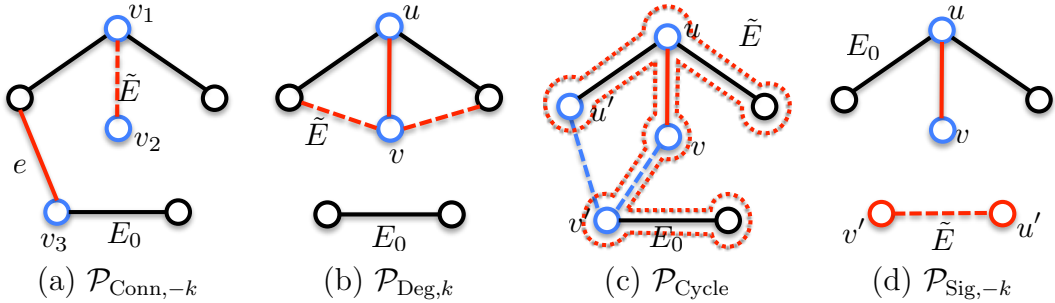

(a) $\mathcal{P}_{\text{Conn},-k}$    (b) $\mathcal{P}_{\text{Deg},k}$    (c) $\mathcal{P}_{\text{Cycle}}$    (d) $\mathcal{P}_{\text{Sig},-k}$

Figure 8: Visualization of the proof of Proposition 3.2.

To prove the other direction, if $e \notin \{(u,v) \in E_0 \mid u \in V_{0\ell}, v \in V_{0\ell'}\}$, we must have $e \in E_{0\ell}$ for some $\ell$. For any $E' \supseteq E$ such that $\mathcal{P}_{\text{Conn},-k}(E') = 1$, since $e \in E_{0\ell}$, $E'\backslash\{e\}$ has the same number of connected subgraphs, therefore $e \notin \mathcal{C}_{\mathcal{P}_{\text{Conn},-k}}(E_0)$. In summary, we prove (3.8).

• **Maximum degree.** We next move to the critical edge set of maximum degree. By (3.3), we have $\mathcal{C}_{\mathcal{P}_{\text{Deg},k}}(E_0) \subseteq E_0^c$. It suffices to prove another direction. Denote the maximum degree of $G_0$ as $k'$. For any $(u,v) \in E_0^c$, since the maximum degree $k' \leq k$, we can select arbitrary $k' - k$ edges from $\{(u, u') \in E_0^c | u' \neq v\}$ and denote the set containing these edges as $\widetilde{E}$. See an illustration of the proof in Figure 8(b). It is easy to check that $E' = E_0 \cup \{(u,v)\} \cup \widetilde{E}$ has $\mathcal{P}_{\text{Deg},k}(E') = 1$ and $\mathcal{P}_{\text{Deg},k}(E'\backslash\{(u,v)\}) = 0$. Therefore, we have $\mathcal{C}_{\mathcal{P}_{\text{Deg},k}}(E_0) = E_0^c$.

• **Acyclic.** Now we turn to show the critical edge set for the cyclic property when $G_0 = (V, E_0)$ is a forest. Similar to the preceding part, it suffices to show $E_0^c \subseteq \mathcal{C}_{\mathcal{P}_{\text{Cycle}}}(E_0)$. For any $(u,v) \in E_0^c$, if $E_0 \cup \{(u,v)\}$ forms a cycle, then $(u,v) \in \mathcal{C}_{\mathcal{P}_{\text{Cycle}}}(E_0)$ by the definition in (3.3). If $E_0 \cup \{(u,v)\}$ is still a forest, we apply the Kruskal's algorithm (Kruskal, 1956) with the initial input $E_0 \cup \{(u,v)\}$ and the weights of all edges are 1. Specifically, the procedure in described in Algorithm 6 and in the algorithm we can choose arbitrary order of adding edges in $(E_0 \cup \{(u,v)\})^c$ in the for loop since the edge weights are the same. We denote the output edge set of Algorithm 6 as $\widetilde{E}$ and $(V, \widetilde{E})$ is a tree by the property of Kruskal's algorithm. We illustrate the graph and the following of the construction in Figure 8(c). We start to construct $E'$ to satisfy (3.3). Since $(V, E')$ is a tree, there are only two cases: (1) neither of $u, v$ is a leaf; and (2) only one of $u, v$ is a leaf. We first consider the case that neither of $u, v$ is a leaf in $\widetilde{E}$. Then, in the graph $(V, \widetilde{E})$, we can choose $u'$ as any neighbor of $u$ but $v$ and choose $v'$ as any neighbor of $v$ but $u$. This is feasible since neither of $u, v$ is a leaf. We can see that $u' \neq v'$, otherwise the triangle $u \to v \to u' \to u$ forms a loop. Let $E' = \widetilde{E} \cup \{(u', v')\}$ and it has a loop $u \to v \to v' \to u' \to u$ and has only this loop (as $\widetilde{E}$ is a tree and adding one edge to it will only form one loop). Therefore, if we delete the edge $(u,v)$ from the loop, the graph $(V, E'\backslash\{(u,v)\})$ has no loop. Therefore, we check (3.3) that $E' = \widetilde{E} \cup \{(u, u')\}$ has $\in \mathcal{P}_{\text{Cycle}}(E') = 1$ and $\mathcal{P}_{\text{Cycle}}(E'\backslash\{(u, u')\}) = 0$. For the second case that only one of $u, v$ is a leaf (we assume the leaf is $u$), the proof is similar to the first one. In the graph $(V, \widetilde{E})$, we find $v'$ as any neighbor of $v$ but $u$. We can check as the first case that $E' = \widetilde{E} \cup \{(u, v')\}$ has $\mathcal{P}_{\text{Cycle}}(E') = 1$ and $\mathcal{P}_{\text{Cycle}}(E'\backslash\{(u, v')\}) = 0$. We can now conclude that $E_0^c = \mathcal{C}_{\mathcal{P}_{\text{Cycle}}}(E_0)$.

---

**Algorithm 6** Kruskal's algorithm for the proof of Proposition 3.2

---

**Input:** $E^{(0)} = E_0 \cup \{(u,v)\}, t = 0$.
  **for** $e \in (E_0 \cup \{(u,v)\})^c$ **do**
    $t \leftarrow t + 1$;
    **if** $\{e\} \cup E^{(t)}$ does not contain a cycle **then**
      $E^{(t)} \leftarrow E^{(t-1)} \cup \{e\}$
    **end if**
  **end for**
**Output:** $E^{(t)}$.

---

• **Singletons.** We finally discuss the critical edge set for isolated nodes. Denote the number

of isolated nodes $|V_{\text{Sig}}| = k' \geq k$. For any $(u,v) \in E_0^c$ such that $u \in V_{\text{Sig}}$ or $v \in V_{\text{Sig}}$, we can check that at least one of $u$ and $v$ is an isolated node. If both $u,v$ are isolated nodes, this implies that $k' > k+1$ and we can select arbitrary $k'-k-1$ nodes from $V_{\text{Sig}}\backslash\{u,v\}$ if only one of $u,v$ is a isolated node, we select arbitrary $k'-k-2$ nodes from $V_{\text{Sig}}\backslash\{u,v\}$. Denote the set of nodes selected as $\widetilde{V}$ and define $\widetilde{E} = \{(u',v') \in E_0^c | u',v' \in \widetilde{V}\}$. We can check that $E' = E_0 \cup \{(u,v)\} \cup \widetilde{E}$ has $\mathcal{P}_{\text{Sig},-k}(E') = 1$ and $\mathcal{P}_{\text{Sig},-k}(E'\backslash\{(u,v)\}) = 0$. Therefore, we have $\mathcal{C}_{\mathcal{P}_{\text{Sig},-k}}(E_0) \supseteq \{(u,v) \in E_0^c | u \in V_{\text{Sig}} \text{ or } v \in V_{\text{Sig}}\}$. The construction is illustrated in Figure 8(d). On the other hand, let the edge $(u,v)$ satisfy $u \notin V_{\text{Sig}}$ and $v \notin V_{\text{Sig}}$. For any $E' \supseteq E_0$ and $\mathcal{P}_{\text{Sig},-k}(E') = 1$, as $E'\backslash\{(u,v)\}$ does not include new isolated nodes, $\mathcal{P}_{\text{Sig},-k}(E'\backslash\{(u,v)\}) = 1$. This implies that $\mathcal{C}_{\mathcal{P}_{\text{Sig},-k}}(E_0) = \{(u,v) \in E_0^c | u \in V_{\text{Sig}} \text{ or } v \in V_{\text{Sig}}\}$.

## S.2 Proofs for Lower Bound of confidence interval length

In this section, we give a general framework of the lower bound of confidence interval length.

### S.2.1 Proof of Theorem 6.1

We define the pre-distance on a graph $G = (V,E)$ between two vertex sets $V_1$ and $V_2$ as $d_G(V_1, V_2) = $ the length of the shortest path on $G$ connecting one of $v_1 \in V_1$ and one of $v_2 \in V_2$, and if there is no such path, we let $d_G(V_1, V_2) = \infty$. In order to prove Theorem 6.1, we need the following lemma.

**Lemma S.2.** Given an invariant, if there exists a construction as follows:

(1) There exists $\overline{G}_L = (V, \overline{E}_L)$ satisfying $\overline{V}_L := V(\overline{E}_L) \subset V = \{1,\ldots,d\}$ and $\mathcal{I}(\overline{G}_L) = I_L^*$.

(2) There exist $N$ disjoint anchor vertex sets $\{A_j\}_{j=1}^N$ such that $\cup_{j=1}^N A_j \subseteq \overline{V}_L$.

(3) There exists $\overline{V}_U \subseteq V \backslash \overline{V}_L$ such that for each vertex set $A_j$, we can find an edge set $E_j$ with $V(E_j) \subseteq A_j \cup \overline{V}_U$ and the graph $G_j = (V, \overline{E}_L \cup E_j)$ has $\mathcal{I}(G_j) = I_U^*$.

(4) There exists a universal constant $R$ independent to the such that for any $1 \leq j \leq N$, $G_j$ is $R$-hollow and $|A_j| \leq R$. Moreover, the number of anchor vertex sets has $N \geq d^\gamma$ with some $\gamma > 0$.

Suppose $\forall 1 \leq j \neq k \leq N$, $d_{\overline{G}_L}(A_j, A_k) = \infty$, $|\overline{V}_U| = o(d^{1/2})$ and $|\overline{V}_L|/d < 1$. there exist constants $C_1$ and $C_2$ such that if

$$\theta \leq C_1\sqrt{\log d/n} \text{ and } \max_j d_{\max}(G_j)\sqrt{\log d/n} \leq C_2, \qquad (S.7)$$

we have the following lower bound on the confidence interval length

$$\liminf_{n\to\infty} \inf_{[\widehat{L},\widehat{U}] \in I(\mathcal{I},\alpha)} \sup_{\Theta \in \mathcal{U}_s(I_L^*, I_U^*;\theta)} \frac{\mathbb{E}_\Theta[\widehat{U} - \widehat{L}]}{\text{Oracle Length}(\Theta)} \geq 1 - 2\alpha.$$

We now prove Theorem 6.1 using Lemma S.2. Let

$$\overline{G}_L = \cup_{j=1}^{N} G_{L,j} \cup G_L \text{ and } \overline{G}_U = \cup_{j=1}^{N} G_{L,j} \cup G_U.$$

We choose anchor vertex sets $A_j = V(G_{L,j})$ for $1 \leq j \leq N$. Since $|V(E_L)| = O(1)$ and $G_{L,j}$'s are isomorphic copies of $G_L$, we have $|A_j| = O(1)$ as well. It is easy to check $\cup_{j=1}^{N} A_j \subseteq \overline{V}_L$. We now construct $G_j$'s as isomorphic copies of $\overline{G}_U$. In order to obtain these copies, we denote the isomorphic map between $G_L$ and $G_{L,j}$ as $\sigma_j : V(G_L) \to V(G_{L,j})$. We extend the domain of $\sigma_j$ to $V$ by letting $\bar{\sigma}_j(u) = \sigma_j(u)$ if $u \in V(G_L)$ and $\bar{\sigma}_j(u) = u$ otherwise. We construct $E_j := \{(\sigma_j(u), \sigma_j(v)) \,|\, (u,v) \in (E_U \backslash E_L))\}$ and we can see the graph $G_j = (V, \overline{E}_L \cup E_j)$ is isomorphic to $\overline{G}_U$. Therefore, $\mathcal{I}(G_j) = I_U^*$ for any $1 \leq j \leq N$. Therefore, by Lemma S.2, we complete the proof.

### S.2.2 Proof of Lemma S.2

The high level idea of the proof is to reduce the length of confidence interval to the bound of $\chi^2$-divergence of two distributions.

In fact, we will prove a stronger lower bound in a smaller parameter space $\mathcal{U}_s(I_L^*, I_U^*; \theta, \mu)$ defined in (S.5). We consider this parameter space because it contains precision matrices $\Theta$ such that $\mathcal{I}(\mathcal{T}_\mu(\Theta)) \geq I_L^*$, for $\mu = C\sqrt{\log d/n}$. Namely, part of entries of precision matrix has strong enough signal strenght.

We can lower bound the length of confidence interval as

$$\inf_{[\widehat{L},\widehat{U}]\in I(\mathcal{I},\alpha)} \sup_{\Theta \in \mathcal{U}_s(I_L^*, I_U^*;\theta)} \frac{\mathbb{E}_\Theta[\widehat{U}-\widehat{L}]}{I_U^* - \mathcal{I}(\mathcal{T}_\mu(\Theta))} \geq \inf_{[\widehat{L},\widehat{U}]\in I(\mathcal{I},\alpha)} \sup_{\Theta \in \mathcal{U}_s(I_L^*, I_U^*;\theta,\mu)} \frac{\mathbb{E}_\Theta[\widehat{U}-\widehat{L}]}{I_U^* - \mathcal{I}(\mathcal{T}_\mu(\Theta))}$$

$$\geq \inf_{[\widehat{L},\widehat{U}]\in I(\mathcal{I},\alpha)} \sup_{\Theta \in \mathcal{U}_s(I_L^*, I_U^*;\theta,\mu)} \frac{\mathbb{E}_\Theta[\widehat{U}-\widehat{L}]}{I_U^* - I_L^*}. \quad (S.8)$$

The following lemma proved in Cai and Guo (2015) reduces the length of confidence interval to the $\chi^2$-divergence. Given two distributions $\mathbb{P}$ and $\mathbb{Q}$, the $\chi^2$-divergence between $\mathbb{P}$ and $\mathbb{Q}$ is defined as

$$D_{\chi^2}(\mathbb{P}, \mathbb{Q}) = \int \left(\frac{d\mathbb{P}}{d\mathbb{Q}}\right)^2 d\mathbb{Q} - 1.$$

**Lemma S.3** (Lemma 1, Cai and Guo (2015)). *Given any monotone invariant $\mathcal{I}$, suppose there exist $\Theta_0, \Theta_1, \ldots, \Theta_N \in \mathcal{U}_s(I_L^*, I_U^*; \theta, \mu)$ satisfying $\mathcal{I}(\Theta_0) = I_L^*$ and $\mathcal{I}(\Theta_j) = I_U^*$ for all $1 \leq j \leq N$. If $\mathbb{P}_L = \mathbb{P}_{\Theta_0}$ and $\overline{\mathbb{P}}_U = \frac{1}{N}\sum_{j=1}^{N} \mathbb{P}_{\Theta_j}$, we have*

$$\inf_{[\widehat{L},\widehat{U}]\in I(\mathcal{I},\alpha)} \sup_{\Theta \in \mathcal{U}_s(I_L^*, I_U^*;\theta,\mu)} \mathbb{E}_\Theta[\widehat{U}-\widehat{L}] \geq (I_U^* - I_L^*)\left(1 - 2\alpha - \sqrt{D_{\chi^2}(\overline{\mathbb{P}}_U, \mathbb{P}_L)}\right).$$

Given two edge sets $E_1, E_2$, we denote the pre-distance on the graph $G$ between $E_1$ and $E_2$ as $d_G(V_1, V_2) = d_G(V(E_1), V(E_2))$. Combining (S.8) with Lemma S.3, the proof of Lemma S.2 can be deduced from the following theorem.

**Theorem S.4.** Given graphs $G_0 = (V, E_0)$ and $\overline{G}_1, \ldots, \overline{G}_N$, we denote the difference edge set as $\boldsymbol{C} = \{E(\overline{G}_j) \backslash E(G_0)\}_{j=1}^N$. Given an edge set $S \in \boldsymbol{C}$, we denote $\mathbf{A}_S$ as the adjacency matrix of $S$. We also define $\mathbf{A}_0$ as the adjacency matrix of $G_0$ and $\mathbf{A}_{S,S'} = \mathbf{A}_0 + \mathbf{A}_S + \mathbf{A}_{S'}$. Let $\mathcal{V}_{S,S'} = \{V(E_0 \cup S) \cap V(S')\} \cup \{V(E_0 \cup S') \cap V(S)\}$. Suppose the following assumption holds:

**A1:** Denote the uniform maximum degree as $\Gamma = \max_{S,S' \in \boldsymbol{C}} \|\mathbf{A}_{S,S'}\|_1$ and uniform spectral norm as $\Lambda = \max_{S,S' \in \boldsymbol{C}} \|\mathbf{A}_{S,S'}\|_2$. We further define

$$\mathcal{R} = \max_{S,S' \in \boldsymbol{C}} \frac{|S \cap S'|}{|\mathcal{V}_{S,S'}|} \quad \text{and} \quad \mathcal{B} = \max_{S,S' \in \boldsymbol{C}} \left( (\Gamma^2 |\mathcal{V}_{S,S'}|) \wedge \Lambda^4 \right).$$

Suppose $S'$ is uniformly sampled from $\boldsymbol{C}$. If for any fixed $S \in \boldsymbol{C}$, $V(S)$ can be split into $\ell$ groups: $V(S) = \cup_{j=1}^\ell V_j(S)$, so that the random variables $\left\{ \mathbb{1}(\mathcal{V}_{S,S'} \cap V_j(S) \neq \varnothing) \right\}_{j \in [\ell]}$ with respect to a uniformly sampled $S'$ from $\boldsymbol{C}$ are negatively associated. In other words, for any pair of disjoint sets $I, J \subseteq [\ell]$ and any pair of coordinate-wise nondecreasing functions $f, g$ we have:

$$\text{Cov}\left( f\left( \left\{ \mathbb{1}(\mathcal{V}_{S,S'} \cap V_j(S) \neq \varnothing) \right\}_{j \in I} \right), g\left( \left\{ \mathbb{1}(\mathcal{V}_{S,S'} \cap V_j(S) \neq \varnothing) \right\}_{j \in J} \right) \right) \leq 0.$$

Denote the largest cardinality of vertices as $V_{\max} = \max_{S \in \boldsymbol{C}} \max_j |V_j(S)|$. Assume the following holds: $\left[ \max_{S \in \boldsymbol{C}} \mathbb{E}_{S'} |\mathcal{V}_{S,S'}| \right]^{-1} \to \infty$,

$$\mu \leq \sqrt{\frac{\mathcal{R}}{\mathcal{B}}} \wedge \frac{1 - C^{-1}}{2\sqrt{2}\Gamma} \quad \text{and} \quad \theta \leq \sqrt{\frac{\log \left( \left[ \max_{S \in \boldsymbol{C}} \mathbb{E}_{S'} |\mathcal{V}_{S,S'}| \right]^{-1} \right)}{4n V_{\max} \mathcal{R}}} \wedge \mu.$$

Under either Setting A or B, for the parameters $\boldsymbol{\Theta}_0 = \mu \mathbf{A}_0$, $\boldsymbol{\Theta}_S = \mu \mathbf{A}_0 + \theta \mathbf{A}_S$ for all $S \in \boldsymbol{C}$, the $\chi^2$-divergence between $\mathbb{P}_L = \mathbb{P}_{\boldsymbol{\Theta}_0}$ and $\overline{\mathbb{P}}_U = \frac{1}{N} \sum_{S \in \boldsymbol{C}} \mathbb{P}_{\boldsymbol{\Theta}_S}$ satisfies

$$\lim_{n \to \infty} D_{\chi^2}(\overline{\mathbb{P}}_U, \mathbb{P}_L) = 0.$$

We defer the proof of Theorem S.4 to Section S.2.4 and continue the proof of Lemma S.2. Our strategy is to show that the conditions in Lemma S.2 satisfy the ones in Theorem S.4 correspondingly. We first bound the norms of the adjacency matrices in Theorem S.4. Let $\overline{G}_L$ and $\{E_j\}_{j=1}^N$ be the construction in Lemma S.2. We have $\Gamma = \max_{1 \leq j,k \leq N} \|\mathbf{A}_{E_j,E_k}\|_1 \leq 2 \max_j d_{\max}(G_j)$. Moreover, since $\{G_j\}_{j=1}^N$ are $R$-hollow, by Definition 6.1, we also have the graphs $(V, E_j)$ and $(V, E_j \cup E_k)$ are $2R$-hollow for any $1 \leq j \neq k \leq N$. According to Lemma S.6, we have

$$\max_{1 \leq j \leq N} \|\mathbf{A}_{E_j}\|_2 \leq 4R \sqrt{\max_j d_{\max}(G_j)} \quad \text{and} \quad \max_{1 \leq j,k \leq N} \|\mathbf{A}_{E_j,E_k}\|_2 \leq 4R \sqrt{2 \max_j d_{\max}(G_j)}. \quad (S.9)$$

Moreover, as $|A_j| \leq R$, since $G_j$ is hollow, we also have $|E_j| \leq R|A_j| \leq R^2$. Now we verify

**A1**. First we describe how to construct $G_0, \overline{G}_1, \ldots, \overline{G}_N$. We begin by splitting the vertices into two parts $\overline{V}_L$ and $V\backslash \overline{V}_L$. Let $G_0 = \overline{G}_L = (V, \overline{E}_L)$. Next, we randomly sample $|\overline{V}_U|$ vertices from the vertex set $V\backslash \overline{V}_L$. Denote as $\widetilde{V}$ the random vertex set we have sampled. We also sample one anchor set $A_{\widetilde{j}}$ where $\widetilde{j}$ follows uniform distribution from $\{1, \ldots, N\}$. We aim to connect $\widetilde{V}$ to $\overline{V}_L$ using the pattern of $G_{\widetilde{j}}$. In specific, since $|\widetilde{V}| = |\overline{V}_U|$, there exists a unique one-to-one mapping $\pi : A_{\widetilde{j}} \cup \overline{V}_U \to A_{\widetilde{j}} \cup \widetilde{V}$ such that $\pi(u) = u$ if $u \in A_{\widetilde{j}}$ and if $u \in \overline{V}_L$, $\pi(u) \in \widetilde{V}$ with $\pi(u) < \pi(v)$ if and only if $u < v$ for any $u, v \in \overline{V}_U$. Given $\widetilde{V}$ and $A_{\widetilde{j}}$, we construct the edge set

$$\widetilde{E} = \{(\pi(u), \pi(v)) \mid (u, v) \in \overline{E}_L \cup E_{\widetilde{j}}\}.$$

Then the graph $\widetilde{G} = (V, \widetilde{E} \cup \overline{E}_L)$ is isomorphic to $G_{\widetilde{j}}$. We construct $\overline{G}_1, \ldots, \overline{G}_N$ for all possible such $\widetilde{G}$ by enumerating $\widetilde{V}$ and $A_{\widetilde{j}}$. The edge set $\boldsymbol{C}$ then contains the edge set $\widetilde{E}$ for all possible vertex set $\widetilde{V} \subseteq V\backslash \overline{V}_L$ and anchor set $A_{\widetilde{j}}$.

We will now prove that the quantity $\mathcal{R}$ in (6.1) is bounded. Due to the construction of $\boldsymbol{C}$ for all $S, S' \in \boldsymbol{C}$ and the fact that $d_{G_0}(E_j, E_k) = \infty$ for all $\forall 1 \leq j \neq k \leq N$, we have $\mathcal{V}_{S,S'} = V(S) \cap V(S')$. Therefore, we obtain that

$$\frac{1}{2} \leq \max_{S \in \boldsymbol{C}} \frac{|S|}{|V(S)|} \leq \max_{S,S' \in \boldsymbol{C}} \frac{|S \cap S'|}{|\mathcal{V}_{S,S'}|} = \max_{S,S' \in \boldsymbol{C}} \frac{|S \cap S'|}{|V(S) \cap V(S')|} \leq R, \tag{S.10}$$

where we used the fact that the graph $(V(S) \cap V(S'), S \cap S')$ is isomorphic to a subset of the $R$-hollow graph $G_j$ for some $1 \leq j \leq N$. Therefore, $1/2 \leq \mathcal{R} \leq R$.

Next we show that there exists a division of $V(S)$ satisfying the negative association. Given any fixed $S \in \boldsymbol{C}$, we split $V(S)$ into $V(S) \cap \overline{V}_L$ and $V(S) \cap \overline{V}_L^c$. Let $S'$ be uniformly sampled from $\boldsymbol{C}$. Since $\mathcal{V}_{S,S'} = V(S) \cap V(S')$, by the construction above, it is simple to see that, the variables

$$\{\mathbb{1}(\mathcal{V}_{S,S'} \cap V' \neq \varnothing) \mid V' = V(S) \cap \overline{V}_L \text{ or } V' = \{v\} \text{ for } v \in V(S) \cap \overline{V}_L^c\}$$
$$= \{\mathbb{1}(V(S) \cap V(S') \cap \overline{V}_L = \varnothing)\} \cup \{\mathbb{1}(v \in V(S)) \mid v \in V(S') \cap \overline{V}_L^c\}.$$

The above variables are negatively associated, since by our construction, $\mathbb{1}(V(S) \cap V(S') \cap \overline{V}_L = \varnothing)$ is independent of the remaining random variables, and the others are induced by a hypergeometric distribution (Joag-Dev and Proschan, 1983). It is further simple to evaluate the expectation:

$$\mathbb{E}_{S'}|\mathcal{V}_{S,S'}| = \mathbb{E}_{S'}|V(S) \cap V(S')|$$
$$= \mathbb{E}_{S'}|V(S) \cap V(S') \cap \overline{V}_L| + \sum_{v \in V(S') \cap \overline{V}_L^c} \mathbb{E}_{S'}[\mathbb{1}(v \in V(S))] \leq \frac{\max_j |A_j|}{N} + \frac{|\overline{V}_U|^2}{d - |\overline{V}_L|}.$$

Since we assume $\max_j |A_j| \leq R$, $N \geq d^\gamma$, $|\overline{V}_U| = o(d^{1/2})$ and $|\overline{V}_L|/d < 1$, we have

$$\log \big[\max_{S \in \boldsymbol{C}} \mathbb{E}_{S'} |\mathcal{V}_{S,S'}|\big]^{-1} \geq C \log d \to \infty.$$

By (S.9) and $\Gamma \leq 2 \max_j d_{\max}(G_j)$, we have

$$\mathcal{B} \leq \Lambda^4 \leq CR^4 \max_j d_{\max}^2(G_j).$$

Combining with $V_{\max} \leq \max_j |A_j| \leq R$ and $1/2 \leq \mathcal{R} \leq R$ by (S.10), we can reduce the rate in **A1** to

$$\sqrt{\frac{\mathcal{R}}{\mathcal{B}}} \wedge \frac{1 - C^{-1}}{2\sqrt{2}\Gamma} \geq \frac{CR^{-3/2}}{\max_j d_{\max}(G_j)}$$

$$\sqrt{\frac{\log \big(\big[\max_{S \in \boldsymbol{C}} \mathbb{E}_{S'} |\mathcal{V}_{S,S'}|\big]^{-1}\big)}{4nV_{\max}\mathcal{R}}} \wedge \mu \geq C\sqrt{\frac{\log d}{n}} \wedge \mu.$$

Therefore, the construction in Lemma S.2 can be reduced to **A1** in Theorem S.4, and therefore for $\boldsymbol{\Theta}_0 = \mu\mathbf{A}_0$, $\boldsymbol{\Theta}_j = \mu\mathbf{A}_0 + \theta\mathbf{A}_{E_j}$ for $1 \leq j \leq N$ we have

$$\lim_{n \to \infty} D_{\chi^2}(\overline{\mathbb{P}}_U, \mathbb{P}_L) = 0, \text{ where } \mathbb{P}_L = \mathbb{P}_{\boldsymbol{\Theta}_0} \text{ and } \overline{\mathbb{P}}_U = \frac{1}{N} \sum_{j=1}^N \mathbb{P}_{\boldsymbol{\Theta}_j}.$$

Combining (S.8) with Lemma S.3 again, we prove Lemma S.2.

### S.2.3 Proofs of Theorem 6.2

The proof of Theorem 6.2 is same as the one of Lemma S.2. We just set $\theta$ in the proof of Lemma S.2 as $\mu = \theta$. We change the parameters in Theorem S.4 into $\boldsymbol{\Theta}_0 = \theta\mathbf{A}_0$, $\boldsymbol{\Theta}_S = \theta\mathbf{A}_0 + \theta\mathbf{A}_S$ for all $S \in \boldsymbol{C}$. Considering two distributions $\mathbb{P}_L = \mathbb{P}_{\boldsymbol{\Theta}_0}$ and $\overline{\mathbb{P}}_U = \frac{1}{N} \sum_{j=1}^N \mathbb{P}_{\boldsymbol{\Theta}_j}$, Theorem S.4 gives us $\lim_{n \to \infty} D_{\chi^2}(\overline{\mathbb{P}}_U, \mathbb{P}_L) = 0$. Applying Lemma S.3, we have

$$\inf_{[\widehat{L},\widehat{U}] \in I(\mathcal{I},\alpha)} \sup_{\boldsymbol{\Theta} \in \mathcal{U}_s(I_L^*, I_U^*; \theta)} \mathbb{E}_{\boldsymbol{\Theta}}[\widehat{U} - \widehat{L}] \geq (I_U^* - I_L^*)\Big(1 - 2\alpha - \sqrt{D_{\chi^2}(\overline{\mathbb{P}}_U, \mathbb{P}_L)}\Big) \to (I_U^* - I_L^*)(1 - 2\alpha).$$

We now start to prove Theorem 6.2.

$$\inf_{\widehat{U} \in U(\mathcal{I},\alpha)} \sup_{\boldsymbol{\Theta} \in \mathcal{U}_s(I_L^*, I_U^*; \theta, \mu)} \mathbb{E}_{\boldsymbol{\Theta}}\big[\widehat{U} - \mathcal{I}(\boldsymbol{\Theta})\big] \geq \inf_{\widehat{U} \in U(\mathcal{I},\alpha)} \sup_{\mathcal{I}(\boldsymbol{\Theta}) = I_L^*} \mathbb{E}_{\boldsymbol{\Theta}}\big[\widehat{U}\big] - I_L^*$$

$$\geq \inf_{\widehat{U} \in U(\mathcal{I},\alpha)} \sup_{\mathcal{I}(\boldsymbol{\Theta}) = I_L^*} (I_U^* \mathbb{P}_{\boldsymbol{\Theta}}\big[\widehat{U} = I_U^*\big] + I_L^* \mathbb{P}_{\boldsymbol{\Theta}}\big[\widehat{U} < I_U^*\big]) - I_L^*$$

$$\geq \inf_{\widehat{U} \in U(\mathcal{I},\alpha)} \sup_{\mathcal{I}(\boldsymbol{\Theta}) = I_L^*} (I_U^* - I_L^*) \mathbb{P}_{\boldsymbol{\Theta}}\big[\widehat{U} = I_U^*\big]. \tag{S.11}$$

We aim to bound $\mathbb{P}_{\boldsymbol{\Theta}}\big[\widehat{U} = I_U^*\big]$ by Theorem S.4. Since $\mathcal{I}$ is hollow, similarly to the proof of Lemma S.2, we can construct $\boldsymbol{\Theta}_0, \boldsymbol{\Theta}_1, \ldots, \boldsymbol{\Theta}_N \in \mathcal{U}_s(I_L^*, I_U^*; \theta, \mu)$ satisfying $\mathcal{I}(\boldsymbol{\Theta}_0) = I_L^*$ and $\mathcal{I}(\boldsymbol{\Theta}_j) = I_U^*$ for all $1 \leq j \leq N$ and the $\chi^2$-divergence between $\mathbb{P}_L = \mathbb{P}_{\boldsymbol{\Theta}_0}$ and $\overline{\mathbb{P}}_U = \frac{1}{N} \sum_{S \in \boldsymbol{C}} \mathbb{P}_{\boldsymbol{\Theta}_S}$ satisfies

$$\lim_{n \to \infty} D_{\chi^2}(\overline{\mathbb{P}}_U, \mathbb{P}_L) = 0. \tag{S.12}$$

Next notice that

$$|\mathbb{P}_L(\widehat{U} < I_U^*) - \overline{\mathbb{P}}_U(\widehat{U} < I_U^*)| \leq \text{TV}(\mathbb{P}_L, \overline{\mathbb{P}}_U) \leq \sqrt{D_{\chi^2}(\mathbb{P}_L, \overline{\mathbb{P}}_U)},$$

where the next to last inequality simply follows by the definition of total variation norm and the last inequality is due to Cauchy-Schwartz. Hence we have:

$$\sup_{\mathcal{I}(\boldsymbol{\Theta}) = I_L^*} \mathbb{P}_{\boldsymbol{\Theta}}\big[\widehat{U} = I_U^*\big] \geq 1 - \mathbb{P}_L(\widehat{U} < I_U^*)$$

$$\geq 1 - \sqrt{D_{\chi^2}(\mathbb{P}_L, \overline{\mathbb{P}}_U)}, -\overline{\mathbb{P}}_U(\widehat{U} < I_U^*)$$

$$\geq 1 - \sup_{\mathcal{I}(\boldsymbol{\Theta}) = I_L^*} \mathbb{P}_{\boldsymbol{\Theta}}(\widehat{U} < I_U^*) - \sqrt{D_{\chi^2}(\mathbb{P}_L, \overline{\mathbb{P}}_U)} \geq 1 - \alpha - \sqrt{D_{\chi^2}(\mathbb{P}_L, \overline{\mathbb{P}}_U)}.$$

Combining the above inequality with (S.12) and (S.11), we prove Theorem 6.2.

### S.2.4 Proof of Theorem S.4

The proof of the theorem mainly relies on the following proposition.

**Proposition S.5.** Let $G_0, \overline{G}_1, \ldots, \overline{G}_N$ and $\boldsymbol{C}$ be the same as Theorem S.4. Define the matrices $\widetilde{\mathbf{A}}_0 = (\mu/\theta)\mathbf{A}_0$, and $\widetilde{\mathbf{A}}_{S,S'} = \widetilde{\mathbf{A}}_0 + \mathbf{A}_S + \mathbf{A}_{S'}$. For all $S, S' \in \boldsymbol{C}$ we denote:

$$\mathcal{V}_{S,S'|S} := V(E_0 \cup S') \cap V(S) = \mathcal{V}_{S,S'} \cap V(S).$$

Then for any collection of vertex buffers $\boldsymbol{V} = \{\mathcal{V}_{S,S'}\}_{S,S' \in \boldsymbol{C}}$ and any of the following two

choices:

**Setting 1:** $\mathcal{H}_{S,S'} = \frac{(|\mathcal{V}_{S,S'|S}| \wedge |\mathcal{V}_{S,S'|S'}|)\|\mathbf{A}_S\|_2\|\mathbf{A}_{S'}\|_2}{\|\widetilde{\mathbf{A}}_{S,S'}\|_2^2}$ and $\mathcal{K}_{S,S'} = 2\|\widetilde{\mathbf{A}}_{S,S'}\|_2$;

**Setting 2:** $\mathcal{H}_{S,S'} = \frac{|\mathcal{V}_{S,S'|S}||\mathcal{V}_{S,S'|S'}|}{\|\widetilde{\mathbf{A}}_{S,S'}\|_1^2}$ and $\mathcal{K}_{S,S'} = 2\|\widetilde{\mathbf{A}}_{S,S'}\|_1$.

when the signal strengths satisfies

$$\theta \leq \mu, \quad \mu \leq \frac{1 - C^{-1}}{2\sqrt{2}\|\mathbf{A}_{S,S'}\|_1}, \tag{S.13}$$

we have for $\mathbf{\Theta}_0 = \mu \mathbf{A}_0$, $\mathbf{\Theta}_S = \mu \mathbf{A}_0 + \theta \mathbf{A}_S$ for all $S \in \boldsymbol{C}$, the $\chi^2$-divergence between $\mathbb{P}_L = \mathbb{P}_{\mathbf{\Theta}_0}$ and $\overline{\mathbb{P}}_U = \frac{1}{N}\sum_{S \in \boldsymbol{C}} \mathbb{P}_{\mathbf{\Theta}_S}$ satisfies

$$D_{\chi^2}(\mathbb{P}_L, \overline{\mathbb{P}}_U) \leq \frac{1}{|\boldsymbol{C}|^2} \sum_{S,S' \in \boldsymbol{C}} \exp\left[n\left(|S \cap S'|\theta^2 + \frac{\mathcal{H}_{S,S'}(\theta\mathcal{K}_{S,S'})^{2(d_{G_0}(S,S')\vee 1+1)}}{2(d_{G_0}(S,S')\vee 1+1)}\right)\right] - 1.$$

We defer the proof of proposition to Section S.2.4.1 and first present the proof of Theorem S.4.

*Proof of Theorem S.4.* Using Proposition S.5 it suffices to control the expression:

$$\frac{1}{|\boldsymbol{C}|^2} \sum_{S,S' \in \boldsymbol{C}} \exp\left[n\left(|S \cap S'|\theta^2 + \frac{\mathcal{H}_{S,S'}(\theta\mathcal{K}_{S,S'})^{2(d_{G_0}(S,S')\vee 1+1)}}{2(d_{G_0}(S,S')\vee 1+1)}\right)\right] - 1,$$

Note that by the definition of $\mathcal{B}$ we have:

$$\max_{S,S' \in \boldsymbol{C}} \left((\|\mathbf{A}_S\|_2\|\mathbf{A}_{S'}\|_2\|\mathbf{A}_{S,S'}\|_2^2) \wedge (|\mathcal{V}_{S,S'}|\|\mathbf{A}_{S,S'}\|_1^2)\right) \leq \mathcal{B}.$$

Therefore using Proposition S.5 it suffices to control:

$$D_{\chi^2}(\mathbb{P}_L, \overline{\mathbb{P}}_U) = \frac{1}{|\boldsymbol{C}|^2} \sum_{S,S' \in \boldsymbol{C}} \exp\left[|\mathcal{V}_{S,S'|S}|n\theta^2\left(\mathcal{R} + \mathcal{B}\mu^2\right)\right].$$

First note that by $|\mathcal{V}_{S,S'|S}| = \sum_{j \in [\ell]} |V_j(S) \cap \mathcal{V}_{S,S'|S}|$, we have:

$$D_{\chi^2}(\mathbb{P}_L, \overline{\mathbb{P}}_U) \leq \frac{1}{|\boldsymbol{C}|^2} \sum_{S,S' \in \boldsymbol{C}} \exp\left[n\theta^2\left(\mathcal{R} + \mathcal{B}\mu^2\right) \sum_{j \in [\ell]} |V_j(S) \cap \mathcal{V}_{S,S'|S}|\right].$$

Denote by $\mathbb{P}_{S'}$ the measure induced by drawing $S'$ uniformly from $\boldsymbol{C}$. Under the assumption: $\mu < \sqrt{\mathcal{R}/\mathcal{B}}$, and using the fact that the random variables $\{\mathbb{1}(\mathcal{V}_{S,S'} \cap V_j(S) \neq \varnothing)\}_{j \in [\ell]}$ are negatively associated for every fixed $S \in \mathbf{C}$, we obtain:

$$\log D_{\chi^2}(\mathbb{P}_L, \overline{\mathbb{P}}_U) \leq \max_{S \in \mathbf{C}} \left[ \sum_{j \in [\ell]} \log \mathbb{E}_{S'} \left[ \exp(2|V_j(S) \cap \mathcal{V}_{S,S'|S}|\mathcal{R}n\theta^2) \right] \right]$$

$$\leq \max_{S \in \mathbf{C}} \left[ \sum_{j \in [\ell]} \log \left[ \exp(2|V_j(S)|\mathcal{R}n\theta^2) \mathbb{P}_{S'}(V_j(S) \cap \mathcal{V}_{S,S'|S} \neq \varnothing) + (1 - \mathbb{P}_{S'}(V_j(S) \cap \mathcal{V}_{S,S'|S} \neq \varnothing)) \right] \right]$$

$$\leq \exp(2\mathcal{R}V_{\max}n\theta^2) \max_{S \in \mathbf{C}} \mathbb{E}_{S'} \sum_{j \in [\ell]} \mathbb{1}(V_j(S) \cap \mathcal{V}_{S,S'|S} \neq \varnothing) \leq \exp(2\mathcal{R}V_{\max}n\theta^2) \max_{S \in \mathbf{C}} \mathbb{E}_{S'} |\mathcal{V}_{S,S'}|,$$

where the expectation $\mathbb{E}_{S'}$ is taken with respect to a uniform draw of $S' \in \mathbf{C}$. The first inequality above is due to negative association, the second inequality is due to $\log(1+x) \leq x$. Recalling that $V_{\max} = \max_{S \in \mathbf{C}} \max_j |V_j(S)|$, we have for values of $\theta$

$$\theta \leq \sqrt{\frac{\log \left( [\max_{S \in \mathbf{C}} \mathbb{E}_{S'}|\mathcal{V}_{S,S'}|]^{1/2}[\max_{S \in \mathbf{C}} \mathbb{E}_{S'}|\mathcal{V}_{S,S'}|]^{-1} \right)}{2nV_{\max}\mathcal{R}}} = \sqrt{\frac{\log \left( [\max_{S \in \mathbf{C}} \mathbb{E}_{S'}|\mathcal{V}_{S,S'}|]^{-1} \right)}{4nV_{\max}\mathcal{R}}},$$

we have $\overline{D}_{\chi^2}(\mathbb{P}_L, \overline{\mathbb{P}}_U) \leq \exp([\max_{S \in \mathbf{C}} \mathbb{E}_{S'}|\mathcal{V}_{S,S'}|]^{1/2})$. The proof is now completed by an application of Lemma S.3. □

### S.2.4.1 Proof of Proposition S.5

The proof of Proposition S.5 follows the same strategy of the proof of Proposition C.1 in Neykov et al. (2016). However, the main difference is that our adjacency matrices $\widetilde{\mathbf{A}}_0 = (\mu/\theta)\mathbf{A}_0$, and $\widetilde{\mathbf{A}}_{S,S'} = \widetilde{\mathbf{A}}_0 + \mathbf{A}_S + \mathbf{A}_{S'}$ have weights on edges while Proposition C.1 in Neykov et al. (2016) handles the adjacency matrices all have the same weight one. This makes our proof different in a few places and to make our proof self-contained, we give a complete proof with both the same and different parts comparing to the proof of Proposition C.1 in Neykov et al. (2016).

In the first step we construct a set of parameters. Define $\boldsymbol{\Theta}_0 = \mathbf{I} + \mu\mathbf{A}_0$, $\boldsymbol{\Theta}_S = \mathbf{I} + \mu\mathbf{A}_0 + \theta\mathbf{A}_S$, $\boldsymbol{\Theta}_{S,S'} = \mathbf{I} + \mu\mathbf{A}_0 + \theta(\mathbf{A}_S + \mathbf{A}_{S'})$, for $S, S' \in \mathbf{C}$ and some $\theta > 0$. For any $S, S' \in \mathbf{C}$ we have:

$$\max(\|\widetilde{\mathbf{A}}_0\|_2, \|\widetilde{\mathbf{A}}_0 + \mathbf{A}_S\|_2, \|\widetilde{\mathbf{A}}_{S,S'}\|_2) \leq \|\widetilde{\mathbf{A}}_{S,S'}\|_1$$
$$\max(\|\widetilde{\mathbf{A}}_0\|_1, \|\widetilde{\mathbf{A}}_0 + \mathbf{A}_S\|_1, \|\widetilde{\mathbf{A}}_{S,S'}\|_1) \leq \|\widetilde{\mathbf{A}}_{S,S'}\|_1,$$

since $\widetilde{\mathbf{A}}_0$, $\widetilde{\mathbf{A}}_0 + \mathbf{A}_S$ and $\widetilde{\mathbf{A}}_{S,S'}$ are symmetric and by Hölder's inequality for any symmetric matrix $\mathbf{A}$ we have $\|\mathbf{A}\|_2 \leq \sqrt{\|\mathbf{A}\|_1\|\mathbf{A}\|_\infty} = \|\mathbf{A}\|_1$.

We can make sure that the matrices $\boldsymbol{\Theta}_0$ and $\boldsymbol{\Theta}_S$ fall into the sets $\mathcal{U}_s(I_L^*, I_U^*; \theta, \mu)$ and $\mathcal{U}_s(I_L^*, I_U^*; \theta, \mu)$ respectively, and in addition the matrix $\boldsymbol{\Theta}_{S,S'}$ is strictly positive definite if $\theta \leq \mu < \frac{1-C^{-1}}{\|\mathbf{A}_{S,S'}\|_1}$. Thus by assumption the graphs $\boldsymbol{\Theta}_0 \in \mathcal{U}_s(I_L^*, I_L^*; \theta, \mu)$ and $\boldsymbol{\Theta}_S \in \mathcal{U}_s(I_L^*, I_U^*; \theta, \mu)$ for all $S \in \mathbf{C}$. In addition this implies that the matrices $\boldsymbol{\Theta}_0, \boldsymbol{\Theta}_S, \boldsymbol{\Theta}_{S,S'}$ are strictly positive definite for any $S, S' \in \mathbf{C}$.

In the second step we bound the $\chi^2$-divergence from below. It suffices to bound:

$$\left(\frac{\det(\mathbf{\Theta}_S)}{\det(\mathbf{\Theta}_0)}\right)^{n/2}\left(\frac{\det(\mathbf{\Theta}_{S'})}{\det(\mathbf{\Theta}_{S,S'})}\right)^{n/2}$$
$$= \exp\left(\frac{n}{2}\sum_{k=1}^{\infty}\frac{(-\theta)^k}{k}\operatorname{Tr}\left(\widetilde{\mathbf{A}}_{S,S'}^k + \widetilde{\mathbf{A}}_0^k - (\widetilde{\mathbf{A}}_0 + \mathbf{A}_S)^k - (\widetilde{\mathbf{A}}_0 + \mathbf{A}_{S'})^k\right)\right),$$

For any $k \in \mathbb{N}$ we have the following bound:

$$\operatorname{Tr}(\widetilde{\mathbf{A}}_{S,S'}^k + \widetilde{\mathbf{A}}_0^k - (\widetilde{\mathbf{A}}_0 + \mathbf{A}_S)^k - (\widetilde{\mathbf{A}}_0 + \mathbf{A}_{S'})^k) \geq 0.$$

Same as the proof of Proposition C.1 in Neykov et al. (2016), we still consider three cases: (1) $k < 2(d_{G_0}(S,S')+1)$, (2) $k < 4$ and $k \geq 2(d_{G_0}(S,S')+1)$ and (3) $k \geq 4$ and $k \geq 2(d_{G_0}(S,S')+1)$. For $k < 2(d_{G_0}(S,S')+1)$, the above is in fact an equality since if the length of any weighted walk contained in $\widetilde{\mathbf{A}}_0, \mathbf{A}_S, \mathbf{A}_{S'}$ or $\widetilde{\mathbf{A}}_{S,S'}$ with length smaller than $2(d_{G_0}(S,S')+1)$ cannot contain two edges in $S$ and $S'$, and hence will cancel out in the expression. On the other hand, when $k < 4$ and $k \geq 2(d_{G_0}(S,S')+1)$, we have

$$\operatorname{Tr}(\widetilde{\mathbf{A}}_0^2 + \widetilde{\mathbf{A}}_{S,S'}^2 - (\widetilde{\mathbf{A}}_0 + \mathbf{A}_S)^2 - (\widetilde{\mathbf{A}}_0 + \mathbf{A}_{S'})^2) \leq 4|S \cap S'|, \tag{S.14}$$

This is because every length 2 closed walk in (S.14), must have one edge in $S$ and one edge in $S'$.

In the third step we check that:

$$\operatorname{Tr}(\widetilde{\mathbf{A}}_0^2 + \widetilde{\mathbf{A}}_{S,S'}^2 - (\widetilde{\mathbf{A}}_0 + \mathbf{A}_S)^2 - (\widetilde{\mathbf{A}}_0 + \mathbf{A}_{S'})^2) \leq \mathcal{H}_{S,S'}\mathcal{K}_{S,S'}^k, \tag{S.15}$$

for any of the two settings below

**Setting 1:** $\mathcal{H}_{S,S'} = \dfrac{(|\mathcal{V}_{S,S'|S}| \wedge |\mathcal{V}_{S,S'|S'}|)\|\mathbf{A}_S\|_2\|\mathbf{A}_{S'}\|_2}{\|\widetilde{\mathbf{A}}_{S,S'}\|_2^2}$ and $\mathcal{K}_{S,S'} = 2\|\widetilde{\mathbf{A}}_{S,S'}\|_2$;

**Setting 2:** $\mathcal{H}_{S,S'} = \dfrac{|\mathcal{V}_{S,S'|S}||\mathcal{V}_{S,S'|S'}|}{\|\widetilde{\mathbf{A}}_{S,S'}\|_1^2}$  and  $\mathcal{K}_{S,S'} = 2\|\widetilde{\mathbf{A}}_{S,S'}\|_1$.

We begin with showing a similar result as Eq.(C.6), Neykov et al. (2016) but for weighted graphs. Given weighted adjacency matrices $\mathbf{A}_1, \ldots \mathbf{A}_j$ and let $w_{ii}$ as the number of weighted closed walks starting and ending at $i$ and the $\ell$-th edge belongs to the $\ell$-th adjacency matrix for $\ell \in [j]$. We have

$$w_{ii} = A_{ii} \leq \|\mathbf{A}\|_2 \leq \prod_{\ell \in [j]}\|\mathbf{A}_\ell\|_2. \tag{S.16}$$

We aim to prove that

$$\text{Tr}(\widetilde{\mathbf{A}}_0^k + \widetilde{\mathbf{A}}_{S,S'}^k - (\widetilde{\mathbf{A}}_0 + \mathbf{A}_S)^k - (\widetilde{\mathbf{A}}_0 + \mathbf{A}_{S'})^k)$$
$$\leq \left[\frac{2\binom{k}{2}(|\mathcal{V}_{S,S'|S}| \wedge |\mathcal{V}_{S,S'|S'}|)\|\mathbf{A}_S\|_2\|\mathbf{A}_{S'}\|_2}{\|\widetilde{\mathbf{A}}_{S,S'}\|_2^2}\right]\|\widetilde{\mathbf{A}}_{S,S'}\|_2^k. \tag{S.17}$$

Since the trace of a weighted adjacency matrix, sums all weighted closed walks of length $k$ in the graph. It suffices to bound the number of closed walks containing both edges in $S$ and $S'$ since a walk is counted on the LHS of (S.17) if and only if it contains both edges in $S$ and $S'$. Recall that $G_S = (V, E_0 \cup S), G_{S'} = (V, E_0 \cup S')$.

Denote $\mathcal{C}_{S,S'}^{(k)} = \{\text{closed walks } \mathcal{C} \text{ of length } k \text{ on } G(\widetilde{\mathbf{A}}_{S,S'})\}$. Given any length $k$ closed walk $\mathcal{C} = v_1 \to v_2 \to \ldots \to v_k$, let $v_t$ be its $t$-th vertex and let $(v_t, v_{t+1})$ be its $t$-th edge. In the special case of $t = 0$, $(v_{k-1}, v_k)$ is the 0-th edge.

For any vertex $v \in \mathcal{V}_{S,S'}$ and any $1 \leq t_1 \neq t_2 \leq k$, we first count the number of closed walks in the set $\mathcal{C}_k(v, t_1, t_2) = \{\mathcal{C} \in \mathcal{C}_{S,S'}^{(k)} |\ v \text{ is the } t_1\text{-th vertex on } \mathcal{C}, \text{ the } t_1\text{-th edge on } \mathcal{C}$ belongs to $S$ and the $t_2$-th edge on $\mathcal{C}$ belongs to $S'\}$. For a closed walk $\mathcal{C} \in G(\mathbf{A}_{S,S'})$ denote with $|\mathcal{C}|_W = \prod_{e \in \mathcal{C}} \{\mathbf{A}_{S,S'}\}_e$ the weighted length of the closed walk. Following (S.16), we have

$$\sum_{\mathcal{C} \in \mathcal{C}_k(v,t_1,t_2)} |\mathcal{C}|_W \leq \|\mathbf{A}_S\|_2 \|\mathbf{A}_{S'}\|_2 \|\widetilde{\mathbf{A}}_{S,S'}\|_2^{k-2}. \tag{S.18}$$

Similarly, we define $\mathcal{C}'_k(v, t_1, t_2) = \{\mathcal{C} \in \mathcal{C}_{S,S'}^{(k)} |\ v \text{ is the } t_1\text{-th vertex on } \mathcal{C}, \text{ the } (t_1 - 1)\text{-th edge}$ on $\mathcal{C}$ belongs to $S$ and the $t_2$-th edge on $\mathcal{C}$ belongs to $S'\}$. We also have

$$\sum_{\mathcal{C} \in \mathcal{C}'_k(v,t_1,t_2)} |\mathcal{C}|_W \leq \|\mathbf{A}_S\|_2 \|\mathbf{A}_{S'}\|_2 \|\widetilde{\mathbf{A}}_{S,S'}\|_2^{k-2}.$$

We notice that the set comprised of closed walks containing $v$ and edges from both $S$ and $S'$ satisfies

$$\mathcal{C}_k(v) = \{\mathcal{C} \in \mathcal{C}_{S,S'}^{(k)} | v \in \mathcal{C}, \mathcal{C} \cap S \neq \varnothing \text{ and } \mathcal{C} \cap S' \neq \varnothing\} \subseteq \bigcup_{1 \leq t_1 \neq t_2 \leq k} \left(\mathcal{C}_k(v, t_1, t_2) \cup \mathcal{C}'_k(v, t_1, t_2)\right).$$

Therefore, we can control the number of such closed walks by

$$\sum_{\mathcal{C} \in \mathcal{C}_k(v,t_1,t_2) \cup \mathcal{C}'_k(v,t_1,t_2)} |\mathcal{C}|_W \leq 2\binom{k}{2}\|\mathbf{A}_S\|_2 \|\mathbf{A}_{S'}\|_2 \|\widetilde{\mathbf{A}}_{S,S'}\|_2^{k-2}.$$

By the definition of $\mathcal{V}_{S,S'}$, each closed walk $\mathcal{C}$ in $\mathcal{C}_{S,S'}^{(k)}$ containing edges from both $S$ and $S'$, has a vertex $v \in \mathcal{C}$ belonging to $\mathcal{V}_{S,S'|S}$, which is also a vertex of an edge in $S$. Therefore, we have the set of closed walks contains edges both from $S$ and $S'$.

$$\mathcal{C}_k = \{\mathcal{C} \in \mathcal{C}_{S,S'}^{(k)} \,|\, \mathcal{C} \cap S \neq \varnothing \text{ and } \mathcal{C} \cap S' \neq \varnothing\} \subseteq \bigcup_{v \in \mathcal{V}_{S,S'|S}} \mathcal{C}_k(v). \tag{S.19}$$

This implies that $\sum_{\mathcal{C} \in \mathcal{C}_k} |\mathcal{C}|_W \leq |\mathcal{V}_{S,S'} \cap V(S)| \cdot 2\binom{k}{2} \|\mathbf{A}_S\|_2 \|\mathbf{A}_{S'}\|_2 \|\widetilde{\mathbf{A}}_{S,S'}\|_2^{k-2}$. By symmetry we immediately obtain:

$$\sum_{\mathcal{C} \in \mathcal{C}_k} |\mathcal{C}|_W \leq (|\mathcal{V}_{S,S'} \cap V(S)| \wedge |\mathcal{V}_{S,S'} \cap V(S')|) \cdot 2\binom{k}{2} \|\mathbf{A}_S\|_2 \|\mathbf{A}_{S'}\|_2 \|\widetilde{\mathbf{A}}_{S,S'}\|_2^{k-2},$$

which completes the proof of (S.17).

For Setting 2, the proof is similar to Setting 1. We have two vertices in $\mathcal{V}_{S,S'}$ and each vertex can go into $\|\widetilde{\mathbf{A}}_{S,S'}\|_1$ weighted vertices, where the weight corresponds to the type of edge we can go on — $\mu$ or $\theta$. Since two of the vertices are on fixed positions we are left with $\|\widetilde{\mathbf{A}}_{S,S'}\|_1^{k-2}$ vertices to complete the path. Since there are at most $k(k-1) \leq 2^k$ possibilities for the positions of two vertices, and there are at least the 2 edges which are coming out of the two fixed vertices have weight $\theta$, we complete the proof for Setting 2.

Finally we complete the proof by showing that if $\{\mathcal{H}_{S,S'}, \mathcal{K}_{S,S'}\}_{S,S' \in \mathbf{C}}$ satisfies either Setting 1 or 2, we can bound the $\chi^2$-divergence as the proposition. In particular, we have

$$\sum_{k=1}^{\infty} \theta^k \operatorname{Tr}\left(\widetilde{\mathbf{A}}_{S,S'}^k + \widetilde{\mathbf{A}}_0^k - (\widetilde{\mathbf{A}}_0 + \mathbf{A}_S)^k - (\widetilde{\mathbf{A}}_0 + \mathbf{A}_{S'})^k\right)/k$$
$$\leq 2|S \cap S'|\theta^2 + \sum_{2|k,\ k \geq 2(d_{G_0}(S,S') \vee 1+1)} \mathcal{H}_{S,S'} k^{-1}(\mathcal{K}_{S,S'})^k$$
$$\leq 2|S \cap S'|\theta^2 + \frac{2\mathcal{H}_{S,S'}(\mathcal{K}_{S,S'}\theta)^{2(d_{G_0}(S,S') \vee 1+1)}}{2(d_{G_0}(S,S') \vee 1+1)}.$$

We therefore complete the proof.

### S.2.5 Auxiliary Results on Graph Spectral

In this section, we prove the following lemma of the operator norm of hollow graph.

**Lemma S.6** (Sparse Graph Lemma). *If we have an $R$-hollow graph $G$ with maximum degree $d_{\max}(G)$ and adjacency matrix $\mathbf{A}$ we have:*

$$\|\mathbf{A}\|_2 \leq 2R\sqrt{d_{\max}(G)}.$$

*Proof.* Using a result of Nash-Williams (1964) which follows by earlier results of Nash-Williams (1961), we know that there exists a decomposition of any sparse graph in at most $R$ forests. Since the graph is of degree at most $d_{\max}(G)$ each of the forests is also of at most $d_{\max}(G)$ degree. We will now argue that for any forest graph $F$ with bounded by $d_{\max}(F)$ degree the largest eigenvalue of its adjacency matrix is bounded by $2\sqrt{d_{\max}(F)}$.

First recall that for two non-intersecting graphs on the same vertex set $G'$ and $G''$ we have $\|\mathbf{A}_{G'\cup G''}\|_2 \leq \|\mathbf{A}_{G'}\| + \|\mathbf{A}_{G''}\|_2$

To show the bound we proceed by decomposing the graph $F$, by greedily constructing star graphs. To describe the procedure it suffices to consider the case where $F$ is a tree since it can be applied to all sub-trees of a forest $F$ in parallel. Starting from a root of $F$ split out the root and all of its children. Recursively moving in a breadth first search manner add edges by the following procedure: do not take any vertex who has already been taken into account, and for any vertex who has not been taken into account take all edges from him to his children. Clearly two passes are enough to partition the forest into two graphs of non-intersecting star-graphs. It is simple to check that a star graph of maximum degree $d_{\max}$ has $\sqrt{d_{\max}}$ as the operator norm of its adjacency matrix. This implies that

$$\|\mathbf{A}_F\|_2 \leq 2\sqrt{d_{\max}(F)},$$

as claimed. Putting everything together we have:

$$\|\mathbf{A}\|_2 \leq 2R\sqrt{d_{\max}(G)}.$$

□

**Remark S.2.1.** Lemma S.6 is developed with the regime $R = o(\sqrt{d_{\max}(G)})$ in mind. Otherwise it is clear that always:

$$\|\mathbf{A}\|_2 \leq d_{\max}(G).$$

## S.3 Power Analysis of the Tests

In this section, we discuss the power analysis for the test $\psi_\alpha$ obtained from Algorithm 2. Under the alternative $H_1 : \mathcal{P}(\Theta) = 1$, we define the parameter space

$$\mathcal{G}_1(\theta; \mathcal{P}) = \left\{ \Theta \in \mathcal{U}_s \,\middle|\, \mathcal{P}(G(\Theta)) = 1 \text{ and } \max_{E' \subseteq E(\Theta), \mathcal{P}(E')=1} \min_{e \in E'} |\Theta_e| \geq \theta \right\}. \tag{S.20}$$

If $\Theta \in \mathcal{G}_1(\theta; \mathcal{P})$, by (S.20), its induced graph $G(\Theta) = (V, E(\Theta))$ must has a sub-edge set $E_0 \subseteq E(\Theta)$ such that $\mathcal{P}(E_0) = 1$ and the minimal signal strength on $E_0$ is larger than $\theta$.

**Theorem S.7** (Power analysis). Suppose $\Theta \in \mathcal{U}_s$ and (3.5) is satisfied. Under the alternative hypothesis $H_1 : \mathcal{P}(G) = 1$, there exists a positive constant $C$ such that for any $\alpha \in (0, 1)$ with $1/(\alpha(d \vee n)) = o(1)$,

$$\lim_{n \to \infty} \inf_{\Theta \in \mathcal{G}_1(C\tau_n; \mathcal{P})} \mathbb{P}_\Theta(\psi_\alpha = 1) = 1, \tag{S.21}$$

where $\tau_n = \sqrt{\log d/n}$.

Notice that the parameter space $\mathcal{G}_1(\theta; \mathcal{P})$ defined in (S.20) is larger than the parameter space such that $\boldsymbol{\Theta}$ has the minimal signal strengths $\theta$ on its support $E(\boldsymbol{\Theta})$. Namely, we have

$$\mathcal{G}_1(\theta; \mathcal{P}) \supseteq \mathcal{G}_1'(\theta; \mathcal{P}) = \Big\{\boldsymbol{\Theta} \in \mathcal{U}_s \,|\, \mathcal{P}(G(\boldsymbol{\Theta})) = 1 \text{ and } \min_{e \in E(\boldsymbol{\Theta})} |\boldsymbol{\Theta}_e| \geq \theta \Big\}. \qquad (\text{S.22})$$

The parameter space $\mathcal{G}_1'(\theta)$ is usually considered for the power analysis of the global hypothesis tests (Han and Liu, 2014) or high dimensional two sample tests (Cai et al., 2013, 2014). For example, Han and Liu (2014) consider the null hypothesis $H_0 : \boldsymbol{\Sigma} = \mathbf{I}_d$ and show their test is powerful if $\min_{j \neq k} |\boldsymbol{\Sigma}_{jk}| > C\sqrt{\log d/n}$ for some constant $C > 0$. Theorem S.7 demonstrates that for the power analysis of graph property test, it suffices to impose the minimal signal strength condition on a subgraph $E'$ of the true support $E(\boldsymbol{\Theta})$ if $\mathcal{P}(E') = 1$. The following proposition gives a concrete characterization of $\mathcal{G}_1(\theta)$.

**Proposition S.8.** Given the graph $G(\boldsymbol{\Theta}) = (V, E(\boldsymbol{\Theta}))$, we set the weights on the edge $(j, k)$ as $|\boldsymbol{\Theta}_{jk}|$ for any $j, k \in V$. We order the weights $|\boldsymbol{\Theta}_e|$ for $e \in V \times V$ as $|\boldsymbol{\Theta}_{e_{[1]}}| \geq \ldots \geq |\boldsymbol{\Theta}_{e_{[d(d-1)/2]}}|$. Let $t^* = \arg\min\{t \,|\, \mathcal{P}(E_{[t]}) = \mathcal{P}(\{e_{[1]}, \ldots, e_{[t]}\}) = 1\}$. If the property $\mathcal{P}$ is monotone, we have

$$\mathcal{G}_1(\theta; \mathcal{P}) = \big\{\boldsymbol{\Theta} \in \mathcal{U}_s \,\big|\, |\boldsymbol{\Theta}_{e_{[t^*]}}| \geq \theta\big\}.$$

*Proof.* By the definition in (S.20), since $\mathcal{P}(E_{[t^*]}) = 1$, we have

$$\mathcal{G}_1(\theta; \mathcal{P}) \supseteq \big\{\boldsymbol{\Theta} \in \mathcal{U}_s \,\big|\, |\boldsymbol{\Theta}_{e_{[t^*]}}| \geq \theta\big\}.$$

It suffices to prove the other direction. If $\boldsymbol{\Theta} \in \mathcal{G}_1(\theta; \mathcal{P})$, there exists $E' \subseteq E(\boldsymbol{\Theta})$ with $\mathcal{P}(E') = 1$ such that $\min_{e \in E'} |\boldsymbol{\Theta}_e| \geq \theta$. Since $\mathcal{P}$ is monotone, we have $\mathcal{P}(E_{[t^*]} \cup E') = 1$. As $\min_{e \in E'} |\boldsymbol{\Theta}_e| \geq \theta$ and $E_{[t^*]}$ contains the top $t^*$ largest weight edges, we have $|\boldsymbol{\Theta}_{e_{[t^*]}}| \geq \theta$. This completes the proof of the proposition. $\square$

**Remark S.3.1.** Proposition S.8 implies that we can add the edges in $E(\boldsymbol{\Theta})$ from the largest to the smallest until we stop at the $t^*$-th step when $\mathcal{P}(E_{[t^*]}) = 1$. For example, for the connectivity property $\mathcal{P}_{\text{Conn}(1)}$, we greedily add edges in order until the graph is connected. Comparing to the procedure of Kruskal's algorithm (Kruskal, 1956), let $E_{\text{MSF}}$ be the edge set of the maximum spanning tree of $G(\boldsymbol{\Theta})$ with $\boldsymbol{\Theta}$ as the weights, we can have

$$\mathcal{G}_1(\theta; \mathcal{P}_{\text{Conn},-1}) = \Big\{\boldsymbol{\Theta} \in \mathcal{U}_s \,\Big|\, \min_{e \in E_{\text{MSF}}} |\boldsymbol{\Theta}_e| \geq \theta \Big\}.$$

This coincides with the alternative set for the connectivity test in Neykov et al. (2016). For the maximum degree property $\mathcal{P}_{\text{Deg},k}$, Neykov et al. (2016) consider the alternative set $\mathcal{G}_1'(\theta)$ in (S.22). Proposition S.8 implies that our test is asymptotically powerful in a larger parameter space. In specific, we order the weights of the neighbor of node $j$ as $|\boldsymbol{\Theta}_{j(1)}| \geq \ldots \geq |\boldsymbol{\Theta}_{j(d)}|$. For for any $k \in [1, d]$, we have

$$\mathcal{G}_1(\theta; \mathcal{P}_{\text{Deg},k}) = \Big\{\boldsymbol{\Theta} \in \mathcal{U}_s \,\Big|\, \mathcal{P}_{\text{Deg},k}(G(\boldsymbol{\Theta})) = 1 \text{ and } \max_{j \in V} |\boldsymbol{\Theta}_{j(k)}| \geq \theta \Big\}.$$

From Theorems 4.1 and S.7, we show that the test $\psi_\alpha$ obtained from Algorithm 1 is optimal for the four examples we discuss in Section 3.2: $\mathcal{P}_{\mathrm{Conn},-k}$, $\mathcal{P}_{\mathrm{Deg},k}$ and $\mathcal{P}_{\mathrm{Sig},-k}$. To achieve this goal, we define the risk of a test $\psi$ as

$$R_\theta(\psi, \mathcal{P}) = \sup_{\boldsymbol{\Theta} \in \mathcal{G}_0} \mathbb{P}_{\boldsymbol{\Theta}}(\psi = 1) + \sup_{\boldsymbol{\Theta} \in \mathcal{G}_1(\theta; \mathcal{P})} \mathbb{P}_{\boldsymbol{\Theta}}(\psi = 0), \qquad (S.23)$$

where the parameter spaces $\mathcal{G}_0$ and $\mathcal{G}_1(\theta; \mathcal{P})$ are defined in (4.3) and (S.20) respectively. We say that the hypotheses $H_0: \boldsymbol{\Theta} \in \mathcal{G}_0$ v.s. $H_{1,\theta}: \boldsymbol{\Theta} \in \mathcal{G}_1(\theta; \mathcal{P})$ are asymptotically separated by a test $\psi$ if $\lim_{n \to \infty} R_\theta(\psi, \mathcal{P}) = 0$. On the other hand, $H_0$ and $H_{1,\theta}$ are *asymptotically inseparable* if

$$\liminf_{n \to \infty} \inf_\psi R_\theta(\psi, \mathcal{P}) = 1.$$

From Theorems 4.1 and S.7, the test $\psi_\alpha$ obtained from Algorithm 1 can asymptotically separate $H_0$ and $H_{1,\theta}$ for any monotone property $\mathcal{P}$ if the signal strength is strong enough. In specific, we have the following corollary directly from (4.2) and (S.21).

**Corollary S.9.** Under the same conditions as Theorem S.7, for any monotone property $\mathcal{P}$, we choose the level of significance $\alpha_n \in (0,1)$ such that $1/(\alpha_n(d \vee n)) = o(1)$ and $\alpha_n = o(1)$. The hypotheses $H_0: \boldsymbol{\Theta} \in \mathcal{G}_0$ v.s. $H_{1,\theta}: \boldsymbol{\Theta} \in \mathcal{G}_1(\theta; \mathcal{P})$ can be asymptotically separated by the test $\psi_{\alpha_n}$ when $\theta \geq C_1 \sqrt{\log d/n}$ for some universal constant $C_1 > 0$.

On the other hand, the lower bound on the signal strength developed in Neykov et al. (2016) shows that the signal strength in Corollary S.9 is optimal for the properties including connectivity, maximum degrees and cyclic graph. In specific, they show that there exists a constant $C_2 > 0$ such that if $\theta \leq C_2 \sqrt{\log d/n}$, and $\log d/n = o(1)$, for the following properties: (1) $\mathcal{P}_{\mathrm{Conn},-k}$ for $k \in [1, d)$; (2) $\mathcal{P}_{\mathrm{Cycle}}$; (3) $\mathcal{P}_{\mathrm{Deg},k}$ for $k \in [1, \sqrt{d})$ and $k \log d/n = o(1)$, the hypotheses $H_0: \boldsymbol{\Theta} \in \mathcal{G}_0$ v.s. $H_{1,\theta}: \boldsymbol{\Theta} \in \mathcal{G}_1(\theta; \mathcal{P})$ are asymptotically inseparable. For $\mathcal{P}_{\mathrm{Sig},-k}$ for $k \in [0, d)$, we can establish the same lower bound. See Section S.3.1.1 in the appendix for details.

### S.3.1 Proofs of Lower Bound of Tests

The lower bounds for the three properties: $\mathcal{P}_{\mathrm{Conn},-k}$, $\mathcal{P}_{\mathrm{Cycle}}$ and $\mathcal{P}_{\mathrm{Deg},-k}$ are proved in Neykov et al. (2016). Following a similar procedure, we will prove the lower bound for $\mathcal{P}_{\mathrm{Sig},-k}$.

For the self-consistency of the paper, we introduce the Le Cam's lemma (LeCam, 1973) first. A standard application of Le Cam's lemma is for the lower bound on when a null hypothesis with a single parameter can be separated from the alternative. Under that case, the Le Cam's lemma selects a set of parameters in the alternatives such that they deviate from the null parameter far away enough. In specific, let $\boldsymbol{\Theta}_0 \in \mathcal{G}_0$ and $\boldsymbol{\Theta}_{(1)}, \ldots, \boldsymbol{\Theta}_{(K)} \in \mathcal{G}_1(\theta)$. Let $\mathbb{P}^n_{\boldsymbol{\Theta}}$ be the product measure of the i.i.d. samples $\{\boldsymbol{X}_1, \ldots, \boldsymbol{X}_n\} \sim N(0, \boldsymbol{\Theta}^{-1})$. We denote the mixture measure $\overline{\mathbb{P}}^n = K^{-1} \sum_{j=1}^K \mathbb{P}_{\boldsymbol{\Theta}_{(j)}}$. Le Cam's lemma gives the following lower bound of the risk (after adapting our notations).

**Lemma S.10.** Let $R_\theta(\psi, \mathcal{P})$ be the risk defined in (S.23), we have

$$\inf_\psi R_\theta(\psi, \mathcal{P}) \geq \mathbb{P}^n_{\Theta_0}(\psi = 1) + \overline{\mathbb{P}}^n(\psi = 0) \geq 1 - \frac{1}{2}\bigg(\int \Big(\frac{d\overline{\mathbb{P}}^n}{d\mathbb{P}^n_{\Theta_0}}\Big)^2 d\mathbb{P}^n_{\Theta_0} - 1\bigg)^{1/2}. \quad (S.24)$$

However, under the property testing problems, we have both composite null hypothesis and composite alternatives. When the null hypothesis contains richer parameter space that the alternative, we should instead choose a parameter from alternatives $\Theta_A \in \mathcal{G}_1(\theta)$. On the other hand, we collate $\Theta_{(1)}, \ldots, \Theta_{(K)} \in \mathcal{G}_0$. Le Cam's lemma also gives the following lower bound of the risk:

$$\inf_\psi R_\theta(\psi, \mathcal{P}) \geq \mathbb{P}^n_{\Theta_A}(\psi = 1) + \overline{\mathbb{P}}^n(\psi = 0) \geq 1 - \frac{1}{2}\bigg(\int \Big(\frac{d\overline{\mathbb{P}}^n}{d\mathbb{P}^n_{\Theta_A}}\Big)^2 d\mathbb{P}^n_{\Theta_A} - 1\bigg)^{1/2}. \quad (S.25)$$

#### S.3.1.1 Proof of Lower Bound on Singletons

For the property of isolated nodes, we consider two cases: (1) $k \in [0, d/2)$ and (2) $k \in (d/2, d)$. As an atom to construct the matrices in (S.25), we consider the precision matrix

$$\Theta_{\text{Chain}(d)} = \begin{pmatrix} 1 & \theta & & & 0 \\ \theta & 1 & \theta & & \\ & \theta & \ddots & \ddots & \\ & & \ddots & \ddots & \theta \\ 0 & & & \theta & 1 \end{pmatrix}. \quad (S.26)$$

The precision matrix above induces a chain with $d$ nodes and $d-1$ edges.

**Case 1:** $k \in [0, d/2)$. For this case, we apply Le Cam's Lemma in (S.25) by selecting a $\Theta_A \in \mathcal{G}_1(\theta)$ and $\{\Theta_{(1)}, \ldots, \Theta_{(K)}\} \in \mathcal{G}_0$. We set

$$\Theta_A = \text{diag}\Big(\underbrace{\Theta_{\text{Chain}(2)}, \ldots, \Theta_{\text{Chain}(2)}}_{\lfloor (d-k)/2 \rfloor \text{ matrices}}, \Theta_{\text{Chain}(d-k-2\lfloor (d-k)/2 \rfloor)}, \mathbf{I}_k\Big) \in \mathbb{R}^{d \times d},$$

where $d - k - 2\lfloor (d-k)/2 \rfloor) = 2$ if $d - k$ is even and 3 otherwise. We also denote $\Theta_{(j)} = \Theta_0 - \theta \mathbf{E}_{2j, 2j+1}$ for $1 \leq j \leq \lfloor (d-k)/2 \rfloor - 1$.

Denote $f_j$ for $1 \leq j \leq \lfloor (d-k)/2 \rfloor - 1$ as the density of $\mathbb{P}^n_{\Theta_{(j)}}$ and $f_0$ for $\mathbb{P}^n_{\Theta_0}$. Similar to the previous proofs, we first bound the square term

$$\int \frac{f_j^2}{f_0} d\mu = \bigg(\frac{\det(\Theta_{(j)})^2}{\det(\Theta_A)\det(2\Theta_{(j)} - \Theta_A)}\bigg)^{n/2} = \bigg(\frac{(1-\theta^2)^{2\lfloor (d-k)/2 \rfloor - 4}}{(1-\theta^2)^{2\lfloor (d-k)/2 \rfloor - 2}}\bigg)^{n/2} = \frac{1}{(1-\theta^2)^n},$$

and the summation

$$\frac{1}{\lfloor (d-k)/2\rfloor^2} \sum_{j=1}^{\lfloor (d-k)/2\rfloor} \int \frac{f_j^2}{f_0} d\mu \leq (\lfloor (d-k)/2\rfloor)^{-1} \exp(1+n\theta^2) \leq \exp(-\log d + n\theta^2),$$

where the last inequality is due to $k \leq d/2$. Therefore, if $\theta \leq C_2 \sqrt{\log d/n}$ for some sufficiently small $C_2 > 0$,

$$\frac{1}{\lfloor (d-k)/2\rfloor^2} \sum_{j=1}^{\lfloor (d-k)/2\rfloor} \int \frac{f_j^2}{f_0} d\mu \to 0.$$

For the cross term, we have

$$\int \frac{f_j f_k}{f_0} d\mu = \left(\frac{\det(\boldsymbol{\Theta}_{(j)})\det(\boldsymbol{\Theta}_{(k)})}{\det(\boldsymbol{\Theta}_A)\det(\boldsymbol{\Theta}_{(j)} + \boldsymbol{\Theta}_{(k)} - \boldsymbol{\Theta}_A)}\right)^{n/2}$$

$$= \left[\frac{(1-\theta^2)^{\lfloor (d-k)/2\rfloor -2}(1-\theta^2)^{\lfloor (d-k)/2\rfloor -2}}{(1-\theta^2)^{\lfloor (d-k)/2\rfloor -1}(1-\theta^2)^{\lfloor (d-k)/2\rfloor -3}}\right]^{n/2} = 1.$$

Therefore, we conclude that

$$\Delta = \frac{1}{\lfloor (d-k)/2\rfloor^2} \sum_{j=1}^{\lfloor (d-k)/2\rfloor} \int \frac{f_j^2}{f_0} d\mu + \frac{2}{\lfloor (d-k)/2\rfloor^2} \sum_{1\leq j<k\leq \lfloor (d-k)/2\rfloor^2} \left[\int \frac{f_j f_k}{f_0} d\mu - 1\right] \to 0.$$

**Case 2:** $k \in [d/2, d)$ We apply the Le Cam's Lemma in (S.24) by choosing a $\boldsymbol{\Theta}_0 \in \mathcal{G}_0$ and $\{\boldsymbol{\Theta}_{(1)}, \ldots, \boldsymbol{\Theta}_{(K)}\} \in \mathcal{G}_1(\theta)$.

We set

$$\boldsymbol{\Theta}_0 = \mathrm{diag}\Big(\mathbf{I}_{k+1}, \underbrace{\boldsymbol{\Theta}_{\mathrm{Chain}(2)}, \ldots, \boldsymbol{\Theta}_{\mathrm{Chain}(2)}}_{\lfloor (d-k-1)/2\rfloor \text{ matrices}}, \boldsymbol{\Theta}_{\mathrm{Chain}(d-k-1-2\lfloor (d-k-1)/2\rfloor)}\Big) \in \mathbb{R}^{d\times d},$$

and $\boldsymbol{\Theta}_{(j)} = \boldsymbol{\Theta}_0 + \theta \mathbf{E}_{2j,2j+1}$ for $1 \leq j \leq \lfloor (k+1)/2\rfloor$.

Using the same notations, we again first bound the square term

$$\int \frac{f_j^2}{f_0} d\mu = \left(\frac{\det(\boldsymbol{\Theta}_{(j)})^2}{\det(\boldsymbol{\Theta}_0)\det(2\boldsymbol{\Theta}_{(j)} - \boldsymbol{\Theta}_0)}\right)^{n/2}$$

$$= \left(\frac{(1-\theta^2)^{2\lfloor (d-k-1)/2\rfloor +2}}{(1-\theta^2)^{\lfloor (d-k-1)/2\rfloor}(1-\theta^2)^{\lfloor (d-k-1)/2\rfloor +1}}\right)^{n/2} = \frac{1}{(1-\theta^2)^{n/2}},$$

and since $k \in (d/2, d)$, we have

$$\frac{1}{\lfloor (k+1)/2 \rfloor^2} \sum_{j=1}^{\lfloor (k+1)/2 \rfloor} \int \frac{f_j^2}{f_0} d\mu \leq (\lfloor (k+1)/2 \rfloor)^{-1} \exp(1 + n\theta^2/2) \leq \exp(-\log d + n\theta^2/2).$$

If $\theta \leq C_2 \sqrt{\log d/n}$ for some sufficiently small $C_2 > 0$,

$$\frac{1}{\lfloor (k+1)/2 \rfloor^2} \sum_{j=1}^{\lfloor (k+1)/2 \rfloor} \int \frac{f_j^2}{f_0} d\mu \to 0.$$

We can calculate the cross term

$$\int \frac{f_j f_k}{f_0} d\mu = \left( \frac{\det(\Theta_{(j)}) \det(\Theta_{(k)})}{\det(\Theta_A) \det(\Theta_{(j)} + \Theta_{(k)} - \Theta_A)} \right)^{n/2}$$

$$= \left[ \frac{(1-\theta^2)^{\lfloor (d-k-1)/2 \rfloor + 1}(1-\theta^2)^{\lfloor (d-k-1)/2 \rfloor + 1}}{(1-\theta^2)^{\lfloor (d-k-1)/2 \rfloor}(1-\theta^2)^{\lfloor (d-k-1)/2 \rfloor + 2}} \right]^{n/2} = 1.$$

Therefore, we have

$$\Delta = \frac{1}{\lfloor (d-k)/2 \rfloor^2} \sum_{j=1}^{\lfloor (d-k)/2 \rfloor} \int \frac{f_j^2}{f_0} d\mu + \frac{2}{\lfloor (d-k)/2 \rfloor^2} \sum_{1 \leq j < k \leq \lfloor (d-k)/2 \rfloor^2} \left[ \int \frac{f_j f_k}{f_0} d\mu - 1 \right] \to 0.$$

### S.3.1.2 Auxiliary Results on Matrix Determinant

**Lemma S.11.** Let

$$r_1 = \frac{1 + \sqrt{1 - 4\theta^2}}{2} \text{ and } r_2 = \frac{1 - \sqrt{1 - 4\theta^2}}{2}. \tag{S.27}$$

For the matrices $\Theta_{\text{Chain}(d)}$ defined in (S.26) and $\mathbf{E}_{j,j+1} = \mathbf{e}_j \mathbf{e}_{j+1}^T + \mathbf{e}_{j+1} \mathbf{e}_j^T$, we have for all $1 \leq j \leq d-1$,

$$\det(\Theta_{\text{Chain}(d)}) = \det(\Theta_{\text{Chain}(d)} - 2\theta \mathbf{E}_{j,j+1}) = \frac{r_1^{d+1} - r_2^{d+1}}{r_1 - r_2}.$$

*Proof.* For a tridiagonal matrix $M_k$:

$$M_k := \begin{pmatrix} a_1 & b_1 & & & 0 \\ c_1 & a_2 & b_2 & & \\ & c_2 & \ddots & \ddots & \\ & & \ddots & \ddots & b_{k-1} \\ 0 & & & c_{k-1} & a_k \end{pmatrix},$$

we expand the determinant on the first row:

$$\det(M_k) = a_k \det(M_{k-1}) - c_{k-1}b_{k-1}\det(M_{k-2}).$$

Plugging in $c_i = b_i = \theta, i \in [k-1]$ and $a_i = 1, i \in [k]$ we obtain:

$$\det(\mathbf{\Theta}_{\text{Chain}(i)}) = \det(\mathbf{\Theta}_{\text{Chain}(i-1)}) - \theta^2 \det(\mathbf{\Theta}_{\text{Chain}(i-2)}). \tag{S.28}$$

Writing out the characteristic polynomial of this recursive relationship, it is easy to check that it has two roots: $r_1$ and $r_2$. Finally checking the initial conditions $\det \mathbf{\Theta}_{\text{Chain}(1)} = 1$, $\det \mathbf{\Theta}_{\text{Chain}(2)} = 1 - \theta^2$ yields the desired formula.

For $\det(\mathbf{\Theta}_{\text{Chain}(d)} - 2\theta \mathbf{E}_{j,j+1})$ observe that $c_j = b_j = -\theta$, which doesn't change the recursive relationship (S.28), and hence the determinant is the same. □

## S.4 Proof of Lemma S.1

We first consider a special case that $\mathcal{I}(\mathcal{T}_\mu(\mathbf{\Theta})) = \mathcal{I}(\varnothing)$. Since $\mathcal{I}$ is monotone, according to Algorithm 1, $\widehat{\mathcal{I}}_L \geq \mathcal{I}(\mathcal{T}_\mu(\mathbf{\Theta})) = \mathcal{I}(\varnothing)$ almost surely and (S.6) is trivially true.

Therefore, in the following of the proof, we consider the case $\mathcal{I}(\mathcal{T}_\mu(\mathbf{\Theta})) > \mathcal{I}(\varnothing)$. This implies that there exists a non-empty edge set $E_0'$ such that

$$E_0' \subseteq E(\mathcal{T}_\mu(\mathbf{\Theta})), \mathcal{I}(E_0') = \mathcal{I}(\mathcal{T}_\mu(\mathbf{\Theta})) \text{ and } \min_{e \in E_0'} |\mathbf{\Theta}_e| > C\sqrt{\log d/n}, \tag{S.29}$$

where the constant $C$ is determined later. By $\mathcal{I}(E_0') > \mathcal{I}(\varnothing)$, similar to (S.4), we claim that $\mathcal{I}(E_0' \cap \mathcal{C}_\mathcal{I}(\varnothing)) = \mathcal{I}(E_0')$. To prove this claim, we find a subgraph $E_0'' \subset E_0'$ such that $\mathcal{I}(E_0'') = \mathcal{I}(E_0')$ and for any $\widetilde{E} \subset E_0''$, we have $\mathcal{I}(\widetilde{E}) < \mathcal{I}(E_0'')$. Such graph $E_0''$ can be constructed by deleting edges from $E_0'$ and making the invariant equal to $\mathcal{I}(E_0')$ until it is impossible to further delete edges without reducing the value of the invariant. By Definition 3.1, we have $E_0'' \subseteq \mathcal{C}_\mathcal{I}(\varnothing)$ which implies $E_0'' \subseteq E_0' \cap \mathcal{C}_\mathcal{I}(\varnothing)$. As $\mathcal{I}(E_0'') = \mathcal{I}(E_0')$, by monotone property, we prove the claim that $\mathcal{I}(E_0' \cap \mathcal{C}_\mathcal{I}(\varnothing)) = \mathcal{I}(E_0')$. Consider the following event

$$\mathcal{E}_1 = \left\{ \min_{e \in E_0' \cap \mathcal{C}_\mathcal{I}(\varnothing)} \sqrt{n}|\widehat{\mathbf{\Theta}}_e^d| > c(\alpha, \mathcal{C}_\mathcal{I}(\varnothing)) \right\}.$$

According to Algorithm 1, the rejected set in the first iteration is

$$E_1 = \{e \in \mathcal{C}_\mathcal{I}(\varnothing) : \sqrt{n}|\widehat{\mathbf{\Theta}}_e^d| > c(\alpha, \mathcal{C}_\mathcal{I}(\varnothing))\}.$$

Under the event $\mathcal{E}_1$, we have $E_0' \cap \mathcal{C}_\mathcal{I}(\varnothing) \subseteq E_1$ and since $\mathcal{I}(E_0' \cap \mathcal{C}_\mathcal{I}(\varnothing)) = \mathcal{I}(E_0')$, we have $\mathcal{I}(E_1) = \mathcal{I}(E_0')$, which makes $\widehat{\mathcal{I}}_L \geq \mathcal{I}(E_0') = \mathcal{I}(\mathcal{T}_\mu(\mathbf{\Theta}))$. Therefore, we have

$$\mathbb{P}(\mathcal{I}(\mathcal{T}_\mu(\mathbf{\Theta})) \leq \widehat{\mathcal{I}}_L) \geq \mathbb{P}(\mathcal{E}_1). \tag{S.30}$$

We will bound $\mathbb{P}(\mathcal{E}_1)$ next. We have

$$\mathbb{P}(\mathcal{E}_1) \geq \mathbb{P}\Big( \min_{e \in E_0'} |\boldsymbol{\Theta}_e| > 2c(\alpha, \mathcal{C}_{\mathcal{I}}(\varnothing)) \text{ and } \max_{e \in V \times V} |\widehat{\boldsymbol{\Theta}}_e^d - \boldsymbol{\Theta}_e| \leq c(\alpha, \mathcal{C}_{\mathcal{I}}(\varnothing))\Big). \tag{S.31}$$

We next estimate the rate of $c(\alpha, \mathcal{C}_{\mathcal{I}}(\varnothing))$. By Lemma E.4. of Neykov et al. (2016), we have with probability $1 - 1/d^3$, there exists a constant $C_1$,

$$\max_{j,k \in [d]} |\widehat{\boldsymbol{\Theta}}_{jk}^d - \boldsymbol{\Theta}_{jk} + \boldsymbol{\Theta}_j^T(\widehat{\boldsymbol{\Sigma}} \boldsymbol{\Theta}_k - \mathbf{e}_k)| \leq C_1 \frac{s \log d}{n}. \tag{S.32}$$

Notice here we change the tail probability to $1 - 1/d^3$ from $1 - 1/d$ in Lemma E.4. of Neykov et al. (2016), which can be easily done by changing the value of $C_1$ because the dimension $d$ only appears in the logarithmic term in (S.32). Similar trick also apply to the tail probabilities below. Since $\boldsymbol{\Theta} \in \mathcal{U}_s$ and $\|\boldsymbol{\Theta}\|_2 \leq \rho$, we have $\|\boldsymbol{\Theta}_j^T \boldsymbol{X} \boldsymbol{X}^T \boldsymbol{\Theta}_k\|_{\psi_1} \leq C_2 \rho^2$. By the maximal inequality (see Lemma 2.2.2 in van der Vaart and Wellner (1996)), we have for some constant $C_3 > 0$

$$\mathbb{P}\Big( \max_{j,k \in [d]} |\boldsymbol{\Theta}_j^T(\widehat{\boldsymbol{\Sigma}} \boldsymbol{\Theta}_k - \mathbf{e}_k)| > C_3 \rho^2 \sqrt{\frac{\log d}{n}} \Big)$$
$$= \mathbb{P}\Big( \max_{j,k \in [d]} \Big| \frac{1}{n} \sum_{i=1}^n (\boldsymbol{\Theta}_j^T \boldsymbol{X}_i \boldsymbol{X}_i^T \boldsymbol{\Theta}_k - \mathbb{E}[\boldsymbol{\Theta}_j^T \boldsymbol{X}_i \boldsymbol{X}_i^T \boldsymbol{\Theta}_k])\Big| > C_3 \rho^2 \sqrt{\frac{\log d}{n}} \Big) \leq 1/d^3.$$

Let $C_0 = C_1 + C_3$ and we have

$$\mathbb{P}\Big( \max_{j,k \in [d]} |\widehat{\boldsymbol{\Theta}}_{jk}^d - \boldsymbol{\Theta}_{jk}| > C_0 \sqrt{\frac{\log d}{n}} \Big) < \frac{2}{d^3}. \tag{S.33}$$

Applying (3.5) and (S.33), for any fixed $\alpha \in (0,1)$ and sufficiently large $d, n$ such that $\mathbb{P}(T_{V \times V} > c(\alpha, V \times V)) > \alpha$ and $2/d^3 \leq \alpha$, we have

$$c(\alpha, \mathcal{C}_{\mathcal{I}}(\varnothing)) \leq c(\alpha, V \times V) \leq C_0 \sqrt{\frac{\log d}{n}}.$$

We set the constant $C = 2C_0$ in $\mu \geq C\sqrt{\log d/n}$ and it follows from (S.29) and (S.33) that

$$\min_{e \in E_0' \cap \mathcal{C}_{\mathcal{I}}(\varnothing)} |\boldsymbol{\Theta}_e| > 2C_0 \sqrt{\frac{\log d}{n}} > 2c(\alpha, \mathcal{C}_{\mathcal{I}}(\varnothing)) \text{ and } \mathbb{P}\Big( \max_{e \in V \times V} |\widehat{\boldsymbol{\Theta}}_e^d - \boldsymbol{\Theta}_e| \leq c(\alpha, \mathcal{C}_{\mathcal{I}}(\varnothing))\Big) > 1 - 2/d^3.$$

Combining the above inequalities with (S.30) and (S.31), we have

$$\mathbb{P}(\mathcal{I}(\mathcal{T}_\mu(\boldsymbol{\Theta})) \leq \widehat{\mathcal{I}}_L) \geq \mathbb{P}(\mathcal{E}_1) > 1 - 2/d^3.$$

Since the right hand side of the above inequality is universal for any $\boldsymbol{\Theta} \in \mathcal{U}_{\mathcal{I}}(I_L^*, I_U^*; \theta)$, we

## S.5 Quantile Estimation for Graphical Models

In this section, we first provide concrete examples satisfying (3.5) under Gaussian graphical model. We then generalize our method to non-Gaussian graphical models.

Let the sample covariance be $\widehat{\Sigma} = n^{-1} \sum_{i=1}^{n} X_i X_i^T$. We use the CLIME estimator proposed by Cai et al. (2011) to estimate $\Theta$. It can in fact be generalized to other estimators and we will discuss the extension in Section 5. In order to achieve asymptotic normality from $\widehat{\Theta}$, Neykov et al. (2015) propose the following debiased estimator

$$\widehat{\Theta}_{jk}^d = \widehat{\Theta}_{jk} - \widehat{\Theta}_j^T(\widehat{\Sigma}\widehat{\Theta}_k - \mathbf{e}_k)/(\widehat{\Theta}_j^T \widehat{\Sigma}_j) \tag{S.34}$$

and show that $\sqrt{n}(\widehat{\Theta}_{jk}^d - \Theta_{jk}) \rightsquigarrow N(0, \mathrm{Var}(\Theta_j^T XX^T \Theta_k))$ under certain regularity conditions. Based on the debiased estimator, we are able to test on multiple edges. For any edge set $E \subseteq V \times V$, we consider the statistic $T_E = \max_{(j,k) \in E} \sqrt{n}(\widehat{\Theta}_{jk}^d - \Theta_{jk})$. We apply the Gaussian multiplier bootstrap (Chernozhukov et al., 2013) by constructing the random variable

$$T^B = \sup_{(j,k) \in E} \frac{1}{\sqrt{n}} \sum_{i=1}^{n} \widehat{\Theta}_j^T (X_i X_i^T \widehat{\Theta}_k - \mathbf{e}_k) \xi_i.$$

and we estimate the conditional quantile by $c(\alpha, E) = \inf\{t \in \mathbb{R} \,|\, \mathbb{P}(|T_{jk}^B| > t \,|\, \{X_i\}_{i=1}^n) \leq \alpha)\}$. Applying the theory developed by Chernozhukov et al. (2013), it can be shown that $c(\alpha, E)$ is a valid quantile estimator. For the CLIME estimator, Neykov et al. (2016) prove that it is true if $\Theta \in \mathcal{U}_s$ and $(\log(dn))^6/n + s^2(\log dn)^4/n = o(1)$. The scaling condition is to show that $c(\alpha, E)$ obtained from the Gaussian multiplier bootstrap is a good quantile estimator of $\max_{e \in E} \sqrt{n}(\widehat{\Theta}_e^d - \Theta_e)$. Similar condition is also considered in Equation (23) in Chernozhukov et al. (2013) for multiple hypothesis testing. The first scaling term $(\log(dn))^6/n = o(1)$ comes from the Berry-Essen bound of the Gaussian approximation for the maximum of a sum of high dimensional random vectors. The second term $s^2(\log dn)^4/n = o(1)$ comes from approximating the statistic $\sqrt{n}(\widehat{\Theta}_e^d - \Theta_e)$ by a summation of i.i.d. random vectors. Similar condition can be found in Condition (M) of Chernozhukov et al. (2013).

For the CLIME estimator in Gaussian graphical model, Neykov et al. (2016) prove that the quantile estimator in Section 3 satisfies (3.5).

In fact, if we apply the theory of Gaussian multiplier bootstrap developed in Chernozhukov et al. (2013), (3.5) is satisfied for many other estimators for Gaussian graphical model. The proof in Neykov et al. (2016) essentially shows that (3.5) holds if we use a precision matrix estimator $\widehat{\Theta}$ in (S.34) satisfying the three conditions: (1) $\|\widehat{\Sigma} - \Sigma\|_{\max} = O_P(\sqrt{\log d/n})$; (2) $\|\widehat{\Theta} - \Theta\|_1 = O_P(s\sqrt{\log d/n})$ and (3) $\|\widehat{\Sigma}\widehat{\Theta} - \mathbf{I}_d\|_{\max} = O_P(\sqrt{\log d/n})$. These conditions (1)-(3) can also be established under the irrepresentability condition for the graphical Lasso estimator (Yuan and Lin, 2007; Ravikumar et al., 2011) and neighborhood selection (Meinshausen and Bühlmann, 2006).

For non-Gaussian graphical models, (3.5) can also established. For example, for the Gaussian copula graphical model, Gu et al. (2015) propose a pseudo maximal log-likelihood estimator based on Kendall's tau statistics and show (3.5). For the semiparametric exponential family graphical model, Yang et al. (2014b) also show (3.5) based on a semi-paramtric estimator. Applying the high dimensional debiasing methods (Zhang and Zhang, 2013; van de Geer et al., 2014; Javanmard and Montanari, 2014; Ning and Liu, 2014) and Gaussian multiplier bootstrap (Chernozhukov et al., 2013), (3.5) can also potentially be shown for many other graphical models, thus Theorem 5.1 shows the proposed skip-down method can also be directly applied.